\documentclass[10pt,a4paper,reqno]{amsart}

\usepackage{xcolor}
\usepackage{tikz}
\usepackage{tikz-cd}
\usepackage{standalone}
\usetikzlibrary{decorations.markings}
\usepackage[pagebackref]{hyperref}
\usepackage{comment}
\usepackage{cleveref}
\numberwithin{equation}{section}

\theoremstyle{plain}
\newtheorem{thm}{Theorem}[section]

\newtheorem{lem}[thm]{Lemma}
\newtheorem{cor}[thm]{Corollary}
\newtheorem{prop}[thm]{Proposition}

\theoremstyle{definition}
\newtheorem{defi}[thm]{Definition}
\newtheorem{rem}[thm]{Remark}

\newtheorem{ex}[thm]{Example}

\newcommand{\IC}{\mathbb{C}}
\newcommand{\IZ}{\mathbb{Z}}
\newcommand{\IR}{\mathbb{R}}
\newcommand{\ff}{\mathfrak f}
\newcommand{\gggg}{\mathfrak g}
\newcommand{\qqet}{\qquad\mbox{and}\qquad}

\DeclareMathOperator{\Ch}{Ch}

\title{Lagrangian Cobordisms in Liouville manifolds}
\author{Valentin Bosshard}

\begin{document}

\begin{abstract}
Floer theory for Lagrangian cobordisms was developed by Biran and Cornea in a series of papers \cite{BC,BC2,BC3} to study the triangulated structure of the derived Fukaya category of monotone symplectic manifolds. This paper explains how to use the language of stops to study Lagrangian cobordisms in Liouville manifolds and the associated exact triangles in the derived wrapped Fukaya category. Furthermore, we compute the cobordism groups of non-compact Riemann surfaces of finite type.
\end{abstract}
\maketitle

\section{Introduction}
The goal of this paper is to give a detailed proof that Lagrangian cobordisms in Liouville manifolds induce cone decompositions in the wrapped Fukaya category. As an application, we compute the cobordism groups of non-compact Riemann surfaces.
\vspace{2mm}

Fukaya \cite{FU} discovered that Lagrangian submanifolds in a symplectic manifold $W$ can be organized as the objects of an $A_\infty$-category, the \emph{Fukaya category}, whose morphisms are intersections between Lagrangian submanifolds. Higher compositions are given by counting pseudo-holomorphic disks with boundary on the Lagrangian submanifolds.

Any $A_\infty$-category can be algebraically completed to its \emph{derived category}. This category is \emph{triangulated}, in particular, all morphisms admit mapping cones. One can try then to decompose objects as iterated mapping cones of other objects and find a minimal generating set of objects. 

However, the definition of a cone is purely algebraic. In the 80's Arnold \cite{ARN1,ARN2} introduced and studied \emph{Lagrangian cobordisms} between Lagrangian submanifolds. Biran and Cornea \cite{BC,BC2} showed that Lagrangian cobordisms in fact give rise to cone decompositions in the Fukaya category. A Lagrangian cobordism is a Lagrangian submanifold living inside the product symplectic manifold $\mathbb C\times W$ which, over the complement of a compact set of $\IC$, only consists of a union of product Lagrangian submanifolds $\bigcup_{j=0}^n(\gamma_j\times L_j)\subset \mathbb C\times W$, where $\gamma_j\subset \IC$ are horizontal rays in $\IC$.

To define the Fukaya category one has to restrict to a set of admissible Lagrangian submanifolds. Biran and Cornea considered closed \emph{monotone} Lagrangian submanifolds in monotone symplectic manifolds. Lagrangian cobordisms are in general flexible \cite{ELI,AUD}, but become rigid when restricted to monotone Lagrangian cobordisms \cite{CHE,BC2}.

To study a Liouville manifold $W$ it is useful to enlarge the Fukaya category by also including exact non-compact Lagrangian submanifolds that are conical at infinity. This $A_\infty$-category is called the \emph{wrapped Fukaya category} $\mathcal W(W)$. Compared to non-wrapped Floer theory, morphisms between transverse Lagrangian submanifolds are not only intersection points but also Reeb orbits at infinity.

Biran and Cornea's cone decomposition result \cite{BC2} has the following analog in a Liouville manifold $(W, \lambda)$:

Put a Liouville structure on $\IC$ which induces the standard orientation. Let $V\subset \IC\times W$ be an exact Lagrangian submanifold which is conical outside a compact set $K\subset\IC\times W$. By making $K$ larger we can assume that $K$ is the product $K=\IC^{int}\times W^{int}$ of Liouville domains $\IC^{int}\subset \IC$ and $W^{int}\subset W$. 
Let $\pi_\IC:\IC\times W\to \IC$ be the projection. For any subset $S\subset \IC$ we denote $V\vert_S=(V\cap \pi_\IC^{-1}(S))\cap (\IC \times W^{int}).$ Note that we take the preimage and additionally intersect with $\IC\times W^{int}$. Because $V$ is conical we do not lose any information by the latter intersection. 

The Lagrangian submanifold $V$ is called a \emph{Lagrangian cobordism} if $V\vert_{\IC\setminus \IC^{int}}$ consists of finitely many \emph{conicalizations} (see Remark \ref{conicalprod} for details) of product Lagrangian submanifolds ${\gamma_0\times L_0},$  $\ldots, \gamma_n\times L_n$, where the rays $\gamma_j$ are invariant under the Liouville flow in $\IC\setminus \IC^{int} $ and the $L_j$ are exact conical Lagrangian submanifolds in $W$, called the \emph{ends} of $V$. 

The first main theorem of this paper is the following result:
\begin{thm}\label{thm1}
Let $V\subset \IC\times W$ be a Lagrangian cobordism with ends $L_0, \ldots, L_n$ (ordered counterclockwise). Then there is an iterated cone decomposition
$$L_{0}\cong[L_{n}\to L_{n-1}\to\cdots\to L_2\to L_1]$$
in the derived wrapped Fukaya category $D\mathcal W(W).$
\end{thm}
We use the notation $[M\to N]$ to describe the $A_\infty$-mapping cone of a homomorphism $M\to N$ between two $A_\infty$-modules $M$ and $N$. Theorem \ref{thm1} is stated in the ungraded setting, where taking mapping cones is associative up to quasi-isomorphism (see \cite{BC3} in section 3.1.1). Adding gradings requires to specify brackets in the iterated mapping cone. A graded cobordism $V$ induces preferred gradings on the ends $L_j$, see section \ref{sec:gradingCob}, and leads to a graded version of Theorem \ref{thm1} which is formulated as Theorem \ref{thm11}.

\vspace{4mm}
One may ask for the minimal set of objects needed to write all objects as iterated cone decomposition. The strongest known generation results are that
\begin{itemize}
    \item thimbles generate the Fukaya category of Lefschetz fibrations \cite{SEI},
    \item Lagrangian cocores generate the wrapped Fukaya category of Weinstein manifolds  \cite{CHA2,GPS2}.
\end{itemize}

Some information on the triangulated structure of an $A_\infty$-category $\mathcal A$ is captured by the \emph{Grothendieck group} $K_0(D(\mathcal A))$ of its derived category $D(\mathcal A)$. This group is the quotient of the free abelian group generated by objects of $D(\mathcal A)$ modulo the relation $A_1-A_2+A=0$ whenever $A$ is the mapping cone of a morphism $A_1\to A_2$ in $D(\mathcal A)$. 

On the other hand Biran and Cornea introduced the  group $\Omega(W)$, called the \emph{Lagrangian cobordism group} of a Liouville manifold $W$. It is the abelian group freely generated by admissible Lagrangian submanifolds in $W$ modulo the relation $L_0=L_1+\cdots+ L_n$ whenever there is a Lagrangian cobordism with ends $L_0,\ldots, L_n$. According to Theorem \ref{thm1}, Lagrangian cobordisms produce cone decompositions in the Fukaya category. So all relations between Lagrangian submanifolds defining the cobordism group also exist in the relations that define the Grothendieck group. In other words, the cobordism group $\Omega(W)$ surjects onto the Grothendieck group $K_0(D\mathcal W(W))$ of the derived wrapped Fukaya category of $W$.

One may wonder under what assumptions the cobordism group is actually isomorphic to its Grothendieck group. The only known cases are closed surfaces: Haug \cite{HAUG} proved that it is indeed an isomorphism in case of the torus, when working with the set of admissible Lagrangian submanifolds being non-contractible embedded curves. Perrier \cite{PER} studied higher genus surfaces and proved partial results for the case of unobstructed immersed curves.

In surfaces with boundary exact conical Lagrangian submanifolds are either embedded arcs connecting boundary points or exact closed curves in the interior. The second result of this text is:
\begin{thm}\label{thm2}
Let $\Sigma$ be a Riemann surface with boundary $\partial \Sigma$, obtained from a closed Riemann surface by removing finitely many open discs. Then there are isomorphisms 
$$\Omega(\Sigma)\cong K_0(D\mathcal W(\Sigma))\cong H_1(\Sigma,\partial \Sigma).$$
\end{thm}
A priori this statement may depend on the Liouville structure that we chose on $\Sigma$. However, on Riemann surfaces with boundary any two Liouville structures which induce the same orientation are equivalent as Liouville manifolds after completing the Liouville flow. Hence the result does not depend on the choice of a Liouville form on $\Sigma$.\\[2mm]
Finally, note that Theorem \ref{thm1} and Theorem \ref{thm2} will be proved in the generalized setting with stops allowed on the Liouville boundary, see Theorem \ref{thm11} and Theorem \ref{thm22}, respectively.

\subsection{Cobordism groups in higher dimensions}
In some cases the full classification of exact conical Lagrangian submanifolds is known. Most examples that we state here come from where the Nearby Lagrangian conjecture (for the zero section) or the analogous Nearby Lagrangian conjecture for fibers in cotangent bundles are known.

Denote $\Omega^c(W)$ the Lagrangian cobordism group of $W$ with admissible Lagrangian submanifolds being closed and exact. The following groups are known:
\begin{itemize}
    \item Gromov \cite{GRO} proved that there are no compact exact Lagrangian submanifolds in $\IR^{2n}$ for $n>0$. Hence: $$\Omega^c(\IR^{2n}) = 0.$$
    \item More generally, for subcritical Weinstein manifolds $W$:
    $$\Omega^c(W)= 0$$
    by general position and by applying Theorem \ref{thm11} to the null-cobordisms coming from Proposition \ref{prop:null}.
    \item The Nearby Lagragian conjecture is known for the cotangent bundle of the sphere $S^2$ (see Hind \cite{HIND}, and \cite{LS,VIT2,RIT} for earlier related results) and the cotangent bundle of the torus $\mathbb T^2$ (\cite{RIZ5}). That is, all closed exact Lagrangian submanifolds are Hamiltonianly isotopic to the zero section. Hence: 
    $$\Omega^c(T^*S^2)\cong \IZ \qqet \Omega^c(T^*\mathbb T^2)\cong \IZ$$ are generated by the zero sections, respectively (they cannot be torsion as the zero section in the Grothendieck group is not torsion and the cobordism group maps to the Grothendieck group of the Fukaya category.)
\end{itemize}

Let $\Lambda$ be a Legendrian submanifold in the Liouville boundary $\partial W$ of a Liouville manifold $W$. Denote $\Omega(W, \Lambda)$ the Lagrangian cobordism group of $W$ with admissible Lagrangian submanifolds being exact and having boundary $\partial L=\Lambda\subset \partial W$ at infinity and relations given by Lagrangian cobordisms having admissible Lagrangian ends.

For the cotangent bundle $T^*\Sigma$ of an open surface $\Sigma$ the Nearby Lagrangian conjecture is known, i.e. every conical exact Lagrangian which is equal to a fiber outside a compact set is Hamiltonianly isotopic to the fiber by a compactly supported Hamiltonian (C\^ot\'e and Rizell \cite{RIZ}; see also \cite{EP} for earlier related work). Let $\Lambda$ be the Legendrian boundary of a fiber in $\partial W$. Hence
$$\Omega(T^*\Sigma, \Lambda)\cong \IZ$$ is generated by the fiber. In subsection \ref{sec:highercob} we discuss some more techniques to find interesting relations in cobordism groups.

\subsection{Strategy of the proofs}
A refinement of the wrapped Fukaya category of a Liouville manifold is the \emph{partially} wrapped Fukaya category. A \emph{stop} is a closed subset of the Liouville boundary $\partial W$ of a Liouville manifold $W$. The name stop already suggests that in the space of morphisms we do not take into account Reeb orbits that intersect the stop. Put differently, we do not wrap Lagrangian submanifolds over stops. This extra structure turns out to be very useful in the context of Lagrangian cobordisms. Given a Lagrangian cobordism $V\subset \IC\times W$ with $n+1$ ends we can put a stop $\ff_n$ on the Liouville boundary $\partial \IC$ which consists of $n+1$ points, one point between each end. 
\vspace{2mm}

In \cite{GPS2} it is shown that for $\IC$ equipped with the stop $\ff_n$ and any stopped Liouville manifolds $(W, \gggg)$ there is a cohomologically full and faithful Künneth $A_\infty$-bifunctor
$$\mathcal W_{\ff_n}(\IC)\times \mathcal W_\gggg(W)\to \mathcal W_\mathfrak h(\IC\times W)$$ where $\mathfrak h$ is the product stop $$\mathfrak h=(\ff_n\times \mathrm {sk}(W))\cup (\ff_n\times \gggg)\cup(\mathrm {sk}(\IC)\times \gggg).$$ The set $\mathrm{sk}(W)$ denotes the skeleton of $W$.  The functor sends a pair $(\alpha,L)$ to the conicalization $\alpha\tilde\times L$ of the product Lagrangian $\alpha\times L$ in the product Liouville manifold $\IC\times W$. We refer to section \ref{ssec:Kunneth} for the precise formulation. 

A Lagrangian cobordism $V$ defines an object in $\mathcal W_\mathfrak h(\IC\times W)$. Let $\mathcal Y(V)$ be the associated Yoneda module. Each arc $\alpha$ in $\mathcal W_{\ff_n}(\IC)$ induces a functor $\mathcal W_\gggg(W)\to \mathcal W_\mathfrak h(\IC\times W)$. Pulling back the Yoneda module $\mathcal Y(V)$ along this functor, yields a $\mathcal W_\gggg(W)$-module $M_{\alpha,V}$ that depends only, up to a quasi-isomorphism, on where (between which stops) the conical ends of the arc $\alpha$ lie. 

This organizes as a functor $$\mathrm{ev}_V:\mathcal W_{f_n}(\IC)^{opp}\to mod_{W_\gggg(W)}$$ sending $\alpha$ to $M_{\alpha, V}$. The partially wrapped Fukaya category $\mathcal W_{f_n}(\IC)^{opp}$ serves as the universal model of cone decompositions. Exact triangles of arcs in $\IC$ give rise to exact triangles of $\mathcal W_\gggg(W)$-modules under the $A_\infty$-functor $\mathrm{ev}_V$, as $A_\infty$-functors respect exact triangles. Picking and isotoping the right arcs in $\IC$ as in \cite{BC2} concludes the proof of Theorem \ref{thm1}.
\vspace{4mm}

For Theorem \ref{thm2} we follow \cite{HKK}. The argument is based on the contractibility of the arc complex for surfaces with boundary \cite{HAR1,HAR2}, it is completely topological. 
A \emph{full arc system} $A$ on a Riemann surface $\Sigma$ with boundary is a set of non-isotopic arcs such that when cutting the surface along all arcs in $A$ we are left with topological disks. In \cite{HKK} it is shown that the Grothendieck group of the derived wrapped Fukaya category of a Riemann surface is generated by the elements of an arc system. In fact, a full arc system generates any exact closed curve and any arc in $\Sigma$ as an iteration of the following operations
\begin{itemize}
    \item Two arcs $\gamma_1$ and $\gamma_2$ are identified if they are isotopic relative $\partial \Sigma,$
    \item An arc $\gamma$ is $0$ if it can be isotoped into the boundary of $\partial \Sigma$,
    \item $\gamma=\gamma_1+\gamma_2$ if $\gamma_1$ and $\gamma_2$ are curves intersecting exactly once and transversely with intersection index $1$, and the resulting surgery is $\gamma$.
\end{itemize}
All these operations can be realized by cobordisms. Thus the Grothendieck group and the cobordism group coincide. Note that it is not necessary to include Hamiltonian isotopies of closed curves, only isotopies of arcs (which are always Hamiltonian). This is possible because the surgery of a closed curve that intersects an arc only once is an arc, i.e. we can present any closed curves as the sum of two arcs in the cobordism group. 

\subsection{Relation to other work}
Independently of Biran and Cornea's work, Nadler and Tanaka proposed a theory of Lagrangian cobordisms in Liouville manifolds in \cite{TAN2}. Tanaka developed this theory further in various follow-up papers \cite{TAN4,TAN5,TAN6,TAN,TAN7}. 

Biran and Cornea's goal was from the beginning motivated by Floer theory culminating in the functorial way that cobordisms induce exact triangles. On the other hand, Nadler and Tanaka a priori view Lagrangian cobordisms as morphisms in a stable $\infty$-category, where mapping cones exist intrinsically. 

This work is Floer-theoretically oriented. That is, we try to redo Biran and Cornea's work in the setting of Liouville manifolds with (possibly non-compact) exact conical Lagrangian and its consequences in wrapped Floer theory. Moreover, the definition of cobordism given here suffices to compute the cobordism group of all stopped Riemann surfaces.

Tanaka's arguments relating Lagrangian cobordisms to Floer theory \cite{TAN5,TAN6} are given without using the language of stops. We incorporate stops in our arguments following \cite{SS,GPS1,GPS2} to advertise this more recent notion. We also show how and when products of Lagrangian submanifolds should be conicalized in Lagrangian cobordism theory.

In our definition, a Lagrangian cobordism consists of conicalized products of Lagrangian submanifolds over the complement of a large compact set of $\IC$. On $\IC$ we put any Liouville structure. Nadler and Tanaka work with $\IC$ as the cotangent bundle $T^*\IR$. Their cobordism have a time direction. Ours have a rotational symmetry (also discussed in \cite{TAN3}). The topology of their Lagrangian cobordism are only not allowed to go arbitrarily far into the positive cotangent direction (they call this  \emph{skeleton-avoiding} or \emph{non-characteristic}). Cobordisms come there with a direction from negative real part to positive real part. So they are empty in the far lower half plane part, consist of two Lagrangian ends on the real line and are arbitrarily wild in the upper half plane of $\IC$.
These cobordism induce not directly cone decompositions, they are regarded as morphisms. Furthermore, their Lagrangian cobordism $\infty$-category has higher morphisms living in $T^*\IR^n\times W$. The papers \cite{TAN5,TAN6} prove that these morphisms induce Floer theoretic morphisms.

To unify these two approaches one may want to jump to section \ref{sec:big} for references on how stabilization of Liouville manifolds, cone decompositions and Lagrangian cobordisms fit into a larger picture.

\subsection{Structure of the text}
In section \ref{sec2} we briefly revise some basic constructions in Liouville manifolds that prepare the definition of Lagrangian cobordisms. In section \ref{sec3} we define Lagrangian cobordisms in the framework of Liouville manifolds and give some examples. Section \ref{sec4} recalls the construction of partially wrapped Fukaya categories, in particular of a disk with stops on the boundary and proves the cone decomposition stated in Theorem \ref{thm1}. The Lagrangian cobordism groups of non-compact Riemann surfaces are discussed in section \ref{sec5}. 

\subsection*{Acknowledgements} I would like to thank my PhD advisor, Paul Biran, for his inspiration and his guidance by always asking the right questions. Further, special thanks to Baptiste Chantraine for our conversation about Lagrangian cobordisms and for pointing out Example \ref{ex:noncob}. I would like to thank the referee for carefully reading this work and for useful comments that helped to improve the paper. The author was partially supported by the Swiss National Science Foundation (grant number  $200021\_204107$).

\section{Some Constructions in Liouville Manifolds}\label{sec2}
A \emph{Liouville domain} $W^{int}$ is a compact exact symplectic manifold $(W^{int}, \omega)$ with contact boundary $\partial W$. More precisely, we fix a 1-form $\lambda$, called \emph{Liouville form}, such that $\omega=d\lambda$. The contact boundary assumption means that the vector field $Z$, called the \emph{Liouville vector field}, uniquely defined by $\iota_Z\omega=\lambda$, points strictly outward the boundary $\partial W$ of $W^{int}$. Denote $\alpha=\lambda\vert_{\partial W}$ the contact form. A \emph{stop} $\ff$ is a closed subset of the boundary $\partial W$. As the name suggests, all the dynamics that we are going to consider avoid the stop.

By flowing inwards from the boundary $\partial W$ along the Liouville vector field $Z$ with negative time $t$ one obtains the negative part of the symplectization $(0,1]\times \partial W\subset W^{int}$ of $(\partial W, \alpha)$. The first coordinate is called the radial coordinate $r=e^t\in (0,1]$. To make the vector field $Z$ complete we glue the positive part $W
^{cone}=[1,\infty)\times \partial W$ of the symplectization (called \emph{conical part}) to $W^{int}$ and denote this (non-compact) manifold by $$W=W^{int}\cup_{\partial W} W^{cone}.$$ On the symplectization part $(0,\infty)\times \partial W$ the Liouville form is given by $\lambda=r\alpha$.
All information on the boundary $\partial W$ should extend to the cone $W^{cone}$ by using the Liouville flow. In presence of a stop $\ff\subset \partial W$ we denote its extension $[1,\infty)\times\ff$ to the cone $W^{cone}$  also simply by $\ff$.

Manifolds that are completions of Liouville domains are called \emph{Liouville manifolds}. It is convenient to be flexible on where the conical part $W^{cone}$ starts. Our main object of interest are stopped Liouville manifolds $(W,\lambda,\ff)$ without fixing the Liouville domain $W^{int}$ and the start of the stop $\ff$. If a property is true on the cone $W^{cone}$ for some choice of Liouville domain $W^{int}\subset W$ then we say the property holds \emph{at infinity}. Note also that the manifold $W$ has no boundary. The \emph{Liouville boundary} $\partial W$ is just a (non-fixed) choice of a contact hypersurface on which the conical part of $W$ is modelled.

The set $$\mathrm {sk}(W)=W\setminus \big((0,\infty)\times \partial W\big)$$ is called the \emph{skeleton} of $W$, i.e. the points that do not escape to infinity under the Liouville flow. The skeleton is compact, possibly singular, but independent of the choice of the Liouville boundary $\partial W$.

\subsection{Linear Hamiltonians}\label{sec:linear}
As the conical part of a Liouville manifold  is the symplectization of the contact manifold, we restrict to Hamiltonian functions which induce contact dynamics at infinity.

Let $(W, \lambda)$ be a Liouville manifold and $(\partial W, \alpha)$ a choice of a contact manifold that models the conical part of $W$. Denote $\xi=\ker \alpha$ the contact structure on $\partial W$.
Let us recall the relation of contact dynamics on the Liouville boundary $\partial W$ to dynamics on the Liouville manifold $W$ at infinity.
The following spaces are in bijection:
\begin{itemize}
    \item Functions $h_t$ on $\partial W$,
    \item Vector fields $v_t$ on $\partial W$ that are lifts of sections of $T(\partial W)/\xi$,
    \item An isotopy $\psi_t$ of $\partial W$ starting at $\psi_0=id_{\partial W}$ such that $\psi_t$ preserves the contact structure $\xi$ for all times $t$.
\end{itemize}
They are called \emph{contact Hamiltonians $h_t$, contact vector fields $v_t$ and contact isotopies $\psi_t$}, respectively. 

The relations among them are given by: 
$$\alpha(v_t)=h_t \qqet \iota_{v_t}d\alpha=dh_t(R)\alpha-dh_t$$ where $R$ is the \emph{Reeb vector field} of $(\partial W, \alpha)$, i.e. the vector field uniquely defined by $\alpha(R)=1$ and $\iota_Rd\alpha=0$. The isotopy $\psi_t$ is the flow of the vector field $v_t$.
As all the $\psi_t$ preserve $\xi=\ker\alpha$, we must have that $\psi_t^*\alpha=f_t\alpha$ for some positive function $f_t:\partial W\to \IR.$ We can recover $f_t$ by $f_0=1$ and $\partial_tf_t=dh_t(R).$ 
\vspace{2mm}

In the situation of a Liouville manifold $(W,\lambda)$ with Liouville flow $Z$ the following spaces are in bijection:
\begin{itemize}
    \item Functions $H_t:W\to \IR$ such that $Z(H_t) = H_t$ at infinity,
    \item Symplectic vector fields $X_t$ on $W$ that commute with $Z$ at infinity,
    \item An isotopy $\varphi_t$ of $W$ starting at $\varphi_0=id_{W}$ such that the $1$-form $\varphi_t^*\lambda-\lambda$ is exact everywhere and vanishes at infinity.
\end{itemize}
They are called \emph{linear Hamiltonians $H_t$, Hamiltonian vector fields of contact type $X_t$ and Hamiltonian isotopies of contact type $\varphi_t$}, respectively. We use the convention that the Hamiltonian vector field $X_t$ associated to a function $H_t:W\to \IR$ is given by $\omega(X_t,\cdot)=-dH_t$ and denote the induced flow on $W$ by $\varphi_t$.

\vspace{4mm}
Any contact Hamiltonian $h_t$ on $\partial W$ induces a linear Hamiltonian $H_t$ on $W$ by defining $H_t=rh_t$ at infinity and extending it to $W^{int}$ by a cutoff function. Actually, the definition of linear Hamiltonian is equivalent to Hamiltonians that have the form $H_t=rh_t$ at infinity by choosing a cone $W^{cone}$ with conical coordinate $r$. Its Hamiltonian vector field $X_t$ at infinity is $X_t(r,z)=(\partial_t f_t,v_t(z))$ where $v_t$ is the contact vector field associated to $h_t$ and the induced Hamiltonian diffeomorphism is of the form $\varphi(r, z)=(rf_t(z),\psi_t(z))$. 

\begin{ex}
Linear Hamiltonians of the form $H_t=r$ at infinity (i.e. $h=1$ in the previous notation), induce the Hamiltonian vector field $X$ which is the Reeb vector field $R$ in each of the contact manifolds $({r}\times\partial W, r\alpha)$ for fixed $r$.
\end{ex} 

In case of a stopped Liouville manifold $(W,\lambda, \ff)$ we demand $H_t$ to vanish on the set $\ff$ at infinity. The corresponding flow $\varphi_t$ is the identity on $\ff$. So if a subset $L\subset W$ does not intersect $\ff$ at infinity then $\varphi_t(L)$ stays disjoint from $\ff$ at infinity for all time $t$.  

\subsection{Conical Lagrangian submanifolds and conicalizations}\label{sec:conical}
Similarly to the choice of linear Hamiltonians, we restrict to Lagrangian submanifolds which are modeled at infinity by Legendrian submanifolds for some choice of a contact Liouville boundary.

In this paper we only consider \emph{exact} Lagrangian submanifolds $L$ in Liouville manifolds $(W, \lambda)$. By exact we mean that there is a function  $f:L\to \IR$, such that $$\lambda\vert_L=df.$$ Such a function $f$ is called a \emph{primitive} of $L$.

Given a choice of a Liouville domain $W^{int}\subset W$ define 
\vspace{-2mm}

$$\partial L=L\cap \partial W, \qquad L^{int}=L\cap W^{int}\qquad \text{and} \qquad L^{cone}=L\cap W^{cone}$$ for any Lagrangian submanifold $L\subset W$.

In the contact manifold $(\partial W,\alpha=\lambda\vert_{\partial W})$ we are interested in submanifolds $\Lambda$ of $\partial W$ which are \emph{Legendrian}, i.e. $T\Lambda\subset T(\partial W)$ lies in the contact distribution $\xi=\ker\alpha$.

The corresponding asymptotic notion for Lagrangian submanifolds $L$ in $W$ is characterized by the following equivalent conditions:
\begin{itemize}
    \item $\lambda\vert_L$ vanishes at infinity,
    \item $L$ is invariant under the Liouville flow at infinity,
    \item any primitive of $L$ is locally constant at infinity.
\end{itemize}
These Lagrangian submanifolds are called \emph{conical}\footnote{To be consistent we should use the term \emph{conical at infinity}. However, for notational simplicity we write \emph{conical} and use the term \emph{strongly conical} \label{strexact} or \emph{strongly exact} if $L$ is invariant under the Liouville flow everywhere.}. One can add to primitives a real number to make it vanish on one specific connected component at infinity. Using a Hamiltonian deformation we can arrange that the Lagrangian submanifold admits a primitive that vanishes at infinity:
\begin{lem}[\cite{GPS2} {\normalfont Lemma 6.1}]\label{lem:GPS}
Any conical exact Lagrangian submanifold can be deformed by a compactly supported Hamiltonian diffeomorphism to a conical exact Lagrangian whose primitive vanishes at infinity, i.e. \!\!admits a compactly supported primitive. Furthermore, if the original Lagrangian does not intersect the skeleton, so does the deformed Lagrangian submanifold.
\end{lem}  
This is done by a compactly supported Hamiltonian which equals $1$ over a large compact set that contains the skeleton.\\[2mm]
We are interested in exact Lagrangian submanifolds $L$ that intersect the Liouville boundary $\partial W$ transversely and such that $\partial L\subset \partial W$ is a Legendrian submanifold. In fact, if $L$ is conical and the Liouville domain large enough then $L^{cone}$ is equal to $[1, \infty)\times\partial L\subset W^{cone}$, the positive conicalization of the Legendrian $\partial L$.

We can deform Lagrangian submanifolds of the former type into conical Lagrangians by making them parallel to the Liouville flow near $\partial W$ and then take their conical extension. We refer to this process as \emph{conicalization}.
\begin{lem}[\cite{AS} {\normalfont Lemma 3.1}]\label{makeconical}
Let $L_0\subset W^{int}$ be an exact Lagrangian submanifold that intersects $\partial W$ transversely such that $\partial L_0$ is Legendrian. Then $L_0$ is Hamiltonianly isotopic (rel $\partial W$) to a Lagrangian submanifold $L_1\subset W^{int}$ which is invariant under the Liouville flow near $\partial W$. Moreover, the Hamiltonian isotopy can be made to have support arbitrarily close to $\partial W$. 
\end{lem}
The lemma ensures that we can extend $L_1\subset W^{int}$ by the Liouville flow conically to a conical Lagrangian submanifold in $W$.
\begin{proof}
Let $f$ be a primitive of $L_0$. As $df$ vanishes on $\partial L_0$, i.e. $f$ is locally constant on $\partial L_0$, there is a global function $h:W^{int}\to\IR$ supported near $\partial L_0$ that vanishes on $\partial W$ and $d(f-h)$ is equal to zero in a neighbourhood of $\partial L_0$ in $W^{int}$. Note that the Liouville vector field for the Liouville structure $\lambda-dh$ is parallel to $L$ near $\partial L_0$. By a standard argument in Liouville manifolds (\cite{CE}, Proposition 11.8) there is a Hamiltonian diffeomorphism $\Phi_t$ of $W^{int}$ satisfying $\Phi_t^*\lambda=\lambda-tdh$. Then $L_1=\Phi_1(L_0)$ is our desired Lagrangian deformation which is constant on $\partial L_0$ as $h$ vanishes on $\partial W$.
\end{proof}
\begin{ex}
An example for a conicalization of an arc living in the disk $\IC^{int}$ is  given in Figure \ref{fig:conicali}. 
\begin{figure}
    \centering
    \begin{tikzpicture}[scale=0.60]

\draw[line width=0.25mm, dashed ] (-2,0) circle (3);
\node at (-2, -2.5) {$\partial \mathbb{C}$};
\draw (-4.83,1) --node[pos=0.5, below=1mm] {$\gamma_0$}  (0.83,1);

\draw[line width=0.25mm, dashed ] (8,0) circle (3);
\node at (8, -2.5) {$\partial \mathbb{C}$};
\draw (3.7,1.5) -- (5.17,1);
\draw (5.17,1) ..  controls (5.6,0.85) and (5.6,1) .. (5.9,1)--(5.9,1) .. controls (6.1,1) and (6.1,1) .. (6.2,1)--node[pos=0.5, below=1mm] {$\gamma_1$} (9.8,1)..  controls (9.9,1) and (9.9,1) .. (10.1,1)--(10.1,1) .. controls (10.4,1) and (10.4,0.85) .. (10.83,1);
\draw (10.83,1) -- (12.3,1.5);
\end{tikzpicture}
    \caption{Drawn on the left is an arc $\gamma_0$ in $\IC^{int}$ which intersects $\partial  \IC$ transversely. Its conicalization $\gamma_1$ is given by flowing the Legendrian $\partial \gamma_0\subset \partial \IC$ out to the cone $\IC^{cone}$ and a small bit to $\IC^{int}$. As we want to deform the curve only at infinity we defined in the construction of the conicalization a function $h$ which is a cutoff (and a shift) of a chosen primitive of $\gamma_0$. Hence we can achieve $\gamma_0=\gamma_1$ on any previously chosen proper compact subset of $\IC^{int}$.}
    \label{fig:conicali}
\end{figure}
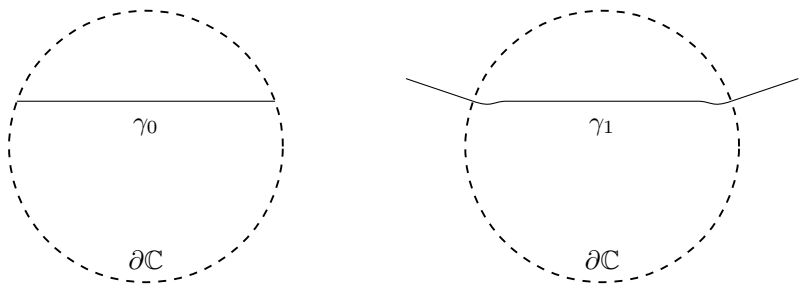
\end{ex}

\begin{rem}
A conicalization is unique in the following sense: If $L_0$ is an exact Lagrangian submanifold and $W^{int}_t$ a smooth family of Liouville domains of $W$ such that $L_0$ and $W^{int}_t$ intersect transversely in Legendrians for all $t$, then the conicalizations for each $t$ constructed in Lemma \ref{makeconical} are unique up to linear Hamiltonian deformations. In case of non-transverse intersections even the topology can change.
\end{rem}

A very similar argument applies when we want to conicalize the product of two conical Lagrangian submanifolds.

Given two Liouville manifolds $(W_1, \lambda_1)$ and $(W_2,\lambda_2)$, the product $L\times K\subset W_1\times W_2$ of two exact conical Lagrangian submanifolds $L\subset W_1$ and $K\subset W_2$ does not need to be conical with respect to the product Liouville structure $\lambda_1+\lambda_2$ on $W=W_1\times W_2$. The problem is that $\lambda_1+\lambda_2$ is not automatically zero when restricted to a neighbourhood of $\partial (L\times K)=\partial L\times K^{int}\cup L^{int}\times \partial K$ in $L\times K$. It is zero if and only if $\lambda_1\vert_L$ and $\lambda_2\vert_K$ both vanish everywhere, not just at infinity, i.e. if they both are strongly exact.

\begin{lem}[\cite{GPS2} {\normalfont Section 6.2}]\label{conicalprod}
Suppose that both $L\subset W_1$ and $K\subset W_2$ are exact Lagrangian submanifolds that admit primitives that vanish at infinity. Then the product Lagrangian submanifold $L\times K$ can be deformed to a conical Lagrangian submanifold $L\tilde\times K\subset W_1\times W_2$. Choose a large Liouville domain $U\subset W_1\times W_2$ and a slightly smaller Liouville domain $U'$ contained in $U$. Then the product $L\times K$ and the conicalization $L\tilde\times K$ agree on $U'$ and are Hamiltonian isotopic in $U\setminus U'$. On $W\setminus U$, the conicalization $L\tilde \times K$ is conical with respect to $\lambda_1+\lambda_2$. This process is also called \emph{conicalization}. Moreover, different choices of Liouville domains $U$ and $U'$ lead to conicalizations which are isotopic by linear Hamiltonian diffeomorphisms. 
\end{lem}

\begin{proof}
Let $f_L$ and $f_K$ be compactly supported primitives of $L$ and $K$, respectively. Let $W_1^{int}\subset W_1$ and $W_2^{int}\subset W_2$ be large Liouville domains. Then choose smooth functions $\psi_L: W_1 \to [0, 1]$ and $\psi_K:W_2 \to [0, 1]$  which vanish in large subsets $U_1\subset W_1^{int}$ and $U_2\subset W_2^{int}$, equal $1$ near $\partial L$ and $\partial K$, respectively, and are invariant under the Liouville flow at infinity. The Liouville domains in the statements are (a smoothening of the cornered manifolds) $U=W^{int}_1\times W^{int}_2$ and $U'=U_1\times U_2$. 

The product  $L\times K$ is conical with respect to 
$$\widetilde \lambda =\lambda_1+\lambda_2-d(f_L\psi_K)-d(f_K\psi_L).$$
One can show that choosing the Liouville domain of $(W_1, \lambda_1)$ and $(W_2, \lambda_2)$ large enough ensures that $W_1^{int}\times W_2^{int}$ is a Liouville domain for both $\lambda_1+\lambda_2$ and $\widetilde \lambda.$ As in Lemma \ref{makeconical}, there is a uniquely defined Hamiltonian isotopy $\Phi_t$ of $W_1\times W_2$ starting at $\Phi_0=\mathrm{id}$ such that $\Phi_1^*(\lambda_1+\lambda_2)=\widetilde \lambda.$ The Lagrangian submanifold $L\tilde\times K=\Phi_1(L\times K)$ is our desired Lagrangian deformation.
\end{proof}

When discussing cobordisms, we will need to deform products of the form $\gamma\times L$, for conical curves $\gamma\subset (\IC, \lambda_\IC)$ and conical Lagrangian submanifolds $L\subset (W, \lambda_W)$, to conical Lagrangian submanifolds $\gamma\tilde \times L$ in $(\IC\times W, \lambda_\IC+\lambda_W)$. Note that the conicalization $\gamma\tilde\times L$ coincides with $\gamma\times L$ over points where the Liouville vector field on $\IC$ is zero.

\begin{ex}
Suppose $\gamma$ is a strongly exact curve, e.g. a straight line through the origin in $\IC$ equipped with the standard radial Liouville structure and suppose we chose a large enough Liouville domain $\IC^{int}\subset \IC$ as in Figure \ref{fig:conicali2}. 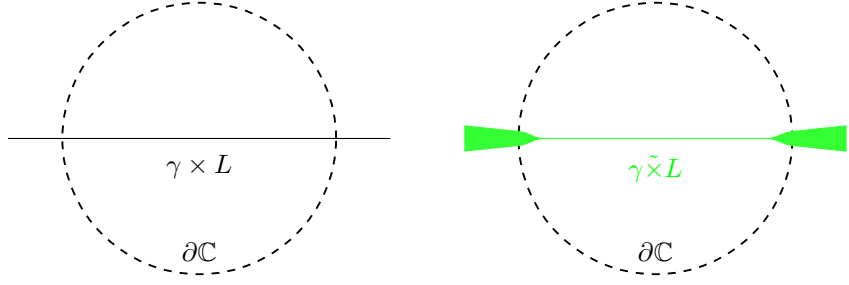
\begin{figure}
    \centering
        \begin{tikzpicture}[scale=0.60]

\draw[line width=0.25mm, dashed ] (-2,0) circle (3);
\node at (-2, -2.5) {$\partial \mathbb{C}$};
\draw (-6.19,0) --node[pos=0.5, below=1mm] {$\gamma\times L$}  (2.19,0);

\draw[line width=0.25mm, dashed ] (8,0) circle (3);
\node at (8, -2.5) {$\partial \mathbb{C}$};
\draw[fill=green!80, draw=green!80] (3.81,0.28)--(5.03,0.15) ..  controls (5.3,0.1) and (5.3,0) .. (5.6,0)--(5.6,0) .. controls (5.3,0) and (5.3,-0.1) .. (5.03,-0.15)--(3.81,-0.28);
\draw[fill=blue!80, draw=green!80] (5.6,0) --node[pos=0.5, below=1mm, text=green] {$\gamma\tilde\times L$} (10.4,0);
\draw[fill=green!80, draw=green!80] (16-3.81,0.28)--(16-5.03,0.15) ..  controls (16-5.3,0.1) and (16-5.3,0) .. (16-5.6,0)--(16-5.6,0) .. controls (16-5.3,0) and (16-5.3,-0.1) .. (16-5.03,-0.15)--(16-3.81,-0.28);\end{tikzpicture}
    \caption{On the left we see the projection of the product Lagrangian submanifold $\gamma\times L$ to $\IC$, i.e. we see the 
    the curve $\gamma$.
    On the right we see the projection of the conicalization $\gamma\tilde\times L$ to $\IC$. As we can make the introduced angle coming from the conicalization in the right figure arbitrarily small, we usually draw conicalizations still as in the left figure.}
    \label{fig:conicali2}
\end{figure}
If $\psi_\gamma:\IC\to[0,1]$ vanishes on $U_1\subset \IC^{int}$ the conicalization $\gamma\tilde\times L$ is equal to $\gamma\times L$ in $U_1\times W$. 
Near $\partial\IC\times W^{int}$ there is a small deformation area which constructs a Legendrian $\partial(\gamma\tilde\times L)\subset \partial\IC\times W^{int}\subset \partial(\IC\times W)$. On $\IC^{cone}\times W$ the conicalization $\gamma\tilde\times L$ is the positive cone of the Legendrian $\partial(\gamma\tilde\times L)$ with respect to the product Liouville flow $\lambda_\IC+\lambda_W.$

By choosing the support of $\psi_\gamma$ close enough to $\gamma_{cone}$ we can have an arbitrarily small angle of the projection of the ends. For the sake of not overloading figures, we usually draw the conicalized curves as in the left picture in Figure \ref{fig:conicali2}, even though one should introduce a small angle in the picture.
\end{ex}

\begin{rem}\label{rem:nonintersect}
Let $V\subset \IC\times W$ be a closed subset which is conical (i.e. invariant under the product Liouville flow outside some product Liouville domain $\IC^{int}\times W^{int}$). Assume that the projection of $V\cap (\IC\times W^{int})$ to $\IC$ is compact. 

Furthermore, suppose $\gamma$ is a curve in $\IC$ not intersecting this projection and $L$ is an exact conical Lagrangian submanifold in $W$. Then by choosing the deformation in the construction of the conicalization $\gamma\tilde\times L$ small enough, the conical set $V$ and the conicalization $\gamma\tilde\times L$ do not intersect.

Similarly, if $\gamma_0$ and $\gamma_1$ are transverse exact conical curves in $\IC$, and $L_0$ and $L_1$ transverse conical exact Lagrangian submanifolds in  $W$, we can ensure that the intersection $(\gamma_0\times L_0)\cap(\gamma_1\times L_1)=(\gamma_0\cap \gamma_1)\times (L_0\cap L_1)$ is in bijection with the intersection of its conicalizations $\gamma_0\tilde\times L_0$ and $\gamma_1\tilde\times L_1$. We will frequently use this procedure to intersect ends of a cobordism with a conicalization $\gamma\tilde\times L$.  
\end{rem}

\begin{rem}
We only defined the conicalization of a product of Lagrangian submanifolds where both Lagrangians admit a primitive that vanishes at infinity. According to Lemma \ref{lem:GPS} this is not too restrictive. We most of the time omit the argument why the Lagrangian submanifolds can be assumed to have this property. For example, even for curves in $\IC$, we always need to deform our curve first before using it for a conicalization. As Lemma \ref{lem:GPS} can be applied to deform the Lagrangian far enough out at infinity we do not run into issues in the form of extra intersection points.
\end{rem}

We end this subsection by emphasizing that we also include exact compact Lagrangian submanifolds $L$ in the definition of conical Lagrangians. Their conical part $L^{cone}$ is just empty. In case of a stopped Liouville manifold $(W, \lambda, \ff)$ we only allow Lagrangian submanifolds which are disjoint at infinity from the stop $\ff$. For example, if $\ff=\partial W$ we restrict our set of admissible Lagrangian submanifolds to exact compact ones.

\subsection{Liouville structures on $\IC$}\label{sec:Liouv}
We are going to define cobordisms living in $\IC\times W$ by requiring them to be exact Lagrangian submanifolds which are conical. This depends on the choice of the Liouville form on $\IC$. We allow any Liouville form on $\IC$ which has compact skeleton and whose differential induces the standard orientation on $\IC$. 
The space of all such \emph{admissible} Liouville structures is convex, thus contractible. The main example to keep in mind is the radial Liouville structure $\frac12(xdy-ydx)$. 

We postpone further discussion to section \ref{sec:Liouvcont}.

\subsection{Lagrangian suspension}\label{sec:Lsusp}
A Hamiltonian $H_t$ on an arbitrary symplectic manifold $(W, \omega)$ with $t\in\IR$ and a Lagrangian submanifold $L$ in $W$ define a \emph{Lagrangian suspension} $V_{L,H}\subset \IC\times W$. It is the image of the map $i_{L,H}:\IR\times L\to \IC\times W$ given by \[(t,w)\mapsto\big(t,H_t(\varphi_t(w)), \varphi_t(w)\big),\] where $\varphi_t$ is the Hamiltonian isotopy generated by $H_t$.
Whenever $L$ is embedded, so is $V_{L,H}$. Moreover $V_{L,H}$ is a Lagrangian submanifold with respect to the symplectic form $dx\wedge dy+\omega$ on $\IC\times W$. 

A standard calculation (see \cite{MS} Proposition 9.3.1) shows that if $(W, \lambda)$ is an exact symplectic manifold and $L$ an exact Lagrangian submanifold of $W$, then so is $V_{L,H}$ in $\IC\times W$.
Concretely, if $f:L\to\IR$ is a primitive of $L$ (i.e. $\lambda\vert_L=df$) then $$F(t, w)=f(w)+\int_0^t\varphi_s^*\big(\lambda(X_s)-H_s\big) ds$$ is a primitive of $V_{L,H}$ with respect to the Liouville form $-ydx+\lambda$ on $\IC\times W$. \footnote{The signs in this formula may differ from the ones used in other texts. Our convention is $\omega=d\lambda$ and $\omega(X_t, \cdot)=-dH_t$. Moreover, note that the Lioville form $-ydx$ on $\IC$ induces the standard orientation on $
\IC$ but is not the canonical Liouville form $ydx$ on $\IC$ when treated as the cotangent bundle $T^*\IR$.} 

The notions in the first two sections \ref{sec:linear}, \ref{sec:conical} are compatible with Lagrangian suspensions:
 
\begin{lem}\label{lemma:susp}
Suppose that $L$ is an exact conical Lagrangian submanifold of $W$ and $H_t:W\to \IR$ a linear Hamiltonian. Then $V_{L,H}$ is an exact conical Lagrangian submanifold. Moreover, if $H_t$ vanishes on a stop $\ff\subset \partial W$ then $V_{L,H}$ is disjoint from $\IC\times\ff$.
\end{lem}

\begin{proof}
We have already argued that the $V_{L,H}$ is exact. The Lagrangian suspension $V_{L,H}$ of linear Hamiltonian $H_t:W\to \IR$, i.e. $H_t=rh_t$ for $r\in [1,\infty)$, and an exact conical Lagrangian $L$, i.e.  $L^{cone}=[1,\infty)\times \partial L$, is given by $$\{(t,r h_t(\psi_t(z)), rf_t(z),\psi_t(z))\mid t\in\IR, r\in [1,\infty), z\in \partial W\}$$
at infinity using the notation from section \ref{sec:linear}. The appearance of the two $r$-factor in the two cone-coordinates $y$ of $\IC$ and $r$ of $W$ shows that the Lagrangian suspension $V_{L,H}$ is invariant under the product Liouville vector field of $\IC\times W$.
\end{proof}

\begin{rem}\label{rem:compactliou}
The Liouville form $-ydx$ on $\IC$ admits only non-compact Liouville domains. This is in contrast to section \ref{sec:Liouv} where we want to regard $\IC$ as a Liouville manifold with compact Liouville domain. Assume that a Hamiltonian $H_t:W\to \IR$ is zero outside some compact time interval $t\in [a,b]\subset \IR$. Then the Lagrangian suspension $V_{L,H}$ of a Lagrangian submanifold $L$ is equal to $(-\infty, a]\times L_a$ for $t\in(-\infty, a]$ and equal to $[b, \infty)\times L_b$ for $t\in[b, \infty)$. Here we write $L_t=\varphi_t(L)$. We can change the Liouville structure $-ydx$ outside $T^*[a,b]\subset \IC$ (see below) such that its Liouville vector field points radially outward at infinity and conicalize the product Lagrangian submanifolds.

Consider a function $\rho:\IR\to \left[\frac 12,1\right]$ which is $1$ on $[a,b]$, equal to $\frac12$ outside a larger compact interval $I$ containing $[a,b]$, and has derivative $|\rho'|<1$ everywhere. Then the $1$-form $\lambda_\IC=(1-\rho(x))xdy-\rho(x)ydx$  defines a Liouville form on $\IC$ with volume form $$d\lambda_\IC= (1-\rho'(x))dx\wedge dy.$$

The Liouville form $\lambda_\IC$ coincides on $T^*[a,b]$ with the (negative of the) cotangent form $-ydx$ and coincides with the radial form $\frac12(xdy-ydx)$ outside $T^*I$. Its skeleton is compact and equal to $\{(x,0)\,\mid\,\rho(x)=1\}$ (in case $0\in [a,b]$). We are thus in the setting of an admissible Liouville structure that we discussed in section \ref{sec:Liouv}.
\end{rem}

\section{Cobordisms in Liouville Manifolds}\label{sec3}
\subsection{The definition of Lagrangian cobordisms}
Exact conical Lagrangian submanifolds $V\subset \IC\times W$ can have conical parts of different forms. 
To illustrate this we consider a Lagrangian suspension of a linear Hamiltonian deformation and a conical exact Lagrangian submanifold. This will motivate the definition of cobordisms in $\IC\times W$.

\begin{ex}\label{ex:susp}
Consider the (time-independent) Hamiltonian function $H:T^*S^1\to \IR$ on the cylinder $T^*S^1$ given by $H(q,p)=p$. It induces the Hamiltonian diffeomorphism $\varphi_t$ of $T^*S^1$ that wraps a fiber around the cylinder with constant velocity. By introducing a non-decreasing smooth function $\rho:\IR\to [0,1]$ which is $0$ for small $t$ and equal to $1$ for large $t$ the Hamiltonian diffeomorphism $\varphi_{\rho(t)}$ is generated by $H_t=\rho'(t)H$. This new Hamiltonian diffeomorphism wraps each fiber around the cylinder exactly once. Looking at one particular fiber $L\subset T^*S^1$ over $q_0\in S^1$, the Lagrangian\begin{figure}
    \centering
    \begin{tikzpicture}
\draw [blue] (0,0)  -- node[pos=0.5, below=1mm] {$L$} (3,0);
\draw [blue](6,0) -- node[pos=0.5, below=1mm] {$L$} (9,0);
\path [fill=green] (3,-3) rectangle (6,3);
\node at (4, 1.5) {V};

\draw[fill=blue!80, draw=blue!80]    (3,0) .. controls (4,0) and (4,1) .. (4.5,1)--(4.5, 1).. controls (5,1) and (5,0) .. (6,0)--(6,0) .. controls (5,0) and (5,-1) .. (4.5,-1)--(4.5, -1).. controls (4,-1) and (4,0) .. (3,0);

\draw[dashed] (4.5,0) ellipse (3.5 and 2.5);

\end{tikzpicture}
    \caption{The green region (plus the real axis) is the projection of the Lagrangian suspension $V_{L,H}\subset \IC\times T^*S^1$ in Example \ref{ex:susp} to the complex plane $\IC$. Clearly, it goes  to infinity in the imaginary and the real direction of the plane $\IC$. To extract more information we focus on the blue region. It is the projection of $V\cap \left(\IC\times T_{|p|\leq 1}^*S^1\right)$ to $\IC$, where $T_{|p|\leq 1}^*S^1$ denotes the disk cotangent bundle of the circle $S^1$. In the blue area only the two real directions go to infinity. The area of the blue region is equal to $1$. Its shape depends on the velocity (which is parametrized by $\rho'$) of how fast the fiber is wrapped around the cylinder.}
    \label{fig:cobex}
\end{figure}
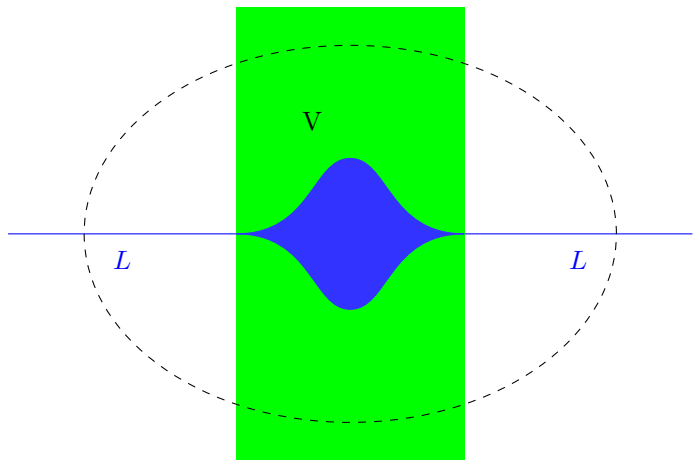 suspension $V_{L,H}$ is the image of the map $\IR^2\to \IC\times T^*S^1$ given by $$(t,p)\mapsto (t,\rho'(t)p, q_0+\rho(t), p).$$ 

\end{ex}
The following observation is true for any Lagrangian suspension of an exact conical Lagrangian submanifold by a linear Hamiltonian: As in Remark \ref{rem:compactliou} we consider a Liouville form that admits compact Liouville domains in $\IC$.
The Liouville boundary $\partial(\IC\times W)$ can be covered by $\partial \IC\times W^{int}$ and $\IC^{int}\times \partial W$. Looking at Figure \ref{fig:cobex} we see that if we choose the Liouville domain in $\IC$ large enough then $\partial \IC\times W^{int}$ is only hit by the part of $V$ where $V$ looks like conicalizations of paths in $\IC$ and the Lagrangian submanifold $L\subset W$.\\[2mm]
This suggests the following definition:

\begin{defi}\label{def:cobordism}
Let $(W,\lambda_W, \ff)$ be a stopped Liouville manifold. A \emph{Lagrangian cobordism} $V\subset \IC\times W$ is an exact conical Lagrangian submanifold with respect to $\lambda_\IC+\lambda_W$ where $\lambda_\IC$ is a Liouville form on $\IC$ that admits a compact Liouville domain and $d\lambda_\IC$ induces the standard orientation. Moreover, we require that there are choices of Liouville domains $\IC^{int}\subset \IC$ and $W^{int}\subset W$ such that
\begin{itemize}
    \item $V$ is conical outside $\IC^{int}\times W^{int}$,
    \item $V\cap \left(\IC^{cone}\times W^{int}\right)$ consists of only of finitely many submanifolds that are conicalizations $\gamma_j\tilde\times L_j$ of product Lagrangian submanifolds of conical rays $\gamma_j$ in $\IC$ and exact Lagrangian submanifolds $L_j\subset W$. These conicalized products are referred to as the \emph{ends} of the cobordism.
    \item $V$ is disjoint from $\IC\times\ff$ except possibly at the ends of the cobordism. Moreover, all the Lagrangian submanifolds $L_j$ that are ends of $V$ do not intersect the stop $\ff$. 
\end{itemize} 
\end{defi}

\begin{rem}
The definition of a Lagrangian cobordism depends on a Liouville form on $\IC$ and a Liouville domain $\IC^{int}\subset \IC$. 

The Liouville structure on $\IC$ is not very important information in the sense that it can be deformed to any other preferred choice of Liouville structure on $\IC$ (see Proposition \ref{prop:changeLiou}).

However, the choice of a Liouville domain can be crucial: The same Lagrangian submanifold $V\subset \IC\times W$ can serve as two Lagrangian cobordisms with a different number of ends. Enlarging the Liouville domain $\IC^{int}$ on $\IC$  can lead to the consequence of not seeing one end. An important example is described by certain stabilization cobordisms that can be viewed as null-cobordisms when enlarging the Liouville domain on one side (more in section \ref{sec:null-cob}). We usually do not want to forget these fake ends so that we could also glue on them.
\end{rem}

\begin{rem}
Our notion of cobordism is somewhat analog to \cite{BC}. The conical property of a cobordism $V$ allows us to focus on $V\cap \left(\IC\times W^{int}\right)$ instead of all of $V$. As noticed in Remark \ref{rem:nonintersect}, to study potential intersections of conical Lagrangian submanifolds and conicalizations of product Lagrangian submanifolds in $\IC\times W$, it is enough to look at the projection of $V\cap \left(\IC\times W^{int}\right)$ to $\IC$. 
\end{rem}

\subsection{Liouville structures on $\IC$ continued}\label{sec:Liouvcont}
Let $(W, \lambda_W, \ff)$ be a stopped Liouville manifold. A Lagrangian cobordism $V\subset \IC\times W$ with exact conical ends should be of the same class of Lagrangian submanifolds in $\IC\times W$.  It should be exact and conical with respect to a Liouville form $\lambda_\IC+\lambda_W$ where the Liouville structure $\lambda_\IC$ on $\IC$ admits a compact Liouville domain and is positively oriented. 

In this section we argue that the choice of a Liouville structure on $\IC$ does not restrict the set of cobordisms. More specifically:
\begin{prop}\label{prop:changeLiou}
Suppose $\lambda_0$ is an admissible Liouville structure on $\IC$, i.e. has compact skeleton and induces the standard orientation. Suppose $V_0$ is a Lagrangian cobordism that is exact and conical with respect to a Liouville structure $\lambda_\IC^0+\lambda_W$ on $\IC\times W$ and let $\lambda_\IC^1$ be another admissible Liouville structure on $\IC$. Then there is a Lagrangian cobordism $V_1$ that is exact and conical with respect to the Liouville structure $\lambda_\IC^1+\lambda_W$ on $\IC\times W$ with the same ends as $V$.

More precisely, $V_0$ and $V_1$ are isotopic through a family of cobordisms that are defined with respect to $\lambda_\IC^t+\lambda_W$, where $\lambda_\IC^t$ is a smooth family of admissible Liouville structures connecting $\lambda^0_\IC$ to $\lambda^1_\IC$.
\end{prop}

Note that we do not claim that this isotopy is unique in any sense. We are solely interested in its existence so that we can switch between Liouville forms chosen for $\IC$, and still have a cobordism with the same ends.
\begin{proof}
The strategy of the proof is to first deform the ends of $V_0$ and then deform the interior of $V_0$.

Fix a Liouville domain $W^{int}$ as in the definition of $V_0$. The space of admissible Liouville structures on $\IC$ is convex. A standard argument as in \cite{CE} Proposition 11.8 shows that by choosing an exhaustion of Liouville domains of $\IC$ and finding contactomorphisms on the corresponding Liouville boundaries, we can produce a diffeomorphism $\psi:\IC\to \IC$ such that
\begin{itemize}
    \item $\psi^*(\lambda_\IC^0)-\lambda_\IC^1$ is exact,
    \item $\psi^*(\lambda_\IC^0)=\lambda_\IC^1$ outside a large compact set $K\subset \IC$,
    \item $\psi=id$ in a smaller compact set $K'\subset K$ which strictly contains the projection $\pi_\IC(V_0\cap \left(\IC\times W^{int}\right))$ when ignoring the projection of the ends.
\end{itemize}
Define $\widetilde \lambda_\IC^0=\psi^*(\lambda_\IC^0)$ and $\widetilde{V_0}=(\psi^{-1}\times id)(V_0)$. The three properties above ensure that
\begin{itemize}
    \item $\widetilde V_0$ is exact with respect to $\widetilde \lambda_\IC^0+\lambda_W$,
    \item $\widetilde V_0$ consists outside $K$ only of ends parametrized by conical paths with respect to $\lambda_\IC^1$,
    \item $\widetilde V_0=V_0$ intersects $K'\times \partial W$ in a Legendrian with respect to $\widetilde\lambda_\IC^0$. 
\end{itemize}
The exact Lagrangian submanifold $\widetilde V_0$ may not be conical with respect to $\widetilde \lambda_\IC^0$ but $\partial \widetilde V_0$ is Legendrian and can be conicalized by Lemma \ref{makeconical}. Denote its conicalization also by $\widetilde V_0$.

\vspace{2mm}
We are in the situation of a Lagrangian cobordism $\widetilde {V_0}$ with respect to $\widetilde\lambda_\IC^0+\lambda_W$ and want to construct a Lagrangian cobordism with respect to the Liouville form $\lambda_\IC^1+\lambda_W$. By construction, we have $\widetilde\lambda_\IC^0=\lambda_\IC^1$ on $\IC^{cone}$. That means the family of Liouville structures $\lambda_t=t\widetilde\lambda_\IC^0+(1-t)\lambda_\IC^1+\lambda_W$ on $\IC\times W$ is constant on $\IC^{cone}\times W$. Again by applying Proposition 11.8 from \cite{CE} there is a isotopy of diffeomorphisms $\Psi_t:\IC\times W\to\IC\times W$ such that 
\begin{itemize}
    \item $\Psi_t^*\lambda_t-\lambda_0$ is exact,
    \item $\Psi_t^*\lambda_t=\lambda_0$ outside a compact set,
    \item $\Psi_t$ is constant on the ends of $\widetilde {V_0}$.
\end{itemize}
That the last property is true can be seen by inspecting the proof of Proposition 11.8 in  \cite{CE}. The contactomorphisms on the Liouville boundary $\partial (\IC\times W)$ which define $\Psi_t$ are constant on the part $\partial \IC\times W^{int}$ which only meets the Lagrangian cobordism $\widetilde{V_0}$ in the ends. Compare to the remark 2.2.3 in \cite{GEI} about the proof of Gray's Stability theorem for contact forms. 

Consequently, $V_t=\Psi_t(\widetilde {V_0})$ is a family of Lagrangian cobordisms with respect to $\lambda_t$. In particular, $V_1$ is a Lagrangian cobordism with respect to $\lambda_1=\lambda_\IC^1+\lambda_W$.
\end{proof}

\subsection{Planar isotopies}\label{sec:planar}
We want to be able to move a Lagrangian cobordism $V\subset (\IC\times W, \lambda_\IC+\lambda_W)$ by a Hamiltonian diffeomorphism of contact type $\varphi_t$ in $(\IC, \lambda_\IC)$. In particular, we want to rotate the ends of cobordisms, to translate cobordisms in $\IC$ and to make two cobordisms disjoint. In general, the image of $(\varphi_t\times id_W)(V)$ is not conical with respect to $\lambda_\IC+\lambda_W$. However, we have the following two lemmas:

\begin{lem}\label{lemma1}
Let $V\subset \IC\times W$ be a Lagrangian cobordism. Suppose $\varphi_t$ is a compactly supported Hamiltonian isotopy of $\IC$. Then there is a Hamiltonian diffeomorphism $\Phi_t$ of contact type of $\IC\times W$ which is equal to $\varphi_t\times id_W$ over any previously fixed compact subset of $\IC\times W$.  
\end{lem}

\begin{proof}
Suppose the function $H_t:\IC\to\IR$ generates $\varphi_t$. Define a contact Hamiltonian on the contact manifold  $\partial(\IC\times W)$ by $H_t$. It is supported in $\IC^{int} \times \partial W\subset \partial(\IC\times W)$. From this construct a linear Hamiltonian $\IC\times W\to\IR$ which is equal to $H_t$ on $\IC^{int}\times W^{int}$ and linear with respect to $\lambda_\IC+\lambda_W.$ Applying the corresponding Hamiltonian diffeomorphism of contact type to $V\subset \IC\times W$ does not move the ends.\hspace{2mm}
\end{proof}

\begin{lem}\label{lemma2}
Let $V\subset \IC\times W$ be a Lagrangian cobordism with ends modelled on conical rays $\gamma_j\subset \IC$ outside some Liouville domain $\IC^{int}\subset\IC$. Suppose $\varphi_t$ is a Hamiltonian isotopy of contact type of $\IC$ supported outside $\IC^{int}$. Then there is a Hamiltonian diffeomorphism $\Phi_t$ of contact type of $\IC\times W$ such that the ends of $\Phi_t(V)$ are modelled on the rays $\varphi_t(\gamma_j)$.
\end{lem}

\begin{proof}
Let $\IC^{int}\times W^{int}$ be a (cornered) Liouville domain as in the definition of the Lagrangian cobordism $V$. The Lagrangian submanifolds $\widetilde{V}_t=(\varphi_t\times id_W)(V)$ are not conical, but intersect $\partial(\IC^{int}\times W^{int})$ in a Legendrian for all $t$. By Lemma \ref{makeconical} for each $t$ we can find a cobordism $V_t$  which is equal to $\widetilde V_t$ in $\IC^{int}\times W^{int}$.  As the $V_t$ are a family of exact Lagrangian submanifolds the isotopy $V_t$ is induced by a Hamiltonian diffeomorphism $\Phi_t$ of contact type starting at $V_0=V$ by Remark \ref{ambient}. Note that the ends of $V_t$ are indeed modelled on the rays $\varphi_t(\gamma_i).$
\end{proof}

For arbitrary diffeomorphism of contact type $\varphi_t$ in $\IC$ we can combine both lemmas. Hamiltonian diffeomorphisms $\Phi_t$ that arise in this way are called \emph{planar isotopies}\footnote{Compare this to the equivalent notion of \emph{horizontal isotopies} in \cite{BC2}. Ends of Lagrangian cobordism can be rotated as long as their order stays the same.}.

Define an equivalence relation on the set of cobordisms:
\begin{defi}
Let $V$ and $V'\subset \IC\times W$ be two Lagrangian cobordisms. By Propositon \ref{prop:changeLiou} we can make them exact and conical in $\IC\times W$ with the standard radial Liouville structure on $\IC$. They are called \emph{planar-isotopic} if they are isotopic by a planar isotopy.
\end{defi}

\subsection{Concatenation of cobordisms}\label{sec:concat}
\begin{figure}
    \centering
    \begin{tikzpicture}

\draw[dashed] (0,0) ellipse (6 and 4);

\draw[fill=blue, draw=blue] (2,0) ellipse (1 and 1.2);
\draw[fill=blue, draw=blue] (-2,0) ellipse (1.2 and 1);
\node at (2, 0) {$V'$};
\node at (-2, 0) {$V''$};

\draw[blue] (-2,0)--node[pos=0.4, below=1mm] {$L_{j_0}$}(2,0);
\draw[blue] (-2,0)--node[pos=0.7, left=1mm] {$K_{1}$}(-2,4);
\draw[blue] (-2,0)--node[pos=0.6, left=1mm] {$K_{m}$}(-2,-4);
\draw[blue] (-2,0)--(-6,2);
\draw[blue] (-2,0)--(-6,-2);
\draw[blue] (2,0)--node[pos=0.7, left=1mm] {$L_{1}$}(2,4);
\draw[blue] (2,0)--node[pos=0.7, left=1mm] {$L_{n}$}(4,-3.7);
\draw[blue] (2,0)--(1,-4.1);
\draw[blue] (2,0)--node[pos=0.8, below=1mm] {$L_{0}$}(6.4,0);

\node[shape = circle,fill = black, inner sep=2pt] at (0,4) {};
\node[shape = circle,fill = black, inner sep=2pt] at (0,-4) {};
\node[shape = circle,fill = black, inner sep=2pt] at (-6,0) {};
\node[shape = circle,fill = black, inner sep=2pt] at (-4,3) {};
\node[shape = circle,fill = black, inner sep=2pt] at (-4,-3) {};
\node[shape = circle,fill = black, inner sep=2pt] at (4,3) {};
\node[shape = circle,fill = black, inner sep=2pt] at (2.3,-3.68) {};
\node[shape = circle,fill = black, inner sep=2pt] at (5.5,-1.6) {};

\draw (1,4.1) ..node[pos=0.2, left=1mm] {$D_{i_0}$}  controls (-1,2) and (1,-2) .. (-1,-4.1);
\end{tikzpicture}
    \caption{Concatenation of two cobordisms $V', V''$ along a common end $L_{j_0}$.}
    \label{fig:cobcon}
\end{figure}
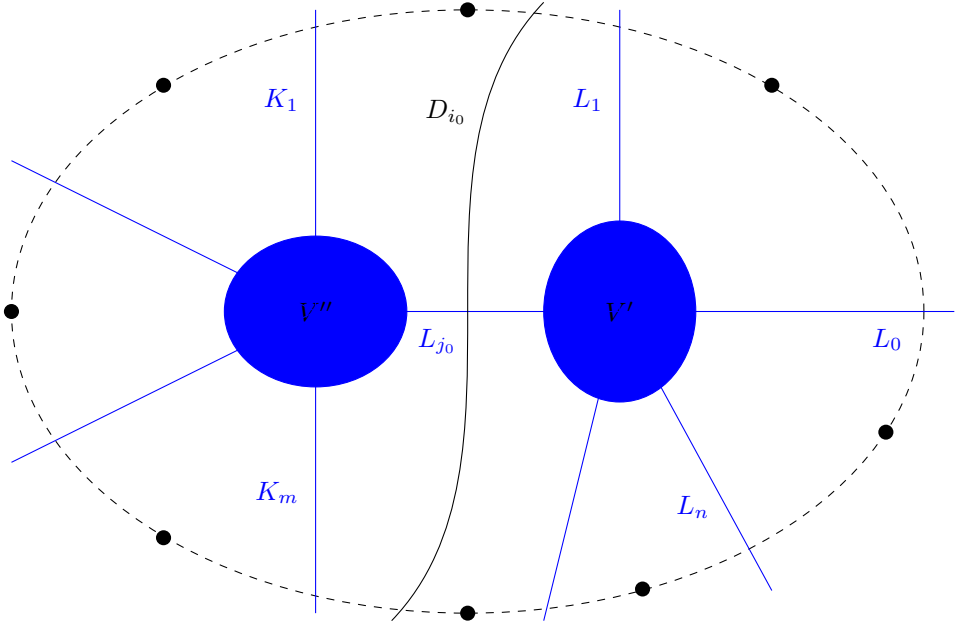
Suppose $V', V''\subset \IC\times W$ are two Lagrangian cobordisms that have a common Lagrangian end $L_{j_0}$. Then they can be concatenated (see Figure \ref{fig:cobcon}) to a Lagrangian cobordism $V\subset \IC\times W$  as follows: Change the Liouville structure on $\IC$ such that $V'$ and $V''$ are defined with respect to the standard radial Liouville structure $\lambda_\IC^0=\frac 12 (xdy-ydx)$ on $\IC$ (using Proposition \ref{prop:changeLiou}). Let $\epsilon>0$. Deform $V'$ and $V''$ separately by planar isotopies such that 
$$V''\cap \left(\IC_{\mathrm{Re}>-\epsilon}\times W\right)= (-\epsilon, \infty)\tilde\times L_{j_0} \qqet V'\cap \left(\IC_{\mathrm{Re}< \epsilon}\times W\right)= (-\infty,\epsilon)\tilde\times L_{j_0}.$$ We then can glue the cobordisms along $(-\epsilon, \epsilon)\tilde\times  L_{j_0}.$

\subsection{Examples of cobordisms}
In this section, we give the three main examples of cobordisms: Suspension cobordisms, surgery cobordisms and null-cobordisms. Moreover, we include an interesting example of non-cobordisms.

\subsubsection{Suspension cobordisms}\label{sec:susp} Let $V_{L,H}$ be a Lagrangian suspension of a conical exact Lagrangian submanifold $L$ and a linear Hamiltonian $H_t$ as defined in section \ref{sec:Lsusp} which vanishes for large $|t|$. We have already explained in Remark \ref{rem:compactliou} and Example \ref{ex:susp} how to change the Liouville form $-ydx$, such that $\IC$ admits a compact Liouville manifold. 

For any choice of Liouville domain $W^{int}\subset W$, the set $\pi_{\IC}(V_{L,H}\cap (\IC\times W^{int}))$ is contained in a compact set except the two ends lying on the negative and the positive part of the real line $\IR\subset \IC$. This is because $$(t,w)\mapsto H_t(\varphi_t(w))$$ is bounded for $(t, w)\in \IR\times W^{int}$. Thus Lagrangian suspensions $V_{L,H}$ as above are examples of Lagrangian cobordisms. 

\begin{rem}\label{ambient}
We often will use the fact that an isotopy of exact conical Lagrangian submanifolds can be realized by an ambient Hamiltonian isotopy of contact type (see for example \cite{OH}). In other words, if two exact conical Lagrangian submanifolds are isotopic through a family of exact conical Lagrangian submanifolds then they are Lagrangian cobordant by a suspension. 
\end{rem}

\subsubsection{Surgery cobordisms}\label{sec:surgery}
Surgery of two transversely intersecting Lagrangian submanifolds goes back to \cite{LS,POL} and was studied extensively in the setting of cobordisms in \cite{BC}. More generally, in \cite{MW} the construction also includes cleanly intersecting Lagrangian submanifolds.

We summarize some relevant results in the case of clean intersections:
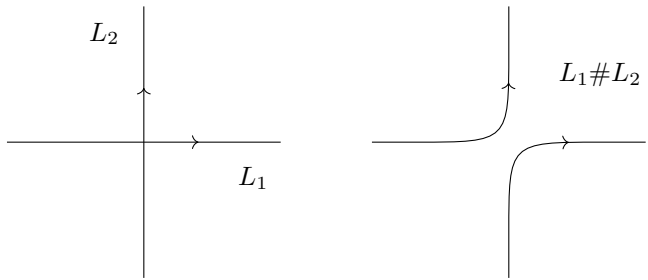
\begin{figure}
    \centering
    \begin{tikzpicture}[ scale=0.6]]
\begin{scope}[decoration={
    markings,
    mark=at position 0.7 with {\arrow{>}}}
    ] 
    
\draw[postaction={decorate}] (-3,0)  -- node[pos=0.9, below=2mm] {$L_1$} (3,0);
\draw[postaction={decorate}] (0,-3)  -- node[pos=0.9, left=2mm] {$L_2$} (0,3);

\draw[postaction={decorate}] (5,0)--(6,0) .. controls (8,0) and (8,0) .. (8,2)--(8,3);
\draw[postaction={decorate}] (8,-3)--(8,-2) .. controls (8,0) and (8,0) .. (10,0)--(11,0);
\node at (10,1.5){$L_1\# L_2$};
\end{scope}

\end{tikzpicture}
    \caption{The lowest dimensional example of surgery: Surgery of two (here oriented) curves.}
    \label{fig:surgery}
\end{figure}

\begin{lem}[\cite{MW} {\normalfont Corollary 2.22 and Lemma 6.2}]
Let $L_1 , L_2 \subset W$ be two exact Lagrangian submanifolds intersecting cleanly along some connected submanifold $D$ of $L_1$ and $L_2$. Then its surgery $L_1 \#_D L_2\subset W$ along $D$ is an exact Lagrangian submanifold.

Moreover, the trace $V\subset \IC\times W$ of this Lagrangian surgery is an exact Lagrangian submanifold.
\end{lem}
When the involved Lagrangian submanifolds are allowed to be conical instead of compact, the trace of the surgery may not be conical at infinity but is a product of possibly intersecting paths in $\IC$ and a conical Lagrangian submanifold. Therefore we can conicalize it by the conicalization Lemma \ref{makeconical}. The trace cobordism then is a Lagrangian cobordism between the three Lagrangians $L_1,L_2$ and $L_1\#_DL_2$.\\[2mm]
Surgery constructions that allow non-compact intersection of two exact conical Lagrangian submanifolds, for example intersections on the conical part, are of particular interest in Liouville manifolds. The construction given in \cite{HIC} fits also into the framework of this article by conicalizing products at infinity.

\subsubsection{Null-cobordisms}\label{sec:null-cob}
A Lagrangian submanifold $L$ which by the Liouville flow can be flowed away from any compact set are zero objects in the wrapped Fukaya category, i.e. all morphism spaces from and to this object are trivial (see Theorem \ref{vanishing}). We will prove that the stabilization cobordism $\IR\tilde\times L$ with two ends can serve as a null-cobordism for $L$. That is, a cobordism that has only one end equal to $L$ by choosing an appropriate Liouville domain on $\IC$.
\begin{prop}\label{prop:null}
Suppose $L$ is an exact conical Lagrangian submanifold which does not intersect the skeleton $\mathrm{sk}(W)$ of the Liouville manifold $W$. Then the Lagrangian submanifold $\IR\tilde\times L$ can serve as a null-cobordism for $L$. 
\end{prop} 
\begin{proof}
The conicalization of $\IR\times L$ must be invariant under the product Liouville flow of $\IC\times W$. After the small deformation region arising in the conicalization process, the Lagrangian $L$ gets pushed away more and more from the skeleton for large $|t|\in \IR\subset \IC$. We can choose a Liouville domain of $\IC$ such that on the negative end the conicalization $\IR\tilde\times L$ does intersect the part $\partial \IC\times W^{int}$ in a conicalization of the product $\gamma\times L$ and on the positive end has no intersection (see Figure \ref{fig:nullcob})\begin{figure}
    \centering
    \begin{tikzpicture}

\draw [blue] (1,0)  -- node[pos=0.5, below=1mm] {$\mathbb{R}\tilde\times L\cap \mathbb{C}\times W^{int}$} (6,0);

%\draw [red] (4.5,3)  .. node[pos=0.2, left=1mm] {$\alpha\tilde\times K$} controls (1,-1) and (1,1) .. (4.5,-3);
%\draw [red] (4.5,3)  .. node[pos=0.2, left=1mm] {$\alpha'\tilde\times K$} controls (8,-1) and (8,1) .. (4.5,-3);

\draw[dashed] (4.5,0) ellipse (2.5 and 1.5);
\end{tikzpicture}
    \caption{The conicalization $\IR\tilde\times L$ of the stabilization $\IR\times L$ of a Lagrangian submanifold $L$ which does not intersect the skeleton $\mathrm{sk}(W)$. Note that we just draw the projection of $(\IR\tilde\times L)\cap \left(\IC\times W^{int}\right)$ to $\IC$ in blue instead of $\IR\tilde\times L$. The Liouville boundary of $\IC$ dotted in black is chosen such that it yields a null-cobordism for $L$, i.e. that the cobordism $\IR\tilde\times L$ just has one end $L$. When we calculate Floer homology of $L$ with another test Lagrangian $K$ the theory in section \ref{sec4} shows that it is equivalent to Floer homology of $\IR\tilde\times L$ with $\alpha\tilde\times K$, where $\alpha$ intersects $\IR\subset \IC$ exactly once. Isotoping the curve $\alpha$ from one side of the chosen Liouville boundary of $\IC$ to the other side eliminates any intersections of $\IR\tilde\times L$ and $\alpha\tilde\times K$.}
    \label{fig:nullcob}
\end{figure}
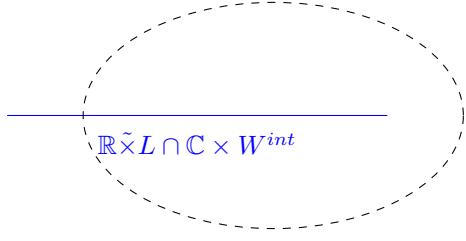.
\end{proof}

In case of a stopped Liouville domain $(W, \ff)$ the above constructed null-cobordism works if $L$ does not intersect the generalized skeleton
$$\mathrm{sk}(W)\cup\bigcup_{t\in \IR}\Psi_t(\ff)\subset W$$ where $\Psi_t$ denotes the Liouville flow. Equivalently, $\Psi_t(L)$ is disjoint from $\mathrm{sk}(W)$ and $\ff$ for all $t$. The Lagrangian $L$ can be flowed away from any compact subset of $W$ by the Liouville flow without flowing over the stop $\ff$.

\subsubsection{Some non-cobordisms}
Suppose $S_t$ is an isotopy of submanifolds in a compact manifold $M$. Then the conormal submanifolds $N^*S_t$ in $T^*M$ form an isotopy of exact conical Lagrangian submanifold which can be realized by an ambient Hamiltonian isotopy (see Remark \ref{ambient}). Hence, if $S_0$ and $S_1$ are smoothly isotopic submanifolds in $M$, the conormals $N^*S_0$ and $N^*S_1$ are Lagrangian cobordant by a Lagrangian suspension.

Trying to generalize this to conormals of  submanifolds in $M$ which are only smoothly cobordant instead of isotopic fails. The following example shows how to naively produce a conical Lagrangian submanifold that does not classify as a cobordism. 

\begin{ex}\label{ex:noncob}
This non-example was communicated to me by Baptiste Chantraine. It constructs an exact conical Lagrangian submanifold which looks almost like a cobordism between the  conormal $N^*S^1\subset T^*S^2$ of a circle $S^1\subset S^2$ and the empty set. The idea of the construction uses the null-cobordism of $S^1\subset S^2$ by a hemisphere in $S^2$. We will give the details of the construction shortly.

However, the Lagrangian submanifold $N^*S^1$ is not a zero object in $\mathcal W(T^*S^2)$ because when pairing with the zero section $0_{S^2}\subset T^*S^2$ wrapped Floer homology does not vanish: $$HW(O_{S^2},N^*S^1)\cong H(S^1).$$ This is proved by Morse theory in \cite{DUR}.
\begin{figure}\label{fig:noncob}
    \centering
    \begin{tikzpicture}
\draw [fill=green!80, draw=green!80] (0,-2)  .. controls (0,0) and (0,0) .. node[pos=0.5, left=1mm] {} (1,0)--(1,0) .. controls (0,0) and (0,0) .. (0,2)--(0,-2);
\draw [green] (0,0)  -- node[pos=0.5, below=1mm]{} (4,0);
\draw [green] (0,-2)  -- node[pos=0.5, left=1mm] {$V$} (0,2);
\end{tikzpicture}
    \caption{A conical exact Lagrangian submanifold in $\IC\times W$ which is not a cobordism as the ends going to negative and positive imaginary infinity are not product Lagrangian submanifolds. See Example \ref{ex:noncob} for the definition of $V$.}
\end{figure}
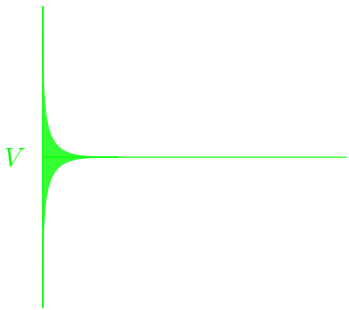

The construction: Fix an embedded copy of $S^1$ in $S^2$. To construct the \emph{inadmissible} cobordism between $N^*S^1$ and the empty set, we start by constructing an embedding $\IR^2\to \IR\times S^2$ in the following way:
map the unit ball $B^2\subset \IR^2$ onto $\{0\}\times B^2\subset\IR\times S^2$ where the latter $B^2$ is meant to be one of the hemispheres cut out by $S^1\subset S^2$. Next, carefully map the rest $\IR^2\setminus B^2$ (a cylinder) onto the cylinder $(0,\infty)\times S^1$ and then smoothen this construction.

Taking the conormal bundle of a submanifold always gives rise to an exact conical Lagrangian submanifold in the cotangent bundle.
Analyzing the conormal $V:=N^*(\IR^2)\subset T^*(\IR\times S^2)=T^*\IR\times T^*S^2$ of the constructed embedding $\IR^2\subset \IR\times S^2$ there is one cylindrical end on the positive real line which is equal to $\IR\times N^*S^1$. The projection $V\cap (\IC\times W^{int})$ to $\IC$ is equal to the union of the imaginary axis, the positive real axis and a smoothing area near the intersection (see Figure \ref{fig:noncob}). However, according to our definitions, it is not a Lagrangian cobordism: the preimage of a point $iy\in \IC$ on the imaginary axis is a hemisphere $B^2\subset S^2$ for any $y\in \IR$, hence  is not the conicalization of a product Lagrangian submanifold (and cannot be conicalized).  

This example generalizes in the following way: Let $C$ be a smooth cobordism in a compact manifold $M$ between submanifolds $S_0$ and $S_1$. Let $\overline{C}$ be the manifold obtained by completing the cobordism by cones $(-\infty, 0)\times S_0$ to the left end and $(0,\infty)\times S_1$ and denote these parts by $C^{c}$. Then we can construct an embedding $\overline{C}\to \IR\times M$ by using the inclusion $C\to \{0\}\times C\subset \{0\}\times M$ and the inclusion $C^{c}\to \IR\times M$ in a smooth manner. This again yields a non-cobordism in $\IC\times T
^*M$ from $N^*S_0$ to $N^*S_1$ which projects on the union of the two axes in $\IC$ (and some smoothening area). 
\end{ex}

\section{(Partially) Wrapped Floer Theory and Cobordisms}\label{sec4}
In this section we recall the relevant Floer theory following \cite{GPS1,GPS2} and add some comments about grading. Then before proving the main theorem \ref{thm11} we discuss the partially wrapped Fukaya category of the disk as it serves as a standard model of iterated mapping cones. 

\subsection{Floer theory of non-compact Lagrangian submanifolds}
The first instance of the open string analog of symplectic homology was when Abbondandolo and Schwarz in \cite{ABSC} defined wrapped Floer homology for fibers in cotangent bundles using a quadratic Hamiltonian. Two possible definitions of wrapped Floer homology and the wrapped Fukaya category of Liouville manifolds, were first sketched in \cite{FSS}. The approach using a direct limit over linear Hamiltonians was systematically introduced in \cite{AS}, the approach over quadratic Hamiltonians in \cite{ABO2}. Auroux \cite{AUR2,AUR} introduced a refinement by not allowing to wrap over certain subsets in the contact boundary of the Liouville manifold. This refined version is called partially wrapped Fukaya category (see also \cite{SYL}). A third method to define the partially wrapped Fukaya category was introduced in \cite{GPS1}. The strategy resembles the linear setup from \cite{AS}, however, its definition does not involve any Hamiltonian functions. It is directly defined over all exact Lagrangian deformations which are positive at infinity and then localizes with respect to Floer continuation elements. Moreover, the setting in \cite{GPS1} also allows certain Liouville manifolds with boundary near which the wrapping stops. These manifolds are called \emph{Liouville sectors}. All these constructions are expected to produce equivalent $A_\infty$-categories. However, this is not yet rigorously established according to the introduction in \cite{CHA2}.

In dimension $2$, the discussion reduces to closed exact curves and arcs in surfaces with boundary. The boundary may be equipped with stops.  Stops  are points on the boundary, over which arcs are not allowed to wrap. The partially wrapped Fukaya categories for some symmetric products of surfaces are calculated in Auroux's papers \cite{AUR2,AUR} by relating the computation to Lefschetz fibrations and Seidel's generation result by thimbles \cite{SEI}. There is another equivalent approach \cite{BOC,HKK}, called the topological Fukaya category, which is, as the name suggests, purely combinatorial and topological. Moreover, \cite{DK} show that triangulations of a surface lead very naturally to the definition of the topological Fukaya category and its derived triangulated structures. Some explicit computation and pictures of partially wrapped Fukaya categories of some Riemann surfaces can be found in \cite{LP}.

For our purpose, we mainly rely on the wrapped Fukaya category and theory developed in \cite{GPS1}. The only explicit calculation we need are exact triangles in the partially wrapped Fukaya category of a disk stopped at finitely many boundary points (see section \ref{sec:wrapC}).

\subsection{Partially wrapped Floer homology}\label{sec:wrappedhomology}
We give an overview of the construction of partially wrapped Floer theory assuming the definition and the properties of Floer homology $HF(L,K)$ of two transverse, exact conical Lagrangian submanifolds $L,K$ in a Liouville manifold $W$ are known. One key difference to the definition of Floer homology in the compact case is that to get Gromov compactness for holomorphic curves with Lagrangian boundary conditions we need to restrict to almost complex structures that are conical at infinity, i.e. almost complex structures $J$ on $W$ satisfying $\lambda\circ J=dr$  at infinity. Below we use the cohomological convention, instead of the homological one used in \cite{GPS1}.

\subsubsection{Partially wrapped Floer homology} A \emph{stop} $\ff$ in a contact manifold $(Y, \alpha)$ is a closed subset of $Y$.
An isotopy of Legendrian submanifolds $\Lambda_t$ in a stopped contact manifold $(Y, \alpha, \ff)$ is called \emph{positive} if $\alpha(\partial_t \Lambda_t)\geq 0$ holds and  $\Lambda_t$ does not intersect $\ff$ for all $t$.
For any fixed Legendrian $\Lambda$ in $Y$ consider the partially ordered set $\Lambda\rightsquigarrow \cdot$ of arbitrary Legendrian isotopies $\Lambda\rightsquigarrow\Lambda^+$ where the order is given by $$\Lambda\rightsquigarrow\Lambda_0\quad<\quad\Lambda\rightsquigarrow\Lambda_1$$ if and only if there is a positive Legendrian isotopy $\Lambda_0\rightsquigarrow\Lambda_1$ such that 
the composite isotopy $\Lambda \rightsquigarrow\Lambda_0\rightsquigarrow\Lambda_1$ is  isotopic to $\Lambda \rightsquigarrow\Lambda_1$ by a path of Legendrian isotopies with fixed ends $\Lambda$ and $\Lambda_1$. 

These notions  translate directly to Liouville manifolds $W$ by restricting to the contact boundary. A \emph{stop} $\ff$ of a Liouville domain $(W^{int}, \lambda)$ is a closed subset in $\partial W$. In Liouville manifolds we should think of the stop $\ff$ as the completion to the cone, i.e. $[1,\infty)\times \ff\subset [1,\infty)\times \partial W = W^{cone}$. A \emph{stopped Liouville manifold} is a triple $(W,\lambda,\ff)$.
An exact Lagrangian isotopy $L\rightsquigarrow L^+$ is called \emph{positive} if at infinity it is a positive Legendrian isotopy that does not intersect the stop $\ff$. More precisely, for all large enough Liouville domains $W^{int}$, the restriction of $L\rightsquigarrow L^+$ to the Liouville boundary $\partial W$ is a positive Legendrian isotopy $\partial L \rightsquigarrow \partial L^+$ that does not intersect the stop $\ff$.

Let us recall some properties of Floer homology $HF(L,K)$ for conical exact Lagrangian submanifolds $L, K$ in a stopped Liouville manifold $(W,\lambda,\ff)$ (see \cite{GPS1} for details). For transverse $L,K$ it is the homology of the chain complex generated by intersection points of $L$ and $K$ with differential counting holomorphic strips with boundary on $L$ and $K$. For non-transverse exact conical Lagrangian submanifolds $L,K$ Floer homology $HF(L,K)$ is defined as $HF(L^+, K)$ where $L^+$ is an unspecified positive Lagrangian pushoff of $L$ (i.e. $L^+$ is a isotopic to $L$ by a small positive exact Lagrangian isotopy and  such that $L^+,K$ are transverse). This construction does not depend on the positive pushoff $L^+$. In particular, the self Floer homology $HF(L,L)$ is unital with a unit coming from a map of associative unital algebras from the singular cohomology $H(L)$ to the self Floer homology $HF(L,L)$. 

For any positive exact Lagrangian isotopy $L_t$ there is a \emph{continuation element} $c(L_t)\in HF(L_1, L_0)$ in Floer homology which for small positive exact Lagrangian isotopies corresponds to the unit in $HF(L_0,L_0)$ and for larger positive isotopies is defined by the composition of continuation elements of small isotopies. For Lagrangian submanifolds $K$ transverse to $L_0$ and $L_1$, the continuation element $c(L_t)\in HF(L_1, L_0)$ can be used to define a \emph{continuation map} $HF(L_0, K)\to HF(L_1, K)$ defined as postcomposing with the continuation element $c(L_t)\in HF(L_1, L_0)$.

We define wrapped Floer homology of exact conical Lagrangian submanifolds $L$ and $K$ in $W$ as the direct limit $$HW(L, K):=\varinjlim_{L \rightsquigarrow L^+}HF(L^+, K)$$
where the direct limit is taken over all exact Lagrangian isotopies $L\rightsquigarrow L^+$, where $L^+$ and $K$ are transverse. The isotopies are partially ordered by the ordering at infinity on Legendrian isotopies $\partial L\rightsquigarrow \partial L^+$, the maps given by continuation maps. 

Wrapped Floer homology can equivalently be defined by negatively wrapping $K$ or even by simultaneously wrapping $L$ positively and $K$ negatively. As a consequence we record that wrapped Floer homology $HW(L_t, K)$ is invariant under exact conical Lagrangian deformations $L_t$. Similarly, we can also vary $K$ through exact conical Lagrangian submanifolds. 

\subsubsection{Vanishing of wrapped Floer homology} In this section we record two important vanishing results for wrapped Floer homology.

Two disjoint Lagrangian submanifolds may not have vanishing wrapped Floer homology. However, if the wrapping can be done in a disjoint way, then  Floer homology vanishes. In more detail: In \cite{GPS1} it is shown that for each Lagrangian submanifold $L$ there is a countable \emph{cofinal} sequence
\begin{align}\label{cofinal}
L = L^{(0)}\rightsquigarrow L^{(1)}\rightsquigarrow L^{(2)}\rightsquigarrow\ldots
\end{align}
in the partially ordered set of Lagrangian isotopies $L\rightsquigarrow L^+$ such that Floer homology can be calculated as  $HW(L,K)=\varinjlim_{j}HF(L^{(j)}, K)$ for every test Lagrangian $K$ which is transverse to all $L^{(j)}$. In case there is a cofinal wrapping of $L$ in which all Lagrangians are disjoint from $K$, the wrapped Floer homology vanishes. Equivalently, $L$ and $K$ are disjoint and there are no Reeb orbits from $\partial L$ to $\partial K$ that avoid the stop $\ff$. 

\vspace{2mm}

Another vanishing property of wrapped Floer homology that does not follow directly from the definition is the following vanishing result stated in \cite{CHA2}.
\begin{thm}[\cite{CHA2} {\normalfont Sections 7 and 10.2}]\label{vanishing}
Let $L\subset W$ be a Lagrangian submanifold that can be moved away from any compact subset of $W$ by the Liouville flow but not flowing over the stop $\mathfrak f$, i.e. $L$ is disjoint from $\mathrm{sk}(W)\cup\bigcup_{t\in \IR}\Psi_t(\ff)$ where $\Psi_t$ denotes the Liouville flow.
Then $L$ has zero wrapped Floer homology with any other Lagrangian submanifold in $W$. 
\end{thm}

This theorem can be deduced by using the null-cobordism from section \ref{sec:null-cob} and the cone decomposition Theorem \ref{thm11}.

\subsection{Brane structures}
The Fukaya category (see next section \ref{sec:part} for the definition of the Fukaya category) will be defined by equipping Lagrangian submanifolds with different extra structure. An object in the Fukaya category consists of tuples $$(L,f,\theta,P)$$ where $L$ is a Lagrangian submanifold, $f:L\to \IR$ is a \emph{primitive} of $L$, $\theta:L\to\IR$ a {grading} and $P$ a \emph{Pin structure}. These structures exist under additional assumptions on $W$ and on $L$.

In section \ref{sec:conical} we already discussed primitives for Lagrangian submanifolds $L\subset W$. For exact Lagrangian submanifolds it is not necessary to work with Novikov coefficients in Floer theory. Moreover, there is a simple argument showing that holomorphic curves with Lagrangian boundary condition cannot escape to infinity. Gradings equip Floer theory with a $\IZ$-grading. We discuss gradings in the next section \ref{sec:grading}. Pin structures allow $\IZ$-coefficients in Floer theory by orienting moduli spaces (see \cite{SEI} for more details). One can additionally add local systems as an extra structure. They are relevant in discussing mirror symmetry or stability conditions \cite{HKK,HEN}.
\subsubsection{Grading}\label{sec:grading}
We discuss first the general case and then restrict to the two-dimensional case to complement section \ref{sec5}, where we discuss compute cobordism groups in dimension two. 
\vspace{2mm}

Let $\Lambda_n$ be the Grassmanian of all linear Lagrangian subspaces in $(\IC^n, \omega_{\mathrm{std}})$. Choosing an orthonormal basis of a Lagrangian subspace $L\subset \IC^n$ gives rise to  a unitary basis of $\IC^n$, i.e. an element in $U_n$. Changing the orthonormal basis of $L$ corresponds to an element in $O_n$. This identifies $\Lambda_n$ with $U_n/O_n$. The square of the determinant map $\det:U_n\to S^1$ descends to a map $\det^2:\Lambda_n\cong U_n/O_n\to S^1$.  The  preferred generator $\mu$ of $H^1(\Lambda_n)\cong S^1$ is called the \emph{Maslov class.}

Given a symplectic manifold $(W^{2n},\omega)$ one has a fiber bundle $\mathcal L^W\to  W$ with fibers $\mathcal L^W_z=\Lambda_n(T_zW, \omega_z)$. 
To calculate determinants, let $J$ be an $\omega$-compatible almost complex structure on $W$. Then $TW$ is a complex $n$-bundle and $\det(W)=\bigwedge ^{n,0}(T^*W,J)$ is a complex line bundle, called the \emph{determinant bundle.}
 The line bundle $\det(W)^{\otimes 2}$ is trivial iff the Chern class $c_1(\det(W)^{\otimes 2})=0\in H^2(W,\IZ)$, or equivalently $2c_1(W)=0\in H^2(W,\IZ)$. If this condition is satisfied,  the line bundle $\det(W)^{\otimes 2}$ admits a non-vanishing section $\Theta:W\to\det(W)^{\otimes 2}$, called a \emph{grading} of $W.$ Moreover, we get a map $\alpha_\Theta:\mathcal{L}^W\to S^1$ by $$\alpha_\Theta(z, \Lambda)=\frac{\Theta_z(v_1\wedge \cdots \wedge v_n)^2}{|\Theta_z(v_1\wedge \cdots \wedge v_n)^2|},$$ called the \emph{squared phase map}.
 
Any Lagrangian embedding $i_L:L\to M$ induces a map $s_L:L\to \mathcal L^W\vert_ L$ called the \emph{Gauss map}, given by $z\mapsto D_zi_L(T_zW)$. A \emph{grading} on $L$ is a map $\theta_L:L\to \IR$ which lifts $\alpha_\Theta\circ s_L:L\to S^1$. Such a lift exists iff $(\alpha_\Theta\circ s_L)_*\pi_1(L)=0\subset \pi_1(S^1),$ i.e. if the Maslov class $\mu_L:=(\alpha_\Theta\circ s_L)^*(\mu)\in H^1(L)$ vanishes.

A grading $\theta_L$ on $L$ can be shifted: Denote $L[1]$ the Lagrangian submanifold $L$ equipped with the grading $\theta_L-1.$

\begin{rem}(Dimension 2)
Let $\Sigma$ be an orientable non-compact surface. To get a grading we are interested in non-vanishing holomorphic quadratic differentials, i.e. sections of $(\bigwedge ^{1,0}T^*W)^{\otimes 2}$. As holomorphic line bundles on non-compact surfaces are trivial, there always exists a non-vanishing quadratic differential $\Theta$. If this section is locally given by $q(z)dz\otimes dz$, it induces a volume form $|q(z)|dx\wedge dy$ and a Riemannian metric $\sqrt{|q(z)|}dx\otimes dy$. 

One can perform a coordinate change $z\mapsto \sqrt{q(z)}dz$ such that the quadratic differential is $dz\otimes dz$. Transition functions on overlaps of two charts are then forced to be \emph{half-translations} $z\mapsto \pm z+c$ for $c\in \IC$. The induced metric is flat. Moreover, it comes with a choice of horizontal foliation given by horizontal lines in a local chart $\IC\cong \IR^2$. See the textbook \cite{STR} for more details.

In fact, there is a bijection between the following structures
\begin{itemize}
    \item complex structures together with a non-vanishing holomorphic quadratic differential
    \item flat Riemannian metrics together with a covariantly constant foliation and an orientation.
\end{itemize}

The stronger assumption of admitting a holomorphic differential leads to the notion of a \emph{translation surface} instead of a \emph{half-translation surface}. On overlapping charts the transitions are $z\mapsto z+c$ for some $c\in \IC$. Moreover, the corresponding horizontal foliation lifts to an oriented foliation, a vector field. To avoid the use of local coefficients with respect to an orientation double cover on $\Sigma$ (see \cite{HKK}) we can restrict to foliations coming from a vector field (see \cite{LP2} for a characterization when a particular foliation lifts).

A grading of a curve $\gamma:L\to \Sigma$ is a map $\theta:L\to \IR$, such that $\Theta(\gamma'(t), \gamma'(t))\in \IR_{>0}e^{2\pi i\theta(t)}$ (when identifying the holomorphic cotangent bundle with the real cotangent bundle). More generally, in terms of a foliation $\eta:\Sigma\to \mathbb P (TX)$, a grading is a homotopy class of paths between $\eta(\gamma(t))$ and $\gamma'(t)$ varying continuously with $t\in L$. The grading induced by a holomorphic quadratic differential is given by the paths from $\eta(\gamma(t))$ to $\gamma'(t)$ by $s\mapsto e^{\pi i\theta(t)s}$ for every $t$ assuming that the foliation $\eta$ is horizontal. A graded curve is oriented by $t\mapsto e^{\pi i\theta(t)}$. Note that the shift $\gamma[1]$ has reversed orientation.

An intersection point $p$ of two transversely intersecting graded curves  $(\gamma_1,\theta_1)$ and $(\gamma_2,\theta_2)$ can be $\IZ$-graded by $$i_p(\gamma_1,\gamma_2) = \lceil\theta_2 (p)- \theta_1 (p)\rceil,$$ where $\lceil\cdot\rceil$ denotes the next higher integer.
\end{rem} 

\subsection{Partially wrapped Fukaya category}\label{sec:part}
The chain level construction of the partially wrapped $A_\infty$-category of a Liouville manifold is more subtle. The strategy pursued in \cite{GPS1} (based on the original definition of the chain level construction in \cite{AS}) is not to define wrapped Floer chain complexes ``$CW(K,L)$'' but rather to collect all the wrapping in a large $A_\infty$-category $\mathcal O$ and then localize with respect to continuation elements.

The details: Choose a countable set $I$ of exact conical Lagrangian submanifolds, all mutually transverse and at least one for each smooth isotopy class. As in the previous section equip the Lagrangian submanifolds with a primitive, grading and a spin structure. Next choose for each $L \in I$ a countable cofinal sequence as in (\ref{cofinal}) realizing the wrapping
$$L = L^{(0)}\rightsquigarrow L^{(1)}\rightsquigarrow L^{(2)}\rightsquigarrow\ldots$$
of positive exact Lagrangian isotopies for $L$. Define the set $\mathcal O := \IZ_{\geq 0} \times I$ that collects all these cofinal sequences and equip $\mathcal O$ with the partial order given by $L^{(i)} < K^{(j )}$
if and only if $i < j.$ By choosing the $L^{(i)}$ to be in general position, we can assume that any finite totally ordered collection
of Lagrangians in $\mathcal O$ is mutually transverse.
To make $\mathcal O$ into an $A_\infty$-category define the morphism sets as 
$$\mathcal O(L_0,L_1)=
\begin{cases}CF(L_0,L_1), & L_0>L_1\\
\IZ, & L_0=L_1\\
0, & \text{otherwise}.\\\end{cases}
$$ The $A_\infty$-operations 
$$\mathcal O(L_{s-1},L_s)\otimes \cdots \otimes \mathcal O(L_{0},L_1)\to \mathcal O(L_0,L_s)$$
are as usual given by counting honest holomorphic polygons whenever $L_0>\cdots>L_s$ and are set to be zero otherwise.
To define the wrapped Fukaya category $\mathcal W_\ff(W)$ we localize $\mathcal O$ at all continuation elements in $HF(L^{(i+1)}, L^{(i)})$. For the definition of localizing an $A_\infty$-category we refer to \cite{GPS1}.

The key observation is that when going from $\mathcal W_\ff(W)$ to its homology category $H\mathcal W_\ff(W)$, the morphisms $(H\mathcal W_\ff(W))(L, K)$ coincide with the wrapped Floer homology $HW(L, K)$ defined in section \ref{sec:wrappedhomology}.

\subsection{Partially wrapped Floer theory of $\IC$}\label{sec:wrapC}
In this section we discuss the partially wrapped Fukaya categories $\mathcal W_{\ff_n}(\IC)$ when viewing $\IC$ as the Liouville completion of a disk equipped with stops $\ff_n$ given by $n+1$ points on the Liouville boundary $\partial \IC\cong S^1$. These $A_\infty$-categories are a prototype for describing exact triangles (for $n=2$) and higher iterated decompositions (see also \cite{NAD}). We use the shorter notation $\mathcal W_n=\mathcal W_{\ff_n}(\IC)$.
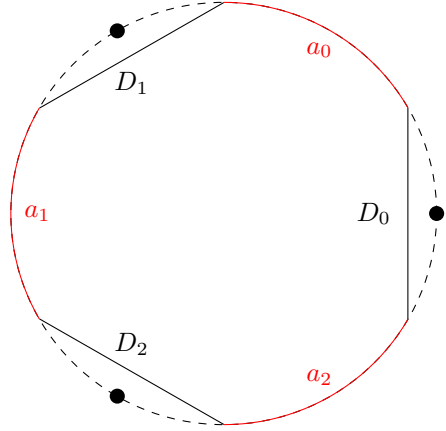
\begin{figure}
    \centering
    \begin{tikzpicture}[scale=0.70]

\draw[dashed] (0,0) ellipse (4 and 4);

\draw (3.46, 2)--node[pos=0.5, left=1mm] {$D_{0}$}(3.46,-2);
\draw[rotate=120] (3.46, 2)--node[pos=0.5, below=1mm] {$D_{1}$}(3.46,-2);
\draw[rotate=-120] (3.46, 2)--node[pos=0.5, above=1mm] {$D_{2}$}(3.46,-2);

\node[red] at (-3.5, 0) {$a_1$};
\node[red] at (1.8, 3.1) {$a_0$};
\node[red] at (1.8, -3.1) {$a_2$};

\draw [red,domain=30:90] plot ({4*cos(\x)}, {4*sin(\x)});
\draw [red,domain=150:210] plot ({4*cos(\x)}, {4*sin(\x)});
\draw [red,domain=270:330] plot ({4*cos(\x)}, {4*sin(\x)});

\node[shape = circle,fill = black, inner sep=2pt] at (-2,-3.46) {};
\node[shape = circle,fill = black, inner sep=2pt] at (4,0) {};
\node[shape = circle,fill = black, inner sep=2pt] at (-2,3.46) {};

\end{tikzpicture}
    \caption{The partially wrapped Fukaya category $\mathcal W_2=\mathcal W_{\ff_2}(\IC)$ as a model for the universal exact triangle. Morphisms $D_j\to D_{j+1}$ are spanned by the Reeb orbit $a_j$ on the Liouville boundary $\partial\IC$.}
    \label{fig:plane}
\end{figure}
\subsubsection{Exact curves} All circles in $\IC$ are non-exact. So the only non-trivial exact curves are non-compact and conical at infinity. We draw them as arcs, i.e. curves in the disk $\IC^{int}$ that connect two points on the Liouville boundary $\partial\IC\setminus \ff_n$. All smooth isotopies of arcs can be realized as exact Lagrangian isotopies according to Remark \ref{ambient}. Such isotopy classes are completely described by the line segments of $\partial \IC\setminus \ff_n$ where the ends of the arc lie. Arcs whose both ends lie in the same segment are zero objects in the partially wrapped Fukaya category $\mathcal W_n$. This can be seen rather directly, as these arcs can be made arbitrarily small and avoid all intersection points with another Lagrangian, or by using Theorem \ref{vanishing}. If the ends of an arc are in two neighbouring boundary segments, i.e. one component cut out by this arc is a topological disk which contains only one stop, then this arc is called a \emph{linking disk}\footnote{ In \cite{GPS2} the notion of a \emph{linking disk} is used in arbitrary dimensions to \emph{link} a locally Legendrian stop. We adopt this notion despite the fact that we are considering a \emph{disk of dimension} $1$.}. Name the linking disks counterclockwise by $D_0,\ldots, D_{n}$ as in Figure \ref{fig:plane}. 
For a morphism from $D_j$ to $D_k$ we first wrap $D_j$ sufficiently in the positive contact direction of the Liouville boundary $\partial \IC$ (but not over a stop) to an arc $D_j^+$. We see that (for $n\geq 2$) the morphism spaces $CW(D_j,D_k)$ are $1$-dimensional for $k=j+1$ and $k=j$, and $0$ otherwise (we use indexing modulo $n+1$). $CW(D_j,D_j)$ is generated by the identity morphism $1_{D_j}$ or the continuation element or by the generator of singular homology $H_0(D_j)$ to list different interpretations. The space $CW(D_j,D_{j+1})$ is generated by the intersection point of $D_j^+$ and $D_{j+1}$ or by the Reeb orbit $a_j$ from $D_j$ to $D_{j+1}$ on the Liouville boundary $\partial \IC$, to list different interpretations again.

\subsubsection{Gradings}\label{sec:gradingC} For the standard complex structure on $\IC$ we use the standard grading given by the squared holomorphic form $dz\otimes dz$ or equivalently by the horizontal foliation given by horizontal lines in $\IC$. Suppose that the linking disks are straight lines in $\IC$ (see Figure \ref{fig:plane} for the case $n=2$). Then a grading of $D_j$ is a constant $d_j$ in $\IR$ and the indices of the morphism $a_j:D_j\to D_{j+1}$ are given by $|a_j|=\lceil d_{j+1}-d_j\rceil$. The only non-trivial higher order $A_\infty$-relations which do not come from strict unitality are $\mu^n(a_{j+n},\ldots, a_j)=1_{D_j}$ for each $j$. This happens only under the constraint $\sum_{j=0}^{n} |a_j|=n-1$. 

\subsubsection{Universal exact triangle} Let $\mathcal W_2=\mathcal W_{\ff_2}(\IC)$ be the disk equipped with $3$ stops on the Liouville boundary. Equip the three linking disks $D_0,D_1,D_2$ with gradings such that the morphisms $a_0,a_1,a_2$ have indices $|a_0|=1$ and $|a_1|=|a_2|=0$. The only non-trivial composition maps are $\mu^3(a_{j+2},a_{j+1}, a_j)=1_{D_j}\in CW(D_j,D_j)$, the identity morphisms. Then the mapping cone of $a_1:D_1\to D_2$ is isomorphic to $D_0$. We write $$D_0\cong [D_1\overset{a_1}\longrightarrow D_2].$$ This exact triangle 
\begin{align}\label{exacttrian}
    D_1\to D_2\to D_0\to D_1[1]
\end{align}
can be used as the prototype of all exact triangles:
\begin{lem}[\cite{SEI} {\normalfont Proposition 3.8, attributed to Kontsevich}]
Let $\mathcal A$ be an $A_\infty$ category. Then any exact triangle in the homological category $H(\mathcal A)$ can be written as the image of an $A_\infty$-functor $F:\mathcal W_2\to \mathcal A$ mapping the exact triangle (\ref{exacttrian}) to the exact triangle in $H(\mathcal A)$.
\end{lem}
Using the previously defined grading and comparing it to the universal exact triangle \ref{exacttrian}, one concludes that the rotated triangles are also exact and are the cones over the two other morphisms $a_2$ and $a_0$: $$D_1[1]\cong [D_2\overset{a_2}{\xrightarrow{\hspace{5mm}}} D_0]\qquad \text{and} \qquad D_2\cong [D_0[-1]\overset{-a_0[-1]}{\xrightarrow{\hspace{8mm}}} D_1].$$

\subsubsection{Universal iterated cone decomposition}\label{sec:universal}
For the definition of an iterated cone decomposition in a triangulated category and some remarks about grading shifts that we chose in this text you may have a look at the appendix \ref{sec:appendix}.

The universal $n$-step iterated cone lives in $\mathcal W_n=\mathcal W_{\ff_n}(\IC)$. Choose gradings on $D_j$ such that the morphisms $a_j:D_j\to D_{j+1}$ have index
\begin{align}\label{eq:gradingC}
    |a_0|=|a_n|=0\qqet |a_{1}|=\ldots =|a_{n-1}|=1.
\end{align} These grading ensure to get the same formulas for the grading shifts as in the cone decomposition given in \cite{HAUG}. 

\begin{figure}
    \centering
    \begin{tikzpicture}[scale=0.70]

\draw[dashed] (0,0) ellipse (4 and 4);

\draw (3.46, 2)--node[pos=0.5, right=5mm] {$D_{0}=\alpha_1$}(3.46,-2);
\draw[rotate=60] (3.46, 2)--node[pos=0.7, above=7mm] {$D_{1}$}(3.46,-2);
\draw[rotate=120] (3.46, 2)--node[pos=0.3, above=7mm] {$D_{2}$}(3.46,-2);
\draw[rotate=180] (3.46, 2)--node[pos=0.5, left=5mm] {$D_{3}$}(3.46,-2);
\draw[rotate=-120] (3.46, 2)--node[pos=0.7, below=7mm] {$D_{4}$}(3.46,-2);
\draw[rotate=-60] (3.46, 2)--node[pos=0.3, below=7mm] {$D_{5}=\alpha_5$}(3.46,-2);

\node[red] at (0, 4.3) {$a_1$};
\node[red] at (0, -4.3) {$a_4$};
\node[red] at (3.7, 2.3) {$a_0$};
\node[red] at (-3.7, 2.3) {$a_2$};
\node[red] at (3.7, -2.3) {$a_5$};
\node[red] at (-3.7, -2.3) {$a_3$};

\node[shape = circle,fill = black, inner sep=2pt] at (-2,-3.46) {};
\node[shape = circle,fill = black, inner sep=2pt] at (-2,3.46) {};
\node[shape = circle,fill = black, inner sep=2pt] at (2,3.46) {};
\node[shape = circle,fill = black, inner sep=2pt] at (2,-3.46) {};
\node[shape = circle,fill = black, inner sep=2pt] at (4,0) {};
\node[shape = circle,fill = black, inner sep=2pt] at (-4,0) {};

\draw (3.46,-2)--node[pos=0.5, right=1mm] {$\alpha_{2}$}(0,4);
\draw (3.46,-2)--node[pos=0.6, right=1.2mm] {$\alpha_{3}$}(-3.46,2);
\draw (3.46,-2)--node[pos=0.5, above=1mm] {$\alpha_{4}$}(-3.46,-2);

\end{tikzpicture}
    \caption{Some arcs in the partially wrapped Fukaya category $\mathcal W_5=\mathcal W_{\ff_5}(\IC)$. The $\alpha_j$ are the intermediate curves that triangulate the hexagon topologically and are intermediate steps in defining an iterated cone decomposing $D_0$.}
    \label{fig:plane2}
\end{figure}
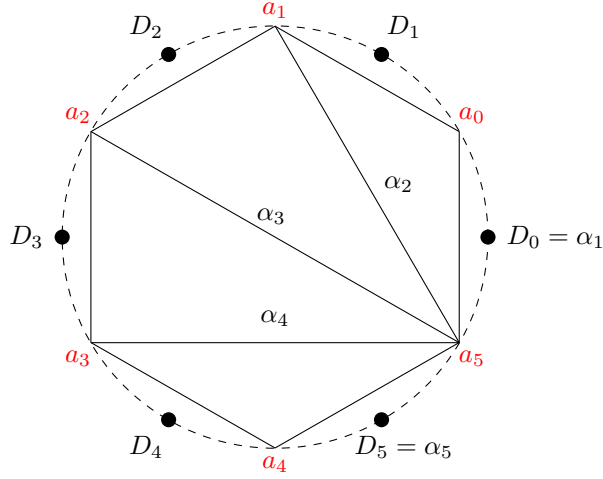

For $j=1,\ldots, n$ define arcs $\alpha_{j}$ with one end between the stop $n$ and the stop $0$, and other end between stop $j-1$ and stop $j$ (see Figure \ref{fig:plane2}).
 In particular, $\alpha_{1}\cong D_{0}$ and $\alpha_n\cong D_{n}$. Equip $\alpha_1$ with the same grading as $D_{0}$, $\alpha_n$ with the same grading as $D_n$. We can grade the other arcs $\alpha_j$ increasing in $j$ compatible with the choice  of the gradings for $\alpha_1$ and $\alpha_n$.

Examining the objects, morphisms and gradings, we get exact triangles
$$\alpha_{j-1}\cong[D_{j-1}[-1]\to \alpha_{j}].$$
Writing this iteratively yields a cone decomposition
$$D_{0}\cong \alpha_1\cong\big[D_1[-1]\to\big[D_2[-1]\to \cdots \to \big[D_{n-2}[-1]\to \big[D_{n-1}[-1]\to D_{n}\big]\big]\cdots\big]\big].$$
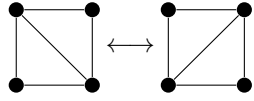
\begin{figure}
    \centering
    \begin{tikzpicture}
\node[shape = circle,fill = black, inner sep=2pt](v1) at (-1,2) {}; 
\node[shape = circle,fill = black, inner sep=2pt] (v2) at (-1,3) {}; 
\node[shape = circle,fill = black, inner sep=2pt] (v4) at (0,3) {};
\node[shape = circle,fill = black, inner sep=2pt] (v3) at (0,2) {};
\draw  (v1) edge node[left] {} (v2);
\draw  (v1) edge node[below] {} (v3);
\draw  (v3) edge node[right] {} (v4);
\draw  (v4) edge node[above] {} (v2);
\draw  (v2) edge (v3);

\node[shape = circle,fill = black, inner sep=2pt](v1) at (1,2) {}; 
\node[shape = circle,fill = black, inner sep=2pt] (v2) at (1,3) {}; 
\node[shape = circle,fill = black, inner sep=2pt] (v4) at (2,3) {};
\node[shape = circle,fill = black, inner sep=2pt] (v3) at (2,2) {};

\node (v6) at (0.5,2.5) {$\longleftrightarrow$};
\draw  (v1) edge node[left] {} (v2);
\draw  (v1) edge node[below] {} (v3);
\draw  (v3) edge node[right] {} (v4);
\draw  (v4) edge node[above] {} (v2);
\draw  (v1) edge (v4);
\end{tikzpicture}\caption{An elementary move (a flip) between two triangulations.}
    \label{fig:elem}
\end{figure}
The geometric interpretation of changing two brackets in the cone decomposition corresponds to a flip (see Figure \ref{fig:elem} and more on this idea in \cite{DK}) in the triangulation of the polygon given by the linking disks. 

Moving all the brackets to the left one obtains
$$D_{0}\cong\big[\big[\big[\cdots\big[D_1[1-n]\to D_2[2-n]\big]\to\cdots \to D_{n-2}[-2]\big]\to D_{n-1}[-1]\big]\to D_{n}\big].$$

Geometrically this corresponds to choosing for $j=1,\ldots, n$ arcs $\beta_j$ with one end between the stop $0$ and the stop $1$, and other end between stop $j$ and stop $j+1$. In particular, $\beta_{1}\cong D_{1}$ and $\beta_n\cong D_{0}$. Equipping $\beta_1$ with the grading of $D_1$, $\beta_n$ with the grading of $D_0$ and $\beta_j$ with an increasing grading between the grading of $D_1$ and $1$ produces exact triangles 
$$\beta_{j}\cong[\beta_{j-1}[-1]\to D_{j}].$$

Finally, note that the linking disks $D_j$ generate the partially wrapped Fukaya category $\mathcal W_n$. 

\subsubsection{Relation to other triangles} The discussed triangles can also be seen as induced by Polterovich surgery of curves: Wrap the ends of the arc $D_{j-1}$ in positive contact direction on $\partial \IC$ such that $\alpha_j$ and $D_{j-1}$ intersect exactly once. Performing surgery of curves (see Figure \ref{fig:surgery}) and forgetting the arc which does not bound a stop (it is a zero object by Theorem \ref{vanishing}) yields the exact triangle. Surgery very close to the Liouville boundary in which one connected component of the resulting surgery is a zero object is called \emph{boundary connected sum}. 

Another point of view is that wrapping an arc over a point, which is a stop, yields an exact triangle. This is the most basic case of wrapping an exact Lagrangian over a locally Legendrian stop (see \cite{GPS2}).

\subsubsection{The big picture}\label{sec:big}
The partially wrapped Fukaya categories $\mathcal W_n$ are a universal model of iterated cone decompositions. This can be compared to Waldhausen's $S_\bullet$-construction. The latter produces a paracyclic object. The articles \cite{NAD,TAN3,DJL} show that this paracyclic system is equivalent to the system produced by the system $\mathcal W_\bullet$. Moreover, \cite{DJL} constructs the higher dimensional geometric realisation $\mathcal W_\bullet^{(d)}$ by partially wrapped Fukaya categories $\mathcal W_{\ff_n}(\mathrm{Sym}^{(d)}(\IC))$ that coherently organise into a paracyclic object and is equivalent to the higher dimensional Waldhausen $S_\bullet^{(d)}$-construction \cite{DYC,POG}. However, the complete description on the symplectic side is not yet written out.

\subsection{Gradings on cobordisms}\label{sec:gradingCob}
Our definition of Lagrangian cobordisms does not a priori define left and right cylindrical ends as in \cite{BC}. We can focus on one end $L_{0}$ to be our preferred end and name the other ends $L_1,\ldots, L_n$ in counterclockwise order which will decompose $L_{0}$ in the wrapped Fukaya category as an iterated cone. We now describe how to get the cone decomposition with correct grading on the ends. Let $\theta_V\colon V\to \IR$ be a grading on the cobordism $V\subset \IC\times W$. Assume that $V$ is connected, otherwise the following discussion can be done separately for each connected component of $V$. An end of a Lagrangian cobordism is the conicalization of a product Lagrangian submanifold $\gamma_j\times L_j$, where the conicalization can start arbitrarily far at infinity. Therefore, we can assume that the complex structure and grading on the product $\IC\times W$ splits as a sum in the relevant region where we will examine intersection points of the cobordism with products $\alpha\times K$ in $\IC\times W$. (For the relevant region compare also to the proof of Proposition \ref{prop:module}.) On the end corresponding to $\gamma_j\times L_j$ the grading $\theta_V$ splits as $\theta_V=\theta_{L_j}+\theta_{\gamma_j}$. Hence
\begin{align}\label{eq:grading}
    &\textit{Given a grading $\theta_V$ on the cobordism $V$, fixing a grading on one end $L_j$ \notag}\\&\textit{induces a grading on $\gamma_j\subset \IC$, and vice versa.}
\end{align}
In \cite{BC,BC2,HAUG,HEN} the ends are chosen to be horizontal. One grades the curves that parametrize the ends by a constant grading $0$ or $1$ (for positive and negative ends, respectively). As we work with ends parametrized by radial curves in $\IC$ we can have a different grading.

Let us record which choices are involved to get a preferred grading on the ends. These are then the correct gradings on the $L_j$ that appear in the cone decomposition of $L_0$ by $L_1,\ldots, L_n$ in Theorem  \ref{thm11}. 
\begin{lem}\label{eq:grading2}
    A grading on the cobordism $V$ and the choice of a preferred end $L_0$ with a preferred grading on $L_0$ induces a preferred grading on all the other ends $L_j$.
\end{lem} 
\begin{proof}
Suppose that $L_0$ is the preferred end with preferred grading $\theta_{L_0}:L_0\to\IR$. This induces a grading on the curve $\gamma_0$ by (\ref{eq:grading}). This in turn fixes a grading on the linking disk $D_0$ by imposing that the intersection of $D_0$ with $\gamma_0$ has index $0$. Given a choice of a grading on $D_0$, we have defined in (\ref{eq:gradingC}), (\ref{eq:gradingCC}) a preferred choice of gradings on the other linking disks $D_j$. These gradings induce gradings on the curves $\gamma_j$ by requiring that the intersection index of $\gamma_j$ and $D_j$ is zero. Finally, again by (\ref{eq:grading}) the gradings on the $\gamma_j$ fixes a grading $\theta_{L_j}$ on the ends $L_j.$
\end{proof}

Alternatively, we can insist on the grading of $\gamma_0$ to be in $[0,1)$, then already the grading on $V$ and the choice of a preferred end $L_0$ induce a preferred grading on all the other ends $L_j$.

\subsection{Cobordisms induce a cone decompositions}\label{sec:conedec}
We work out the idea of Sheridan and Smith \cite{SS} (Section 6.2) to produce exact triangles from a cobordism using stops on the complex plane. In our setting, we need some extra argument how to deal with conicalizations, as cobordisms are not honest products in $\IC\times W$ in the complement of large compact set of $\IC$. 

Our aim is to prove the graded version of Theorem \ref{thm1} in the generalized setting of a stop $\gggg$ on $W$:
\begin{thm}\label{thm11}
Let $(V, \theta_V)$ be a graded connected Lagrangian cobordism in $\IC\times W$. Name one end $(L_0, \theta_{L_0})$ and denote the other ends in counterclockwise order by $L_1, \ldots, L_n$ with induced gradings as discussed in section \ref{sec:gradingCob}. Then there is a cone decomposition
$$L_{0}\cong\big[L_n[-1]\to\big[L_{n-1}[-1]\to\cdots \to \big[L_3[-1]\to \big[L_2[-1]\to L_1\big]\big]\cdots\big]\big]$$ in the derived partially wrapped Fukaya category $D\mathcal W_{\gggg}(W)$.
\end{thm}

The proof will be carried out in section \ref{proof} below.

\subsubsection{Künneth formula of wrapped Fukaya categories}\label{ssec:Kunneth} If $W_1^{int}$ and $W_2^{int}$ are Liouville domains of stopped Liouville manifolds $(W_1, \ff)$ and $(W_2, \gggg)$, respectively, then $W_1^{int}\times W_2^{int}$ is a cornered Liouville domain of $W_1\times W_2$. When we smoothen the corners the Liouville boundary $\partial (W_1\times W_2)$ can be covered by $\partial W_1\times W_2^{int}$, $W_1^{int}\times \partial W_2$ and a part $\partial W_1\times \partial W_2\times \IR$ where the $\IR$-factor parametrizes the smoothening. We define the \emph{product stop} $\mathfrak h$ to be $$\mathfrak h=(\mathfrak f\times \mathrm {sk}(W_2))\cup(\mathrm {sk}(W_1)\times \gggg)\cup (\mathfrak f\times \gggg\times \IR).$$

\begin{ex}(Stabilization)\label{ex:stabs}
Consider a stopped Liouville manifold $(W, \gggg)$ and take its product with the stopped Liouville manifold $(\IC, \pm i\infty)$. The wrapped Fukaya category of the latter Liouville manifold is $\mathcal W_1$ in the notation of section \ref{sec:wrapC}. Recall that a stop should be regarded as lying in the ideal boundary of $\IC$. That is, the stop $\pm i\infty$ should be thought of stopping the dynamics arbitrarily far away on the imaginary axis. The product stop in the Liouville boundary of $\IC\times W$ is then 
$$\mathfrak h=(\{\pm i\infty\}\times \mathrm {sk}(W))\cup(\{0\}\times \gggg)\cup (\{\pm i\infty\}\times \gggg\times \mathbb{R}).$$

When we consider a Lagrangian suspension $V_{L,H}$ of a linear Hamiltonian $H_t$ and an exact conical Lagrangian submanifold $L$ we see that $V_{L,H}$ does not intersect $\{\pm i\infty\}\times \mathrm{sk}(W)$ when thinking about the stop $\{\pm i\infty\}\subset \IC$ as living  far enough outside on the imaginary axis because $V_{L,H}\cap \left(\IC\times W^{int}\right)$ is bounded in the imaginary direction. Moreover, if $L_t$ denotes the Lagrangian isotopy generated by the Hamiltonian flow of $H_t$ and $L_t$ never intersects the stop $\gggg$ in $W$, then $V_{L,H}$ avoids the product stop $\mathfrak h$. Thus $V_{L,H}$ is an object in the wrapped Fukaya category $\mathcal W_\mathfrak h(\IC\times W)$.
\end{ex}

Let us recall the following Künneth formula for partially wrapped Fukaya categories:
\begin{thm}[\cite{GPS2}]\label{thm:Kunn}
For any two Liouville manifolds $W_1, W_2$ equipped with stops $\mathfrak f$ and $\gggg$, respectively, there is a homologically full and faithful $A_\infty$-bifunctor 
$$\mathcal K:\mathcal W_\mathfrak f(W_1)\times \mathcal W_\gggg(W_2)\to \mathcal W_\mathfrak h(W_1\times W_2).$$ It sends a pair $(L,K)$ to the conicalization $L\tilde\times K$ in the product Liouville manifold $W_1\times W_2$.
\end{thm}

\subsubsection{Some algebra}\label{sec:algebra} Let $\mathcal A, \mathcal B, \mathcal C$ be $A_\infty$-categories and suppose $\mathcal F:\mathcal A\times \mathcal B\to\mathcal C$ is an $A_\infty$-bifunctor. Given an object $C$ in $ \mathcal C$ denote $\mathcal Y(C)\in mod_{\mathcal C}$ the image of the Yoneda embedding of $C$ into $mod_{\mathcal C}$. Recall that a $\mathcal C$-module can be identified with an $A_\infty$-functor $\mathcal C^{opp}\to Ch$, where \emph{opp} indicates the opposite category and $Ch$ the $A_\infty$-category of chain complexes. Composing the original functor with the Yoneda module $\mathcal Y(C)$ induces an $A_\infty$-bifunctor 
$\mathcal A^{opp}\times \mathcal B^{opp}\to Ch,$ which in turn induces an $A_\infty$-functor $\mathcal A^{opp}\to mod_{\mathcal B}.$

From another point of view: for every object $A$ in $\mathcal A$ there is an $A_\infty$-functor given by $\mathcal F_A:=\mathcal F(A,-):\mathcal B\to \mathcal C$. Given an object $C$ in $\mathcal C$, the Yoneda module $\mathcal Y(C)$ pulls back to a module $F_A^*\mathcal Y(C)\in mod_{\mathcal B}$ for every $A$. Moreover, the map $A\mapsto F_A^*\mathcal Y(C)$ is a $A_\infty$-functor $\mathcal A^{opp}\to mod_{\mathcal B}.$

\subsubsection{The cone decomposition}\label{proof}
This section contains the proof of Theorem \ref{thm11}. Suppose $(W, \gggg)$ is a stopped Liouville manifold and $V\subset \IC\times W$ is a Lagrangian cobordism with $n+1$ ends $L_0, \ldots, L_{n}$ modeled on conical paths $\gamma_0, \ldots, \gamma_{n}$ in $\IC$. We separate the ends by choosing $n+1$ stops (denoted by $\ff_n$) on a large enough Liouville boundary of $\IC$, one between every two ends (see Figure \ref{fig:cob} for the case $n=2$). Denote by $D_0,\ldots, D_{n}$ the linking disks in $\IC$ of the stops ordered counterclockwise as in section \ref{sec:wrapC} where the $j$-th end of the cobordism lies between the stops with linking disks $D_{j}$ and $D_{j+1}$ (use indexing modulo $n+1$). The Lagrangian cobordism $V$ is an object in the $A_\infty$-category $W_{\mathfrak h}(\IC\times W)$.

As argued in the previous section \ref{sec:algebra}, the Künneth bifunctor $$\mathcal K:\mathcal W_{\mathfrak f_n}(\IC)\times \mathcal W_{\gggg}(W)\to \mathcal W_{\mathfrak h}(\IC\times W)$$
together with the object $V$ induces an $A_\infty$-functor
\begin{align*}
    \mathrm{ev}_V:\mathcal W^{opp}_{\mathfrak f_n}(\IC)\to mod_{\mathcal W_{\gggg}(W)}.
\end{align*}

The image $\mathrm{ev}_V(\alpha)$ of an arc $\alpha$ in $\IC$ is the pullback of the Yoneda module of the cobordism $V$ by the functor induced by $\mathcal K(\alpha, -)$. Geometrically, $\mathrm{ev}_V(\alpha)$ evaluates the cobordism $V$ at those ends which intersect the arc $\alpha$. Let us denote this $\mathcal W_{\gggg}(W)$-module by $M_{\alpha, V}=\mathrm{ev}_V(\alpha)=\mathcal K(\alpha, -)^*\mathcal Y(V)$. The module $M_{\alpha, V}$ evaluated at a test Lagrangian $K$ in $W$ is a chain complex whose homology is the wrapped Floer homology $HW(\alpha\tilde\times K,V).$

The opposite category $\mathcal W_\ff^{opp}(\widetilde{W}, \lambda)$ of a Liouville manifold $\widetilde{W}$ is canonically isomorphic to the wrapped Fukaya category $\mathcal W_\ff(\widetilde W,-\lambda)$. That is, we have the same Lagrangian submanifolds, the $A_\infty$-operations come from holomorphic polygons with reversed orientation, and the wrapping at infinity is in the opposite direction. 

We labelled the linking disks $D_j$ in the counterclockwise direction. Denote $a'_j:D_{j+1}\to D_{j}$ the morphism in $\mathcal W^{opp}_{\mathfrak f_n}(\IC)$ given by the clockwise Reeb chord from $D_{j+1}$ to $D_{j}$. Grading the linking disks $D_j$ analogeously to (\ref{eq:gradingC}) to get indices
\begin{align}\label{eq:gradingCC}
    |a_0'|=|a_n'|=0\qqet |a_1'|=\cdots=|a_{n-1}'|=1
\end{align}we get an iterated cone decomposition:
$$D_{0}\cong\big[D_n[-1]\to\big[D_{n-1}[-1]\to\cdots \to \big[D_2[-1]\to D_1\big]\cdots\big]\big]$$
in $\mathcal W^{opp}_{\mathfrak f_n}(\IC)$. As $\mathrm{ev}_V$ is an $A_\infty$-functor it respects exact triangles and hence cone decompositions, i.e.
\begin{align}\label{eq:cone}
    M_{V, D_0}\cong\big[M_{V,D_{n}}[-1]\to\big[M_{V,D_{n-1}}[-1]\to\cdots \to \big[M_{V,D_{2}}[-1]\to M_{V,D_{1}}\big]\cdots\big]\big]
\end{align} in $mod_{\mathcal W_{\gggg}(W)}$.

To get the desired cone decomposition Theorem $\ref{thm11}$, we are left to prove that, when we apply the functor $\mathrm{ev}_V$ to the linking disks $D_j$, we get exactly the corresponding ends $L_j$. This also proves that the cone decomposition (\ref{eq:cone}) lives in the derived wrapped Fukaya category $D\mathcal W_{\gggg}(W)$.
\begin{prop}\label{prop:module}
There are quasi-isomorphisms $$M_{D_j,V}\cong \mathcal Y(L_j)$$ of $\mathcal W_g(W)$-modules for all $j=0,\dots, n$, where $\mathcal Y(L_j)$ is the image of $L_j$ under the Yoneda embedding $\mathcal W_{\mathfrak g}(W)\to mod_{\mathcal W_{\mathfrak {g}}(W)}$.
\end{prop}

The proof of Proposition \ref{prop:module} follows closely the main argument that proves the cohomological version of the K\"unneth formula (Theorem \ref{thm:Kunn}). For the sake of readability we include a summary of this result and its proof. This is followed by the necessary modifications to our setting that prove Proposition \ref{prop:module}.
\begin{prop}[{\normalfont Cohomological K\"unneth formula \cite{GPS2}}]\label{prop:Kunn}
Given stopped Liouville manifolds $(\IC,\ff)$, $(W,\gggg)$ and Lagrangian submanifolds $\gamma\subset \IC$, $L\subset W$, there is a chain level map inducing isomorphisms 
$$HW(D,\gamma) \otimes HW(K,L)\to HW(D\tilde\times K, \gamma\tilde\times L)$$ for all Lagrangian submanifolds $D\subset \IC,K\subset W$.
\end{prop}

\begin{proof}
The proof is a summary of the sections 6.4-6.5 in \cite{GPS2}.

Start with generic compatible almost complex structures $j_\IC$ on $\IC$ and $J_W$ on $W$ which are both of contact type. Then we have tautological isomorphisms of chain complexes 
\begin{equation}\label{eq111}
    CF(D, \gamma;j_\IC)\otimes CF(K,L;J_W)=CF(D\times K, \gamma\times L;j_\IC\oplus J_W)
\end{equation} as holomorphic maps project to holomorphic maps in each factor of $\IC\times W$.

There is an isomorphism 
\begin{equation}\label{eq222}
CF(D, \gamma;j_\IC)\otimes CF(K,L;J_W)\cong CF(D\tilde\times K, \gamma\tilde\times L_j;J)
\end{equation}
whenever we modify the two involved product Lagrangian submanifolds to be conical in $\IC\times W$ and deform $j_\IC\oplus J_W$ outside a large enough compact set to another almost complex structure $J$, one of contact type. Holomorphic strips can not travel between a small compact set (that in particular contains all intersections of $\gamma_j\times L_j$ and $D_j\times K$) and the modified region due to energy bounds if the modified region is far enough away at infinity. For more details, see the monotonicity and bounded geometry arguments from Proposition 3.19 in \cite{GPS1}.

As the involved choices of $J$ are contractible, standard continuation arguments in Floer theory show that there is an isomorphism $$CF(D, \gamma)\otimes CF(K,L) \cong CF(D\tilde\times K, \gamma\tilde\times L)$$ in the homotopy category $H^0\Ch$ of chain complexes.

Moreover, this isomorphism respects multiplication and hence is compatible with taking the direct limit to wrap the first arguments in $CF(-,-)$. This results in a map
$$HW(D, \gamma)\otimes HW(K,L) \to HW(D\tilde\times K, \gamma\tilde\times L).$$
Actually this is an isomorphism because if we conicalize the product of a cofinal wrapping for $D$ in $\IC$ and a cofinal wrapping of $K$ we get a cofinal wrapping for $D\tilde\times K$ in $\IC\times W$.
\end{proof}

Proposition \ref{prop:module} can now be proved by combining the previous argument with the main argument of Biran and Cornea \cite{BC,BC2} of projecting holomorphic disks to unbounded open subsets of $\IC$.

\begin{proof}[of Proposition \ref{prop:module}]
According to Lemma A.2 from \cite{GPS2} it is enough to check the statement $$\mathcal Y(L_j)\cong M_{D_j,V}$$ pointwise, i.e. we check that there is a chain level map inducing isomorphisms 
$$HW(K,L_j)\to HW(D_j\tilde\times K, V)$$ for all test Lagrangian submanifolds $K\subset W$. Geometrically, if we take the product of a Lagrangian $K$ with the linking disk $D_j$ we only see the $j$-th end $L_j$ of the cobordism $V$.

The definition of a Lagrangian cobordism $V$ provides Liouville domains $\IC^{int}\subset \IC$ and $W^{int}\subset W$ such that the projection of $V\cap \IC\times W^{int}$ intersects $\partial\IC\times W^{int}$ only in (the conicalization of) product Lagrangian submanifolds (see Figure \ref{fig:cob} for the projection to $\IC$).  

The $j$-th end of the cobordism $V$ is modelled on the ray $\gamma_j$ in $\IC$. By wrapping $D_j$ enough in the Reeb direction of $\lambda_\IC$ at infinity we can ensure that $D_j$ and $\gamma_j$ have exactly one intersection point. Equation (\ref{eq111}) then writes as
$$CF(K,L_j;J_W)\cong CF(D_j, \gamma_j;j_\IC)\otimes CF(K,L_j;J_W)=CF(D_j\times K, \gamma\times L_j;j_\IC\oplus J_W).$$
By subsection \ref{sec:gradingCob} the intersection of $D_j$ and $\gamma_j$ does not introduce a grading shift.

Note that if $D_j$ does intersect the projection of $V\cap (\IC\times W^{int})$ to $\IC$ only once and over the curve $\gamma_j$ then in general  $D_j\times K$ and $V$ does not have to intersect only once and over $\gamma_j$. However, we can ensure that the conicalization $D_j\tilde\times K$ and $V$ only intersect once and over $\gamma_j$.

Let $B\subset \IC$ be the open subset as defined in Figure \ref{fig:cob} and denote $\overline{B}=B\times W^{int}$. More explicitly, let $z_j$ be the unique intersection point of $\gamma_j$ and  $D_j$. Let $z_j\in U\subset \IC$ be a small neighbourhood. Then choose an unbounded open set $$B\subset (\IC\setminus (D_j\cup\pi_\IC(V\cap\IC\times W^{int})))\cup U$$ such that $B$ touches each of the four connected component of $\IC$ that the curves $\gamma_j$ and $D_j$ define at the intersection point $z_j$. 
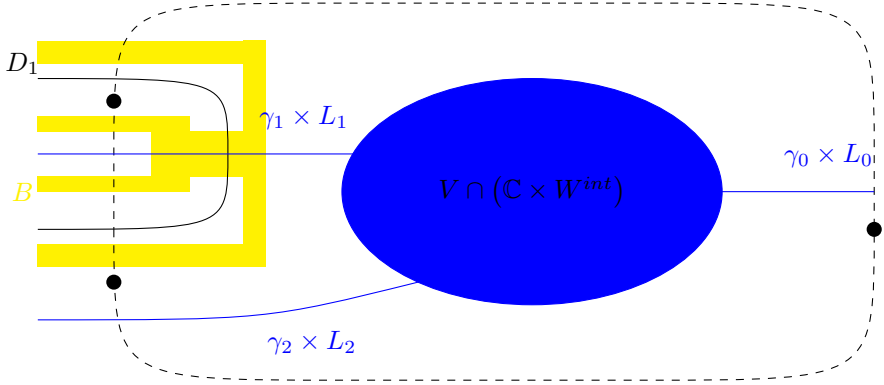
\begin{figure}
    \centering
    \begin{tikzpicture}
\path [fill=yellow] (0.5,2.7) rectangle (2,3.3);
\path [fill=yellow] (-1,4.2) rectangle (2,4.5);
\path [fill=yellow] (-1,1.8) rectangle (2,1.5);
\path [yellow, fill] (1.7,1.5) rectangle (2,4.5);
\path [fill=yellow] (-1,2.7) rectangle (1.,2.5);
\path [fill=yellow] (-1,3.3) rectangle (1.,3.5);

\draw[dashed]   (0,2) .. controls (0,0) and (0,0) .. (3,0) -- (7,0) .. controls (10, 0) and (10, 0) .. (10, 2) -- (10,3) .. controls (10,5) and (10,5) .. (7,5) -- (3,5) .. controls (0,5) and (0, 5) .. (0,3) -- cycle ;

\node[shape = circle,fill = black, inner sep=2pt] at (0,1.3) {};
\node[shape = circle,fill = black, inner sep=2pt] at (0,3.7) {};
\node[shape = circle,fill = black, inner sep=2pt] at (10,2) {};

\draw (1.5,3) ..  controls (1.5,4) and (1.5,4) .. (-1,4);
\draw (1.5,3) ..  controls (1.5,2) and (1.5,2) .. (-1,2);

\node[text=yellow] at (-1.2, 2.5){$B$};
\node at (-1.2, 4.2) {$D_1$};

\draw [blue] (-1,0.8)  .. node[pos=0.7, below=2mm] {$\gamma_2\times L_2$} controls (2,0.8) and (2,0.8) .. (4,1.3);
\draw [blue](-1,3) -- node[pos=0.7, above=2mm] {$\gamma_1\times L_1$}  (4,3);
\draw [blue](7,2.5) -- node[pos=0.8, above=2mm] {$\gamma_0\times L_0$} (10,2.5);
\draw[fill=blue, draw=blue] (5.5,2.5) ellipse (2.5 and 1.5);
\node at (5.5, 2.5) {$V\cap \left(\mathbb C\times W^{int}\right)$};

\end{tikzpicture}
    \caption{In blue, the projection of a Lagrangian cobordism and in yellow, an example of an open set $B$ used to prove that holomorphic curves connecting points in the intersection of $V$ and $D_2\times K$ stay in the same fiber over $\IC$.}
    \label{fig:cob}
\end{figure}
Let $J$ be an almost complex structure on $\IC\times W$ which is equal to $j_\IC\oplus J_W$ in $\overline B$ and conical in the product $\IC\times W$ outside a small neighbourhood around $\overline B.$

A standard argument as in \cite{BC} (Lemma 4.2.1) and \cite{BC2} (Proposition 3.3.1) proves that all holomorphic disks lie in the fiber over the point $z_j\in \IC$. Thus we get isomorphic chain complexes as in equation (\ref{eq222}): $$CF(D_j\times K, \gamma_j\times L_j;j_\IC\oplus J_W)\cong CF(D_j\tilde\times K, \gamma_j\tilde\times L_j;J)\cong CF(D_j\tilde\times K, V;J).$$ Indeed, after replacing $\gamma_j\tilde\times L_j$ with $V$, still all holomorphic disks need to lie in the fiber over $z_j$. 

By the same reasoning as before we can change the almost complex structure $J$ far out at infinity and make $J$ conical everywhere in $\IC\times W$. As the involved choices of $J$ (and $B$) are contractible, standard continuation arguments in  Floer theory show that there is an isomorphism 
$$CF(K,L_j)\cong CF(D_j\tilde \times K, V)$$ in the homotopy category $H^0\Ch$ of chain complexes.
Now wrapping the first argument in $CF(-,-)$ yields
$$HW(K,L_j)\to HW(D_j\tilde\times K, V)$$
which is an isomorphism arguing exactly as in the proof of Proposition  \ref{prop:Kunn}.
\end{proof}

\subsection{Category of cobordisms}\label{sec:functorial} Biran and Cornea \cite{BC2} introduced a functor that describes how cobordisms induce cone decompositions in a functorial way. Below  we  recall this functor in our framework. We refer to \cite{BC2} and the appendix \ref{sec:appendix} for details about cone decompositions (like composing cone decompositions).

The objects of the \emph{cobordism category} $Cob_\ff(W)$  of a stopped Liouville manifold $(W, \lambda, \ff)$ are  tuples $(L_1,\ldots, L_n)$ of exact conical Lagrangian submanifolds $L_j$ that avoid $\ff$ and are equipped with some extra structure. A morphism $[V]:(L')\to (L_1, \ldots, L_n)$ from a single-elemented tuple $(L')$ is a planar isotopy class of a Lagrangian cobordism $V$ that has ends $L', L_1,\ldots, L_n$ (in counterclockwise order) and avoids $\ff$ as described in the definition of Lagrangian cobordisms. This implicitly contains the choice of a special end $L'$ which will be decomposed by the other ends $L_1, \ldots, L_n$.
A morphism $[V]:(L_1', \ldots, L'_{m})\to (L_1, \ldots, L_n)$ from a general tuple is an ordered tuple $([V_1],\ldots, [V_{m}])$ of morphisms $[V_j]:(L_j')\to (L^j_1, \ldots, L^j_{n_j})$ where 
$$(L_1, \ldots, L_n)=(L_1^1,\ldots, L_{n_1}^1, L_1^2, \ldots,L_{n_2}^2, \ldots, L_1^{m}, \ldots,L_{n_m}^{m})$$ as ordered tuples.
By the discussion in section \ref{sec:planar} we can assume that the cobordisms $V_1,\ldots, V_m$ live in $\IC\times W$ with the same Liouville structure on $\IC$, are disjoint and if projected to $\IC$ ordered from lower half plane to upper half plane in $\IC$. Hence, $V$ is a cobordism on its own with ends $$L'_{m}, \ldots, L'_1, L_1, \ldots,L_n.$$ 
The composition of two morphisms $[V]$ and $[V']$ is given by gluing ends, see section \ref{sec:concat} about concatenating cobordisms. Notice that in compositions we always glue only one end of a Lagrangian cobordism inside the total Lagrangian cobordism $V'$ to an end in a Lagrangian cobordism in $V$, so composition is well-defined. As explained in \cite{BC2} we can make this category into a strict monoidal one by quotiening by void Lagrangians and void cobordisms.

Next we recall the category $T^S\mathcal A$ of \emph{triangle resolutions} of a triangulated category $\mathcal A$. An object in $T^S\mathcal A$ is a tuple $(x_1,\ldots, x_n)$ of objects $x_j$ in $\mathcal A$. 
A morphism $(x)\to (y_1,\ldots y_n)$ from a single-elemented tuple $(x)$ is a cone decompositions $$\eta=[y_n\to\cdots\to y_1]$$
plus an isomorphism $x\cong \eta$ modulo equivalences of cone decompositions (see \cite{BC2}).
Note that we do not write the brackets and gradings $$x\cong \big [y_n[-1]\to\cdots \big[y_3[-1]\to \big[y_2[-1]\to y_1\big]\big]\cdots]$$ in this section for notational convenience. See \ref{sec:appendix} for more on brackets and grading shifts.

A morphism $(x_1, \ldots x_{m})\to (y_1,\ldots y_n)$ from a general tuple is an ordered tuple of morphisms $(x_j)\to (y^j_1, \ldots, y^j_{n_j})$ where 
$$(y_1, \ldots, y_n)=(y_1^1,\ldots, y_{n_1}^1, y_1^2, \ldots,y_{n_2}^2, \ldots, y_1^{m}, \ldots,y_{n_m}^{m})$$ as ordered tuples.

\begin{prop}
There exists a functor 
$$\mathcal G: Cob_\ff(W)\to T^SD\mathcal W_\ff(W)$$ that maps tuples $(L_1,\ldots, L_n)$ of Lagrangians to  $(L_1,\ldots, L_n)$. A simple morphisms coming from a Lagrangian cobordism $V$ is mapped to the morphism $\mathcal[V]:(L_0)\to (L_1, \ldots, L_n)$ in $T^SD\mathcal W_\ff(W)$ given by the cone decomposition that is induced by $V$ and was described in Theorem \ref{thm11}.
\end{prop}

\begin{proof}
The details of the proof can be found in \cite{BC2}. Gluing two cobordisms basically boils down to gluing the two Fukaya categories
$$\mathcal W_{\ff_n}(\IC) \qqet \mathcal W_{\ff_m}(\IC)$$ to  get $\mathcal W_{\ff_{n+m-1}}(\IC)$ for any $m,n\geq 0$.

Let us illustrate how  $\mathcal G$ is compatible with gluing two cobordisms at one end.
Let $V':(L_{0})\to (L_1, \ldots, L_n)$ and $V'':(L_{j_0})\to (K_1, \ldots, K_{m})$ be cobordisms for some $j_0\in\{1, \ldots n\}$ and let $V$ be the cobordism obtained by gluing (see section \ref{sec:concat}) along the ends $L_{j_0}$.
The counterclockwise labelled ends of $V$ are $$L_1, \ldots, L_{j_0-1}, K_1, \ldots, K_{m}, L_{0}, L_{j_0+1}\ldots L_n,$$ see Figure \ref{fig:cobcon} on \cpageref{fig:cobcon}.

For the cone decompositions of $V'$ and $V''$, we separate the ends by exactly one stop between each end. However, we can put more stops in between without changing the cone decomposition: Forget $V''$ in Figure \ref{fig:cobcon} for a moment and just consider $V'$. Let $D_{j_0}$ in $\IC$ be an arc that intersects $V$ (and hence also $V'$) only in one end, the end $L_{j_0}$. The arc $D_{j_0}$ is not a linking disk by our standard definition, because on its left side it might contain more than one stop coming from the cobordism $V''$. However, we get the same cone decomposition $\mathcal G(V')$ as if we had only one stop.

The same line of arguments applies to the cone decomposition of $L_{j_0}$ by the cobordism $V''$. 

The cone decomposition of $L_0$ by the glued cobordism $V$ is 
$$L_{0}\cong[L_n\to\cdots \to L_{j_0+1}\to K_{m}\to \cdots\to K_1\to L_{j_0-1}\to\cdots\to L_1]$$
which is the concatenation in $T^S\mathcal W_\ff(W)$ of the cone decompositions 
$$L_{0}\cong[L_n\to\cdots\to L_1]\qqet L_{j_0}\cong [K_m\to \cdots\to K_1]$$
that are induced by $V'$ and $V''$. The equivalence is induced by isotoping $D_{j_0}$ in $\IC$ to intersect the other ends $K_1,\ldots, K_{m}$ of $V''$ instead of only intersecting the end $L_{j_0}$.
\end{proof}

\section{Cobordism Groups}\label{sec5}
In this section we prove Theorem \ref{thm22}. Let us start by recalling the definition of the Lagrangian cobordism group.
\begin{defi}\label{def:cobgroup}
The \emph{Lagrangian cobordism group} $\Omega_\ff (W)$ is the abelian group generated by all exact conical Lagrangian submanifolds (equipped with a primitive, a Pin structure, a grading) modulo the relations $L_{0}=L_1+\cdots +L_n$ whenever there is a graded cobordism $V$ with graded ends $L_0,L_1\ldots, L_{n}$ as in Theorem \ref{thm11}, where $V$ induces the gradings on the $L_j$ as explained in section \ref{sec:gradingCob}.
\end{defi}

As an immediate consequence of Theorem \ref{thm11} we have
\begin{cor}
There is a surjective group homomorphism
\begin{align*}
    \Omega_\ff (W ) \to K_0 (\mathcal W_\ff(W))
\end{align*}
induced from $L\mapsto L$.
\end{cor}

The relation in Definition \ref{def:cobordism} may not look invariant under choosing another end of the cobordism as $L_0$ and rotating the enumeration of the other ends. However, changing the preferred end leads to different gradings the other ends, eventually resulting in the same relation (see subsection \ref{sec:gradingCob} for the definition of the induced gradings and the appendix \ref{sec:appendix} for the gradings after rotations). 

A related group is the following. The group $\mathcal L_\ff(W)$ is generated by Lagrangian submanifolds of $W$ (equipped with extra structures)  modulo the relations 
\begin{itemize}
    \item $L_0=L_1$ whenever there is an isotopy of conical exact Lagrangian submanifolds $L_t$ connecting $L_0$ to $L_1$,
    \item $L=L_2+L_1$ whenever $L_1$ and $L_2$ intersect transversely in exactly one point (with intersection index $1$) and $L$ is the Polterovich surgery of $L_1$ and $L_2$,
    \item $L=0$ whenever $L$ can be flowed away from any compact set by the Liouville flow (without flowing over the stop $\ff)$. In other words $L$ does not intersect the generalized skeleton (see subsection \ref{sec:null-cob}).
\end{itemize}
As these three relations can be represented by cobordisms as described in sections \ref{sec:susp}, \ref{sec:surgery} and \ref{sec:null-cob} we get two surjective group morphisms
\begin{align}\label{Groth}
    \mathcal L_\ff(W)\to\Omega_\ff (W ) \to K_0 (\mathcal W_\ff(W)).
\end{align}
We will see in section \ref{sec:surfaces} that for surfaces with boundary all of these maps are isomorphisms. One may conjecture that this continues to hold in the more general case of Weinstein manifolds.

\subsection{Cobordism group of the cylinder}\label{sec:cylinder}
As a warm-up let us calculate the cobordism group of the cylinder $T^*S^1$. We claim that all the maps in (\ref{Groth}) are isomorphisms. 

As the cylinder is a cotangent bundle its Grothendieck group is isomorphic to $\mathbb Z$ and is generated by a fiber (see \cite{ABO3}). On the other hand, there are four smooth isotopy classes of embedded non-oriented curves and arcs. They can be represented by the zero section, a fiber, an arc $\gamma_+$ with both ends at $+\infty$ in the cotangent coordinate and an arc $\gamma_-$ with both ends at $-\infty$. In $T^*S^1$ any closed exact Lagrangian submanifold is in fact exact Lagrangian isotopic to the zero section. For arcs in a surface, smooth isotopies and exact Lagrangian isotopies are the same as $H^1(\IR)=0.$ This classifies conical exact Lagrangian submanifolds up to exact Lagrangian isotopy. 

The arcs $\gamma_-$ and $\gamma_+$ are both null-cobordant as Proposition \ref{prop:null} shows. Moreover, when we perform a surgery between the zero section and a fiber we get a Dehn twist which is Hamiltonianly isotopic through a linear Hamiltonian to the fiber again. Consequently, the cobordism group is at most one dimensional, generated by a fiber. As $K_0 (\mathcal W(T^*S^1 ))\cong \mathbb Z$ is generated by the fiber the maps in (\ref{Groth}) must all be isomorphisms.

\subsection{Surfaces with boundary}\label{sec:surfaces}
Before looking at surfaces with boundary, let us state the known results that relate the Grothendiek group of the Fukaya category of closed surfaces to its cobordism group.

The cobordism group of the torus $\mathbb{T}^2$ has been calculated in \cite{HAUG} (a further discussion of higher genus surfaces can be found in \cite{PER}). When the Lagrangian submanifolds are non-contractible embedded circles in $\mathbb T^2$ then 
$$\Omega(\mathbb {T}^2)\cong K_0( \mathrm{Fuk}(\Sigma))\cong H_1(\mathbb T^2)\oplus \IR/\IZ.$$ The latter identification is not canonical. The extra term $\IR/\IZ$ roughly arises as flux by fixing a circle and measure the area swept out when isotoping to another circle. This term will not arise in the case of surfaces with boundary as we will restrict to exact curves. Exact Lagrangian isotopies of circles do always sweep out zero area. 

Suppose now that $\Sigma$ is a closed Riemann surface with finitely many open disks removed, called a Riemann surface of finite type with boundary $\partial \Sigma$. Moreover, let $\ff$ be a finite subset of $\partial \Sigma$, where arcs are stopped. 
The wrapped Fukaya category of punctured surfaces has been calculated in \cite{ABETC,BOC} and was generalized in \cite{HKK} to include stops on the boundary. 
Open Riemann surfaces are Stein (see \cite{CE}), so we can equip each surface $\Sigma$ with a Weinstein structure and speak about exact Lagrangian submanifolds. Recall that the space of Liouville structures on an oriented punctured surface is convex. By \cite{CE} any two complete Liouville structures are exact symplectomorphic. So it is enough to calculate the partially wrapped Fukaya category with respect to one particular complete Liouville structure. 

A \emph{full arc system} $A$ is a finite set of pairwise disjoint, non-isotopic arcs in $\Sigma$ such that $A$ cuts $\Sigma$ into topological disks that each contain at most one point of the stop $\ff$. In \cite{HKK} it is described how to deduce from the contractibility of the arc complex (see \cite{HAR1,HAR2,HAT}) that any full arc system generates the partially wrapped Fukaya category $\mathcal W_\ff(\Sigma)$ and that there is an isomorphism 
$$K_0(\mathcal W_\ff(\Sigma))\cong H_1(\Sigma, \partial \Sigma\setminus \ff).$$ 

Following the argument in \cite{HAT} we can show that the relations that define  $\mathcal L_\ff(\Sigma)$ are enough to recover the relations in the Grothendieck group of the partially wrapped Fukaya category $\mathcal W_\ff(\Sigma)$. In particular, the cobordism group and the Grothendieck group agree.

\begin{thm}\label{thm22}
We have isomorphisms of groups
$$\mathcal L_\ff(\Sigma)\cong\Omega_\ff(\Sigma)\cong K_0(\mathcal W_\ff(\Sigma))\cong  H_1(\Sigma, \partial \Sigma\setminus \ff).$$
sending an arc to its class in relative homology.
Moreover, any full arc system $A$ in $\Sigma$ generates each of these groups.
\end{thm}

\begin{proof}
The proof is purely topological and combinatorial. There is a well-defined surjective map of groups 
\begin{align}\label{eq:homo}
    \mathcal L_\ff(\Sigma)\to H_1(\Sigma, \partial \Sigma\setminus \ff)
\end{align} taking an arc or an exact closed curve to its homology class. Arcs that can be flowed away from any compact set by the Liouville flow are isotopic to an arc in $\partial \Sigma\setminus \ff$ and are therefore zero in relative homology. Note also that surgery of two Lagrangian submanifolds intersecting transversely at one point is equivalent to their sum in the first homology group.

Fix a full arc system $A$ on $(\Sigma, \ff)$. As the complement of $A$ in $\Sigma$ consists only of disks, $A$ can be shown to generate $H_1(\Sigma, \partial \Sigma\setminus \ff)$ by considering cellular homology and relations $L_1+\cdots +L_n=0$ for any disk which is bounded by $L_1,\ldots, L_n$ (and some boundary arcs in $\partial \Sigma\setminus \ff$). Furthermore, all the relations in cellular homology can be realized by surgeries close to the boundary $\partial \Sigma$ (see Figure \ref{fig:concat}).
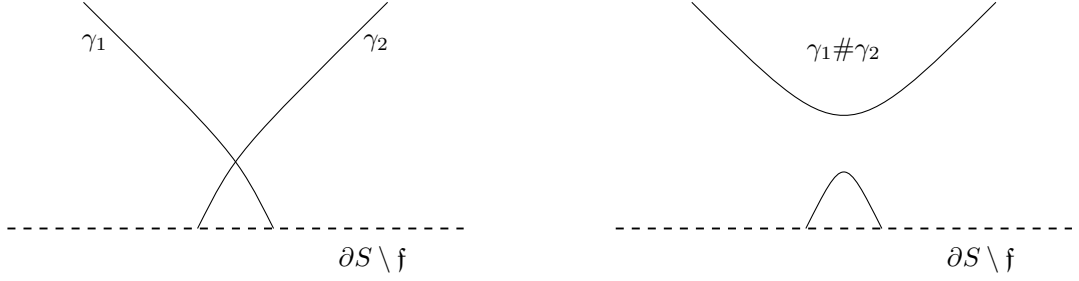
\begin{figure}
    \centering
    \begin{tikzpicture}
\draw[line width=0.25mm, dashed ] (-3,0) -- node[pos=0.8, below=1mm] {$\partial S\setminus \mathfrak f$}(3,0);
\draw (-2,3) .. node[pos=0.1, left=1mm] {$\gamma_1$} controls (0,1) and (0,1) .. (0.5,0);
\draw (2,3) .. node[pos=0.1, right=1mm] {$\gamma_2$} controls (0,1) and (0,1) .. (-0.5,0);

\draw[line width=0.25mm, dashed ] (5,0) -- node[pos=0.8, below=1mm] {$\partial S\setminus \mathfrak f$}(11,0);
\draw (6,3) .. node[pos=0.5, above=6mm] {$\gamma_1\#\gamma_2$} controls (8,1) and (8,1) .. (10,3);
\draw (7.5,0) .. node[pos=0.1, right=1mm] {} controls (8,1) and (8,1) .. (8.5,0);
\end{tikzpicture}
    \caption{Surgery of two arcs $\gamma_1, \gamma_2$ that intersect transversely and only once and with one connected component of $\Sigma\setminus (\gamma_1\cup \gamma_2)$ being a topological disk that does not contain a stop. This surgery can be viewed as a concatenation of arcs because the resulting arc that bounds the disk is a zero object. This is also known as \emph{boundary connected sum} of arcs.}
    \label{fig:concat}
\end{figure}

The theorem follows as soon as we prove that a full arc system $A$ in $\Sigma$ also generates $\mathcal L_\ff(\Sigma)$. 
First, let us show that any embedded arc can be obtained by using the operations in the construction of $\mathcal L_\ff(\Sigma)$. By definition of a full arc system, for any point in the stop $\ff$ we can do surgery of the arcs, that bound the disk which contains the point, to get a linking disk for this point. So we already constructed all linking disks.
This argument reduces the analysis to a surface with boundary but without any stops.

The polygons in the complement of the arc system can be divided into triangles by adding more non-isotopic arcs. These extra arcs are concatenations of arcs in the arc system and can be realized as surgeries near the boundary of $\partial \Sigma$ (see Figure \ref{fig:concat}). We end up with a triangulation $\Delta$ of $\Sigma$ in the general sense of degenerate triangles, i.e. that vertices of triangles can coincide, as can pairs of edges. Let now $\gamma$ be any arc. Cut $\Sigma$ along $\gamma$ and let $A'$ be a full arc system for the surface $\Sigma\setminus \gamma$. We can assume that $A'$ has no ends lying on $\gamma$ by sliding the ends of arcs in $A'$ that lie on $\gamma$ to the ends of $\gamma$. Then $A'\cup \{\gamma\}$ is a full arc system for $\Sigma$ that contains $\gamma$. We can complete this full arc system to a triangulation $\Delta'$ which contains $\gamma$ as argued before. 

According to \cite{HAT} any two triangulations $\Delta, \Delta'$ of a surface can be related by finitely many elementary moves (see Figure \ref{fig:elem} on  \cpageref{fig:elem}).
An elementary move can be obtained with our operations. Forget the diagonal arc and do the boundary connected sum of the two arcs that yields the other diagonal arc. Consequently, any arc $\gamma$ can be obtained by starting with an arc system $A$ and by performing finitely many surgeries of two transversely intersecting arcs that only intersect once.

We address now closed exact curves in the interior of $\Sigma$. We only need to find an arc that intersects the closed curve exactly once and transversely. Then after performing surgery we get another arc (possibly isotopic to the first arc). This will prove that we can present the closed circle by means of two arcs in $\mathcal L_\ff(\Sigma)$.

An exact closed curve $\gamma$ cannot separate a surface into two parts where one part is compact and has only $\gamma$ as its boundary (otherwise this part would have zero area by Stokes theorem). If $\gamma$ separates $\Sigma$ into two connected components, then find an arc going from the boundary in one connected component to $\gamma$ and an arc going from $\gamma$ to the other connected component. Its concatenation intersects $\gamma$ exactly once.
In case the closed curve is non-separating, then by starting at one boundary component, and crossing the closed curve transversely one can go back to a boundary component.

This proves that (\ref{eq:homo}) is an isomorphism.
\end{proof}

\begin{rem}
In section \ref{sec:surfaces} and hence Theorem \ref{thm22} we assumed that the stop $\ff$ is a finite subset of $\partial \Sigma$. The cases where the stop $\ff$ contains entire boundary circles is not covered by Theorem \ref{thm22} but interesting as well. For example, the case $\ff=\partial \Sigma$ is equivalent to allow only compact exact curves in the Fukaya category, no arcs. This case is not covered by Theorem \ref{thm22} and needs to be analyzed more carefully.

In the case of the cylinder $\Sigma=T^*S^1$ with stop $\ff=\partial \Sigma$ the statement $$\mathcal L_\ff(\Sigma)\cong\Omega_\ff(\Sigma)\cong K_0(\mathcal W_\ff(\Sigma))\cong  H_1(\Sigma, \partial \Sigma\setminus \ff)\cong H_1(\Sigma)\cong \IZ$$
is still true with generator the zero section.

However, already for the pair of pants $\Sigma$ the three isotopy classes coming from the three boundary components $\partial \Sigma$ contain each a representative of an exact circle. But these three circles can never intersect pairwise in only one transverse intersection, we do not get a relation by surgery. Hence 

$$\IZ^3\cong\mathcal L_\ff(\Sigma)\not\cong H_1(\Sigma)\cong \IZ^2.$$
\end{rem}

\vspace{4mm}
Sometimes it is helpful to have a full arc system with minimal number of arcs. That is, we would like to have an arc system that divides a stopped surface into a minimal number of topological disks. 

\begin{cor}\label{prop:minimal} For any connected stopped surface $(\Sigma, \ff)$ with $\ff\neq\emptyset$ there is a full arc system $A$ such that cutting the surface along the arcs in $A$ results in $|\ff|$ topological disks with each disk containing exactly one point in $\ff$. In this case 
$$\Omega_\ff(\Sigma)\cong\IZ^{|A|}.$$
We call such a full arc system $A$ \emph{minimal}.
\end{cor}

To construct a minimal arc system, note that any surface with boundary can be obtained by identifying sides starting from one polygon. Add one arc for each two identified sides plus add linking disks for all but one stop.

Using these minimal full arc systems we can conclude the following:

\begin{cor}\label{cor:stopremoval}(Stop removal)
Suppose $(\Sigma, \ff)$ is a stopped surface with stop $\ff\neq \emptyset$. Let $p\in \partial \Sigma\setminus \ff$. Then $$\Omega_{\ff\cup\{p\}}(\Sigma)\cong\Omega_\ff(\Sigma)\oplus \IZ.$$
On the other hand if $\ff=\emptyset$ then adding a point $p\in \partial \Sigma$ to $\ff$ does not change the cobordism group:
$$\Omega_{\{p\}}(\Sigma)\cong \Omega(\Sigma).$$
\end{cor}

The second assertion follows because any full arc system for $\Sigma$ stays a full arc system for $\Sigma$ stopped at a point which is not lying on any of the arcs in the system.

\subsection{Gluing surfaces with boundary}
Wrapped Fukaya categories are expected to be well-behaved under gluing. The major known result is the homotopy pushout formula of \cite{GPS2} where a Liouville sector is divided along a hypersurface into two Weinstein sectors. A special case of this is the operation of attaching Weinstein handles. For surfaces, more concrete tools are at hand. In dimension two, gluing of two surfaces can be understood from the point of view of gluing the skeletons. This is entirely topological. The works of \cite{DK,DYC2,SP} illustrate how gluing wrapped Fukaya category in dimension two boils down to gluing graphs along half-infinite edges. This corresponds to attaching a one-handle connecting two points in the boundary of stopped surfaces. We have seen in section \ref{sec:surfaces} that for surfaces with boundary the cobordism group and the Grothendieck group of the wrapped Fukaya category both agree with the first homology group. The statement of the following section are compatible with these isomorphisms. The given arguments are a version of relative cellular homology using disjoint 1-cells, namely the arcs. 

\subsubsection{Linear gluing}\label{sec:lineargluing}
We use the following model for attaching a Weinstein one-handle in dimension two: A one-handle is described by the unit disk $D\subset \IC$ with stop given by two points $\pm i$ on the boundary. Suppose $\Sigma$ is a (possibly disconnected) surface with boundary $\partial \Sigma$ and stop $\ff \subset \partial \Sigma$ containing at least two points $q_\pm$. We can deform the Liouville vector field on each of the surfaces $D$ and $\Sigma$ such that it is parallel to the boundary in a neighbourhood of $\pm i$ and $q_\pm$, respectively. Then we can identify a small interval in $\partial D$ containing $+i$ with a small interval in $\partial \Sigma$ containing $q_+$, and similarly for $-i$ and $q_-$. The intervals in the resulting stopped surface $(\Sigma', \ff'=\ff\setminus \{q_\pm\}$) are called the \emph{cocores} of the one-handle. As they are Hamiltonianly isotopic we speak of one cocore of the one-handle.

The effect of a handle attachment on cobordism groups (and analogously on the Grothendieck group as well as the relative homology) is the following: 
\begin{prop}\label{handle}
Suppose the connected stopped surface $(\Sigma', \ff')$ arises by attaching a one-handle to stopped surface $(\Sigma, \ff)$ at the points $q_1,q_2\in \ff$.

Then if $\ff'\neq \emptyset$ there is an exact sequence
$$0\to \Omega_{\pm i}(D)\to \Omega_\ff(\Sigma)\to \Omega_{\ff'}(\Sigma')\to 0,$$ where the second map comes from the inclusion of stopped surfaces. Denoting by $\gamma$ an arc generating $\Omega_{\pm i}(D)$, then $\gamma$ is mapped to a linking disk $\gamma_-$ of $q_-$ minus a linking disk $\gamma_+$ of $q_+$. The second map sends both of $\gamma_\pm$ to the cocore $\gamma'$ of the handle attachment.

In the case $\ff'=\emptyset$ the inclusion $(\Sigma, \ff)\to (\Sigma',\emptyset )$ induces an isomorphism
$$\Omega_\ff(\Sigma)\to \Omega(\Sigma').$$
\end{prop}

\begin{proof}
We will give a proof by using minimal full arc systems as defined in Corollary \ref{prop:minimal}. By abuse of notation we identify arcs in $(\Sigma, \ff)$ and $(\Sigma', \ff')$ which is possible whenever the arcs do not intersect the cocore $\gamma'$. In particular, to conclude the statement for non-empty stop $\ff'$  we construct a minimal arc systems for $(\Sigma', \ff')$ such that adding one more arc chosen from $\{\gamma_\pm\}$ is a minimal arc systems for $(\Sigma, \ff)$. The stopped surface $(\Sigma'\setminus \gamma', \ff')$ is homeomorphic to $(\Sigma, \ff\setminus\{q_\pm\})$ by definition of the handle attachment.

Suppose first that $\ff'\neq \emptyset$ and $\Sigma$ is connected. In this case $(\Sigma'\setminus \gamma', \ff')$ admits a minimal full arc system  $A'$ by Corollary \ref{prop:minimal}. Note that $A'\cup \{\gamma'\}$ is a minimal full arc system for $(\Sigma', \ff')$ and $A'\cup \{\gamma_\pm\}$ is a minimal full arc system for $(\Sigma, \ff)$. 

Again suppose that $\ff'\neq \emptyset$ but now that $(\Sigma, \ff)$ has two connected components $(\Sigma_-, \ff_-\cup\{q_-\})$ and $(\Sigma_+, \ff_+\cup\{q_+\})$. Suppose $\ff_-, \ff_+$ are both non-empty. Then $(\Sigma_-, \ff_-)$ and $(\Sigma_+, \ff_+)$  both admit minimal full arc systems $A_-$ and $A_+$, respectively. In this case $A_-\cup A_+\cup \{\gamma'\}$ is a minimal full arc system for $(\Sigma', \ff')$ and $A_-\cup A_+\cup \{\gamma_\pm\}$ is a minimal full arc system for $(\Sigma, \ff)$. We are missing the case when one of $\ff_-$, $ \ff_+$ is empty, say $\ff_+=\emptyset$. Then we can find a minimal full arc system $A_-$ for $(\Sigma_-, \ff_-)$ and a minimal full arc system $A_+$ for $(\Sigma_+, \{q_+\})$. Then $A_-\cup A_+$ is a minimal full arc system for $(\Sigma', \ff')$ and $A_-\cup A_+\cup \{\gamma_-\}$ is a minimal full arc system for $(\Sigma, \ff)$.

Finally, we reduce the case $\ff'=\emptyset$ to the previous cases by enlarging the stop $\ff'$ (and hence $\ff$) by one point $\{p\}$. We have an exact sequence $$0\to \Omega_{\pm i}(D)\to \Omega_{\ff\cup\{p\}}(\Sigma)\to \Omega_{\{p\}}(\Sigma')\to 0.$$ But according to Corollary \ref{cor:stopremoval} there are isomorphisms $\Omega_{\ff}(\Sigma)\cong\Omega_{\ff\cup\{p\}}(\Sigma)\oplus \IZ$ and $\Omega_{\{p\}}(\Sigma')\cong\Omega(\Sigma')$ which are compatible with all maps as we only deal with inclusions of stopped surfaces.
\end{proof}

\subsubsection{Circular gluing}
Another way to decompose surfaces with boundary are pair of pants decompositions. Put differently, we are interested in what happens when we glue two surfaces along a circle. \cite{SP,LEE} show that Wrapped Fukaya categories behave well under these decompositions.

We take here the following model: let $(\Sigma, \ff)$ be a (possibly disconnected) stopped surface and $C_-, C_+\subset \partial \Sigma\setminus \ff$ two disjoint circles. Gluing these two circles can be done in the following three steps: Add two points $q_-\in C_-,q_+\in C_+$ to the set of stops $\ff$. Connect $q_-$ and $q_+$ by attaching a one-handle as described in section \ref{sec:lineargluing} and call the resulting circle $C$. Glue in a disk by identifying its boundary with $C$.

To prove the circular gluing result we need the following lemma:
\begin{lem}\label{delete}
Suppose a stopped surface $(\Sigma', \ff'=\ff)$ arises by gluing a disk $D$ to a boundary circle $C\subset \partial \Sigma\setminus \ff$ of a stopped surface $(\Sigma, \ff)$. Then  $$\Omega_\ff(\Sigma)\cong \Omega_\ff(\Sigma')\oplus \IZ.$$
\end{lem}
The effect of this operation on the skeleton of the surface is collapsing a circle to a point.
\begin{proof}
By adding a point to $\ff'$ we can assume without loss of generality that $\ff'\neq \emptyset$. Let $A'$ be a minimal full arc system for $(\Sigma', \ff')$. Thinking of $\Sigma$ as a subset of $\Sigma'$ we can assume that $A'$ is disjoint from the disk $D$ as $D$ is contractible. Cutting $\Sigma$ along the arcs in $A'$ leaves topological discs containing exactly one point of the stop, except of an annulus with one point of the stop (the component that contains $D$). Then adding one more arc which cuts this annulus into a disk to $A'$ is a minimal full arc system for $(\Sigma, \ff)$.
\end{proof}

The discussion at the beginning of this subsection concludes the following statement for circular gluing by first applying Corollary \ref{cor:stopremoval}, then Proposition \ref{handle} and finally Lemma \ref{delete}.

\begin{prop}
Suppose a connected stopped surface $(\Sigma', \ff'=\ff)$ arises by identifying two disjoint boundary circles $C_-, C_+\subset \partial \Sigma\setminus \ff$  of a stopped surface $(\Sigma, \ff)$.

If either $\Sigma$ is connected or $\Sigma$ has two connected components and each of them has a non-empty stop then $$\Omega_{\ff}(\Sigma)\cong \Omega_{\ff'}(\Sigma').$$
If $\Sigma$ has two connected components with at most one component having a non-empty stop, then $$\Omega_{\ff}(\Sigma)\cong \Omega_{\ff'}(\Sigma')\oplus \IZ.$$
\end{prop}

\subsection{More on Cobordism groups in higher dimensions}\label{sec:highercob}
In this section we discuss two more techniques to get relations in cobordism groups in higher dimensions and give some examples. 

\subsubsection{Relations induced by surgeries and Hamiltonian isotopies}
Let $M=S^n$ and let $F_q=T_q^*S^n$ be the fiber over a point $q\in S^n$. The Dehn twist $\tau(F_q)$ of $F_q$ is given by the time-$1$ flow of the fiber $F_{-q}$ with Hamiltonian flow coming from the Hamiltonian $H(q,p)=\rho(|p|)$ where $\rho:[0,\infty)\to\IR$ is a smooth increasing function with $\rho(0)=\rho'(0)=0$ and $\rho(r)=\frac r2$ for $r>>0$. Note that the Hamiltonian $H$ is linear at infinity. 

An alternative description of the Dehn twist $\tau(F_q)$ is given by doing Polterovich surgery of the fiber $F_q$ with the zero section $S^n\subset T^*S^n$.

These two different descriptions produce the following identities in the cobordism group:
$$F_q=F_{-q}=\tau(F_q)=F_q+S^n.$$
We conclude:
\begin{prop}
The zero section $[S^n]$ is zero in $\Omega(T^*S^n)$. 
\end{prop}

Let $S\subset M$ be a closed submanifold of a closed submanifold $M$. Recall that the conormal bundle $$N^*S=\{(q,\xi)\in T^*_qM \text{ where } \xi\vert_{T_qS}=0\}\subset T^*M$$ describes a strongly exact (see page \pageref{strexact}) Lagrangian submanifold of $T^*M$. 

\begin{prop}
Let $(S_t)_{t\in [0,1]}$ be a smooth family of spheres in a manifold $M$ such that $S_t$ is embedded for $t\in [0,1)$ and $S_1$ a point in $S_t\subset M$ for all $t$. Then the class of the conormal $[N^*S_0]$ is zero in $\Omega(T^*M)$.
\end{prop}
\begin{proof}
Denote $S_1=\{q\}$ and $F=T_q^*M$ the fiber over $q$. The Lagrangian submanifolds $F$ and $N^*S_t$ intersect cleanly. Doing Polterovich surgeries of $F$ and $N^*S_t$ that depend smoothly on $t$ produces a smooth family of exact conical Lagrangian submanifolds connecting the surgery $F\# N^*S$ to $F$. Consequently, $F\# N^*S$ and $F$ are isotopic by a linear Hamiltonian isotopy. We get the relation $F+ N^*S=F$ in the cobordism group which concludes the proof.
\end{proof}

Note that the argument via surgery was necessary: $(N^*S_t)_{t\in[0,1]}$ is not a smooth family of conical exact Lagrangian submanifolds, otherwise $[N^*S_0]=[F]$ would hold. But $[F]\neq 0$ as it is non-zero in the Grothendieck group of the wrapped Fukaya category. In other words, the ball coming from the contraction of $S_t$ to the point $S_0$ is not inducing a Lagrangian cobordism. Compare to the Example \ref{ex:noncob}.

\subsubsection{K\"unneth formula for cobordism groups}
Cobordism groups satisfy a K\"unneth formula. If $(W_1,\ff),(W_2,\gggg)$ are two stopped Liouville manifolds then there is a group homomorphism
\begin{equation}\label{Kun}
    \Omega_\ff(W_1)\times \Omega_\gggg(W_2)\to \Omega_{\mathfrak h}(W_1\times W_2)
\end{equation}
given by $(L,K) \mapsto L \tilde\times K$, where the product stop $\mathfrak h$ is given as in Example $\ref{ex:stabs}$. Indeed, for any Lagrangian cobordism $V\subset \IC\times W_1$ with ends $L_j$ and a Lagrangian submanifold $K\subset W_2$ the conicalized product $V\tilde \times K\subset \IC\times W_1\times W_2$ is a Lagrangian cobordism with ends $L_j\tilde\times K$.

Applying (\ref{Kun}) to the zero section in the cylinder $T^*S^1$ (Example \ref{sec:cylinder}) we conclude:
\begin{prop}
The zero section $[\mathbb T^n]$ is zero in $\Omega(T^*\mathbb T^n)$. 
\end{prop}

\appendix
\section{The Grothendieck group and cone decompositions}\label{sec:appendix}
We collect here some remarks about triangulated categories, the Grothendieck group and iterated cone decompositions. More can be found in \cite{WEI},\cite{SEI}, \cite{BC2} (Section 2.6) and \cite{HEN} (Appendix A).

The \emph{Grothendieck group} $K_0(\mathcal A)$ of a triangulated $\mathcal A$ is the group freely generated by objects in $\mathcal A$ modulo the relation $A_1-A_2+ A=0$ whenever there is an exact triangle 
$$A_1\to A_2\to A\to A[1]$$ in $\mathcal A$. Every morphism $A_1\overset{a}{\to} A_2$ can be completed to an exact triangle $A_1\to A_2\to A\to A[1]$ where $A$ is unique up to isomorphism. In this case, we call $A$ the \emph{cone} over $a:A_1\to A_2$ and write $A\cong [A_1\overset{a}{\to} A_2]$.

An ($n$-step) \emph{iterated cone decomposition} of an object $A$ in $\mathcal A$ is a sequence of exact triangles
\begin{align}\label{eq:exact}
A_{j}[-1]\to X_{j-1}\to X_j\to A_{j}    
\end{align}
for $j=1, \ldots,n$ with objects $X_1=0, X_2,\ldots, X_{n-1}, X_n=A$ in $\mathcal A$.

See \cite{BC2} for a more detailed discussion, e.g. when two cone decompositions are isomorphic. Let us record here some simple facts about cone decompositions related to grading which directly follow from the axioms of triangulated categories and the given definitions above:
\begin{enumerate}
    \item An iterated cone decomposition (\ref{eq:exact}) induces a relation $A=A_1+\ldots +A_n$ in the Grothendieck group $K_0(\mathcal A)$. 
    \item Substituting isomorphic objects we can write the exact triangles (\ref{eq:exact}) as cones $X_j\cong [A_{j}[-1]\to X_{j-1}]$ and the iterated cone can be written as
\begin{align}\label{eq:ex3}
A\cong\big[A_n[-1]\to\big[A_{n-1}[-1]\to\cdots \to \big[A_2[-1]\to A_1\big]\cdots\big]\big].
\end{align}
    \item Rotating exact triangles and using the isomorphic decomposition 
    \begin{align}\label{eq:rot}
        [[X\to Y]\to Z]\cong [X[1]\to[Y\to Z]]
    \end{align} the iterated cone decomposition (\ref{eq:ex3}) can be rotated as
    \begin{align}
        A_n\cong\big[A_{n-1}\to\big[A_{n-2}\to\cdots \to \big[A_1\to A\big]\cdots\big]\big].
    \end{align}
    The other rotations of the iterated cone decompositions yield grading shifts with less insight. 
    \item If one of the objects in $A_j$ admits itself an iterated cone decomposition 
    $$A_j\cong\big[B_{m}[-1]\to\big[B_{m-1}[-1]\to\cdots \to \big[B_2[-1]\to B_1\big]\cdots\big]\big]$$
    then we can combine it with (\ref{eq:ex3}) and get an iterated cone decomposition
    \begin{align*}
    A\cong\big[A_n[-1]\to\cdots \to [A_{j+1}[-1]\to \big[B_{m}[-1]&\to\cdots \to \big [B_1[-1]\to [A_{j-1}[-1] \cdots \to A_1\big]\cdots\big]\big]\big]\big]    
    \end{align*}
    using (\ref{eq:rot}) several times.
\end{enumerate}

\bibliography{bibi}

\providecommand{\etalchar}[1]{$^{#1}$}
\providecommand{\bysame}{\leavevmode\hbox to3em{\hrulefill}\thinspace}
\providecommand{\noopsort}[1]{}
\providecommand{\mr}[1]{\href{http://www.ams.org/mathscinet-getitem?mr=#1}{MR~#1}}
\providecommand{\zbl}[1]{\href{http://www.zentralblatt-math.org/zmath/en/search/?q=an:#1}{Zbl~#1}}
\providecommand{\jfm}[1]{\href{http://www.emis.de/cgi-bin/JFM-item?#1}{JFM~#1}}
\providecommand{\arxiv}[1]{\href{http://www.arxiv.org/abs/#1}{arXiv~#1}}
\providecommand{\MR}{\relax\ifhmode\unskip\space\fi MR }
\providecommand{\MRhref}[2]{%
  \href{http://www.ams.org/mathscinet-getitem?mr=#1}{#2}
}
\providecommand{\href}[2]{#2}
\begin{thebibliography}{CGGR19}

\bibitem[AS06]{ABSC}
\bgroup\scshape{}A.~Abbondandolo\egroup{} and
  \bgroup\scshape{}M.~Schwarz\egroup{}, On the {F}loer homology of cotangent
  bundles,  \emph{Communications on Pure and Applied Mathematics: A Journal
  Issued by the Courant Institute of Mathematical Sciences} \textbf{59} no.~2
  (2006), 254--316.

\bibitem[Abo10]{ABO2}
\bgroup\scshape{}M.~Abouzaid\egroup{}, A geometric criterion for generating the
  {F}ukaya category,  \emph{Publications math{\'e}matiques de l'IH{\'E}S}
  \textbf{112} (2010), 191--240.

\bibitem[Abo11]{ABO3}
\bgroup\scshape{}M.~Abouzaid\egroup{}, A cotangent fibre generates the {F}ukaya
  category,  \emph{Advances in Mathematics} \textbf{228} no.~2 (2011), 894 --
  939.

\bibitem[AAE{\etalchar{+}}13]{ABETC}
\bgroup\scshape{}M.~Abouzaid\egroup{}, \bgroup\scshape{}D.~Auroux\egroup{},
  \bgroup\scshape{}A.~Efimov\egroup{}, \bgroup\scshape{}L.~Katzarkov\egroup{},
  and \bgroup\scshape{}D.~Orlov\egroup{}, Homological mirror symmetry for
  punctured spheres,  \emph{Journal of the American Mathematical Society}
  \textbf{26} no.~4 (2013), 1051--1083.

\bibitem[AS10]{AS}
\bgroup\scshape{}M.~Abouzaid\egroup{} and \bgroup\scshape{}P.~Seidel\egroup{},
  An open string analogue of {V}iterbo functoriality,  \emph{Geometry \&
  Topology} \textbf{14} no.~2 (2010), 627--718.

\bibitem[Arn80a]{ARN1}
\bgroup\scshape{}V.~I. Arnol'd\egroup{}, Lagrange and {L}egendre cobordisms.
  {I},  \emph{Akademiya Nauk SSSR. Funktsional'ny Analiz i ego Prilozheniya}
  \textbf{14} no.~3 (1980), 1--13, 96.

\bibitem[Arn80b]{ARN2}
\bgroup\scshape{}V.~I. Arnol'd\egroup{}, Lagrange and {L}egendre cobordisms.
  {II},  \emph{Akademiya Nauk SSSR. Funktsional' ny Analiz i ego Prilozheniya}
  \textbf{14} no.~4 (1980), 8--17, 95.

\bibitem[Aud85]{AUD}
\bgroup\scshape{}M.~Audin\egroup{}, Quelques calculs en cobordisme lagrangien,
  in \emph{Annales de l'institut Fourier}, \textbf{35}, 1985, pp.~159--194.

\bibitem[Aur10a]{AUR2}
\bgroup\scshape{}D.~Auroux\egroup{}, {F}ukaya categories and bordered
  {H}eegaard-{F}loer homology,  in \emph{Proceedings of the International
  Congress of Mathematicians 2010 (ICM 2010) (In 4 Volumes) Vol. I: Plenary
  Lectures and Ceremonies Vols. II--IV: Invited Lectures}, World Scientific,
  2010, pp.~917--941.

\bibitem[Aur10b]{AUR}
\bgroup\scshape{}D.~Auroux\egroup{}, Fukaya categories of symmetric products
  and bordered {H}eegaard-{F}loer homology, 2010. \arxiv{1001.4323}.

\bibitem[BC13]{BC}
\bgroup\scshape{}P.~Biran\egroup{} and \bgroup\scshape{}O.~Cornea\egroup{},
  {L}agrangian cobordism. {I},  \emph{Journal of the American Mathematical
  Society} \textbf{26} no.~2 (2013), 295--340.

\bibitem[BC14]{BC2}
\bgroup\scshape{}P.~Biran\egroup{} and \bgroup\scshape{}O.~Cornea\egroup{},
  {L}agrangian cobordism and {F}ukaya categories,  \emph{Geometric and
  functional analysis} \textbf{24} no.~6 (2014), 1731--1830.

\bibitem[BC17]{BC3}
\bgroup\scshape{}P.~Biran\egroup{} and \bgroup\scshape{}O.~Cornea\egroup{},
  Cone-decompositions of {L}agrangian cobordisms in {L}efschetz fibrations,
  \emph{Selecta Mathematica} \textbf{23} no.~4 (2017), 2635--2704.

\bibitem[Boc16]{BOC}
\bgroup\scshape{}R.~Bocklandt\egroup{}, Noncommutative mirror symmetry for
  punctured surfaces,  \emph{Transactions of the American Mathematical Society}
  \textbf{368} no.~1 (2016), 429--469.

\bibitem[CGGR19]{CHA2}
\bgroup\scshape{}B.~Chantraine\egroup{}, \bgroup\scshape{}P.~Ghiggini\egroup{},
  \bgroup\scshape{}R.~Golovko\egroup{}, and \bgroup\scshape{}G.~D.
  Rizell\egroup{}, Geometric generation of the wrapped fukaya category of
  weinstein manifolds and sectors, 2019. \arxiv{1712.09126}.

\bibitem[Che97]{CHE}
\bgroup\scshape{}Y.~V. Chekanov\egroup{}, Lagrangian embeddings and
  {L}agrangian cobordism,  \emph{Translations of the American Mathematical
  Society-Series 2} \textbf{180} (1997), 13--24.

\bibitem[CE12]{CE}
\bgroup\scshape{}K.~Cieliebak\egroup{} and
  \bgroup\scshape{}Y.~Eliashberg\egroup{}, \emph{From Stein to {W}einstein and
  back: symplectic geometry of affine complex manifolds}, \textbf{59}, American
  Mathematical Soc., 2012.

\bibitem[CR20]{RIZ}
\bgroup\scshape{}L.~C{\^o}t{\'e}\egroup{} and \bgroup\scshape{}G.~D.
  Rizell\egroup{}, Symplectic rigidity of fibers in cotangent bundles of
  {R}iemann surfaces, 2020. \arxiv{2004.04233}.

\bibitem[DRGI16]{RIZ5}
\bgroup\scshape{}G.~Dimitroglou~Rizell\egroup{},
  \bgroup\scshape{}E.~Goodman\egroup{}, and \bgroup\scshape{}A.~Ivrii\egroup{},
  Lagrangian isotopy of tori in {$S^2\times S^2$} and {$\mathbb{C}P^2$},
  \emph{Geometric and Functional Analysis} \textbf{26} no.~5 (2016),
  1297--1358.

\bibitem[Dur17]{DUR}
\bgroup\scshape{}J.~Dureti{\'c}\egroup{}, Piunikhin-{S}alamon-{S}chwarz
  isomorphisms and spectral invariants for conormal bundle,  \emph{Publications
  de l'Institut Mathematique} \textbf{102} no.~116 (2017), 17--47.

\bibitem[Dyc17]{DYC2}
\bgroup\scshape{}T.~Dyckerhoff\egroup{}, $\mathbb{A}^{1}$-homotopy invariants
  of topological {F}ukaya categories of surfaces,  \emph{Compositio
  Mathematica} \textbf{153} no.~8 (2017), 1673--1705.

\bibitem[Dyc21]{DYC}
\bgroup\scshape{}T.~Dyckerhoff\egroup{}, A categorified {D}old-{K}an
  correspondence,  \emph{Selecta Mathematica. New Series} \textbf{27} no.~2
  (2021), 14.

\bibitem[DJL21]{DJL}
\bgroup\scshape{}T.~Dyckerhoff\egroup{}, \bgroup\scshape{}G.~Jasso\egroup{},
  and \bgroup\scshape{}Y.~Lekili\egroup{}, The symplectic geometry of higher
  {A}uslander algebras: symmetric products of disks,  \emph{Forum of
  Mathematics. Sigma} \textbf{9} (2021), Paper No. e10, 49.

\bibitem[DK18]{DK}
\bgroup\scshape{}T.~Dyckerhoff\egroup{} and
  \bgroup\scshape{}M.~Kapranov\egroup{}, Triangulated surfaces in triangulated
  categories,  \emph{Journal of the European Mathematical Society (JEMS)}
  \textbf{20} no.~6 (2018), 1473--1524.

\bibitem[EP96]{EP}
\bgroup\scshape{}Y.~Eliashberg\egroup{} and
  \bgroup\scshape{}L.~Polterovich\egroup{}, Local {L}agrangian {$2$}-knots are
  trivial,  \emph{Annals of Mathematics. Second Series} \textbf{144} no.~1
  (1996), 61--76.

\bibitem[Eli84]{ELI}
\bgroup\scshape{}Y.~Eliashberg\egroup{}, Cobordisme des solutions de relations
  diff{\'e}rentielles,  in \emph{South Rhone seminar on geometry, Travaux en
  Cours, 1984}, 1984.

\bibitem[{F}uk93]{FU}
\bgroup\scshape{}K.~{F}ukaya\egroup{}, Morse homotopy, {$A_\infty$}-category
  and {F}loer homologies,  in \emph{Proceeding of Garc Workshop on Geometry and
  Topology}, Seoul National Univ, 1993.

\bibitem[FSS08]{FSS}
\bgroup\scshape{}K.~{F}ukaya\egroup{}, \bgroup\scshape{}P.~Seidel\egroup{}, and
  \bgroup\scshape{}I.~Smith\egroup{}, The symplectic geometry of cotangent
  bundles from a categorical viewpoint,  in \emph{Homological mirror symmetry},
  Springer, 2008, pp.~1--26.

\bibitem[GPS19]{GPS2}
\bgroup\scshape{}S.~Ganatra\egroup{}, \bgroup\scshape{}J.~Pardon\egroup{}, and
  \bgroup\scshape{}V.~Shende\egroup{}, Sectorial descent for wrapped {F}ukaya
  categories, 2019. \arxiv{1809.03427}.

\bibitem[GPS20]{GPS1}
\bgroup\scshape{}S.~Ganatra\egroup{}, \bgroup\scshape{}J.~Pardon\egroup{}, and
  \bgroup\scshape{}V.~Shende\egroup{}, Covariantly functorial wrapped {F}loer
  theory on {L}iouville sectors,  \emph{Publications Math\'{e}matiques.
  Institut de Hautes \'{E}tudes Scientifiques} \textbf{131} (2020), 73--200.

\bibitem[Gei08]{GEI}
\bgroup\scshape{}H.~Geiges\egroup{}, \emph{An introduction to contact
  topology}, \textbf{109}, Cambridge University Press, 2008.

\bibitem[Gro85]{GRO}
\bgroup\scshape{}M.~Gromov\egroup{}, Pseudo holomorphic curves in symplectic
  manifolds,  \emph{Inventiones Mathematicae} \textbf{82} no.~2 (1985),
  307--347.

\bibitem[HKK17]{HKK}
\bgroup\scshape{}F.~Haiden\egroup{}, \bgroup\scshape{}L.~Katzarkov\egroup{},
  and \bgroup\scshape{}M.~Kontsevich\egroup{}, Flat surfaces and stability
  structures,  \emph{Publications math{\'e}matiques de l'IH{\'E}S} \textbf{126}
  no.~1 (2017), 247--318.

\bibitem[Har85]{HAR1}
\bgroup\scshape{}J.~L. Harer\egroup{}, Stability of the homology of the mapping
  class groups of orientable surfaces,  \emph{Annals of Mathematics}
  \textbf{121} no.~2 (1985), 215--249.

\bibitem[Har86]{HAR2}
\bgroup\scshape{}J.~L. Harer\egroup{}, The virtual cohomological dimension of
  the mapping class group of an orientable surface,  \emph{Inventiones
  mathematicae} \textbf{84} no.~1 (1986), 157--176.

\bibitem[Hat91]{HAT}
\bgroup\scshape{}A.~Hatcher\egroup{}, On triangulations of surfaces,
  \emph{Topology and its Applications} \textbf{40} no.~2 (1991), 189--194.

\bibitem[Hau15]{HAUG}
\bgroup\scshape{}L.~Haug\egroup{}, The {L}agrangian cobordism group of {$T^2$},
   \emph{Selecta Mathematica. New Series} \textbf{21} no.~3 (2015), 1021--1069.

\bibitem[Hen20]{HEN}
\bgroup\scshape{}F.~Hensel\egroup{}, Stability conditions and {L}agrangian
  cobordisms,  \emph{Journal of Symplectic Geometry} \textbf{18} no.~2 (2020),
  463--536.

\bibitem[Hic20]{HIC}
\bgroup\scshape{}J.~Hicks\egroup{}, Tropical lagrangian hypersurfaces are
  unobstructed,  \emph{Journal of Topology} \textbf{13} no.~4 (2020),
  1409--1454.

\bibitem[Hin12]{HIND}
\bgroup\scshape{}R.~Hind\egroup{}, {L}agrangian unknottedness in {S}tein
  surfaces,  \emph{Asian Journal of Mathematics} \textbf{16} no.~1 (2012),
  1--36.

\bibitem[LS91]{LS}
\bgroup\scshape{}F.~Lalonde\egroup{} and \bgroup\scshape{}J.-C.
  Sikorav\egroup{}, Sous-vari{\'e}t{\'e}s {L}agrangiennes et {L}agrangiennes
  exactes des fibr{\'e}s cotangents,  \emph{Commentarii mathematici Helvetici}
  \textbf{66} no.~1 (1991), 18--33.

\bibitem[Lee15]{LEE}
\bgroup\scshape{}H.~Lee\egroup{}, \emph{Homological mirror symmetry for open
  Riemann surfaces from pair-of-pants decompositions}, Ph.D. thesis, University
  of California, Berkeley, 2015.

\bibitem[LP18]{LP}
\bgroup\scshape{}Y.~Lekili\egroup{} and
  \bgroup\scshape{}A.~Polishchuk\egroup{}, Auslander orders over nodal stacky
  curves and partially wrapped {F}ukaya categories,  \emph{Journal of Topology}
  \textbf{11} no.~3 (2018), 615--644.

\bibitem[LP20]{LP2}
\bgroup\scshape{}Y.~Lekili\egroup{} and
  \bgroup\scshape{}A.~Polishchuk\egroup{}, Derived equivalences of gentle
  algebras via {F}ukaya categories,  \emph{Mathematische Annalen} \textbf{376}
  no.~1-2 (2020), 187--225.

\bibitem[MW18]{MW}
\bgroup\scshape{}C.~Y. Mak\egroup{} and \bgroup\scshape{}W.~Wu\egroup{}, Dehn
  twist exact sequences through {L}agrangian cobordism,  \emph{Compositio
  Mathematica} \textbf{154} no.~12 (2018), 2485--2533.

\bibitem[MS17]{MS}
\bgroup\scshape{}D.~McDuff\egroup{} and \bgroup\scshape{}D.~Salamon\egroup{},
  \emph{Introduction to symplectic topology}, Oxford University Press, 2017.

\bibitem[Nad15]{NAD}
\bgroup\scshape{}D.~Nadler\egroup{}, Cyclic symmetries of {$A_n$}-quiver
  representations,  \emph{Advances in Mathematics} \textbf{269} (2015),
  346--363.

\bibitem[NT20]{TAN2}
\bgroup\scshape{}D.~Nadler\egroup{} and \bgroup\scshape{}H.~L. Tanaka\egroup{},
  A stable {$\infty$}-category of {L}agrangian cobordisms,  \emph{Advances in
  Mathematics} \textbf{366} (2020), 107026, 97.

\bibitem[Oh15]{OH}
\bgroup\scshape{}Y.-G. Oh\egroup{}, \emph{Symplectic topology and {F}loer
  homology}, \textbf{2}, Cambridge University Press, 2015.

\bibitem[PS19]{SP}
\bgroup\scshape{}J.~Pascaleff\egroup{} and
  \bgroup\scshape{}N.~Sibilla\egroup{}, Topological {F}ukaya category and
  mirror symmetry for punctured surfaces,  \emph{Compositio Mathematica}
  \textbf{155} no.~3 (2019), 599--644.

\bibitem[Per19]{PER}
\bgroup\scshape{}A.~Perrier\egroup{}, Lagrangian cobordism groups of higher
  genus surfaces, 2019. \arxiv{1901.06002}.

\bibitem[Pog17]{POG}
\bgroup\scshape{}T.~Poguntke\egroup{}, Higher segal structures in algebraic
  $k$-theory, 2017. \arxiv{1709.06510}.

\bibitem[Pol91]{POL}
\bgroup\scshape{}L.~Polterovich\egroup{}, The surgery of {L}agrange
  submanifolds,  \emph{Geometric \& Functional Analysis GAFA} \textbf{1} no.~2
  (1991), 198--210.

\bibitem[Rit09]{RIT}
\bgroup\scshape{}A.~F. Ritter\egroup{}, Novikov-symplectic cohomology and exact
  lagrangian embeddings,  \emph{Geometry \& Topology} \textbf{13} no.~2 (2009),
  943--978.

\bibitem[Sei08]{SEI}
\bgroup\scshape{}P.~Seidel\egroup{}, \emph{{F}ukaya categories and
  {P}icard-{L}efschetz theory}, \textbf{10}, European Mathematical Society,
  2008.

\bibitem[SS21]{SS}
\bgroup\scshape{}N.~Sheridan\egroup{} and \bgroup\scshape{}I.~Smith\egroup{},
  Lagrangian cobordism and tropical curves,  \emph{Journal f{\"u}r die reine
  und angewandte Mathematik (Crelles Journal)} \textbf{2021} no.~774 (2021),
  219--265.

\bibitem[Str84]{STR}
\bgroup\scshape{}K.~Strebel\egroup{}, Quadratic differentials,  in
  \emph{Quadratic Differentials}, Springer, 1984, pp.~16--26.

\bibitem[Syl19]{SYL}
\bgroup\scshape{}Z.~Sylvan\egroup{}, On partially wrapped {F}ukaya categories,
  \emph{Journal of Topology} \textbf{12} no.~2 (2019), 372--441.

\bibitem[Tan16a]{TAN5}
\bgroup\scshape{}H.~L. Tanaka\egroup{}, The {F}ukaya category pairs with
  {L}agrangian cobordisms, 2016. \arxiv{1607.04976}.

\bibitem[Tan16b]{TAN6}
\bgroup\scshape{}H.~L. Tanaka\egroup{}, The {F}ukaya category pairs with
  {L}agrangian cobordisms exactly, 2016. \arxiv{1609.08400}.

\bibitem[Tan18]{TAN}
\bgroup\scshape{}H.~L. Tanaka\egroup{}, Surgery induces exact sequences in
  {L}agrangian cobordisms, 2018. \arxiv{1805.07424}.

\bibitem[Tan19]{TAN3}
\bgroup\scshape{}H.~L. Tanaka\egroup{}, Cyclic structures and broken cycles,
  2019. \arxiv{1907.03301}.

\bibitem[Tan20]{TAN7}
\bgroup\scshape{}H.~L. Tanaka\egroup{}, Generation for lagrangian cobordisms in
  weinstein manifolds, 2020. \arxiv{1810.10605}.

\bibitem[Tan21]{TAN4}
\bgroup\scshape{}H.~L. Tanaka\egroup{}, In simply connected cotangent bundles,
  exact {L}agrangian cobordisms are h-cobordisms,  \emph{Advances in Geometry}
  \textbf{21} no.~1 (2021), 1--4.

\bibitem[Vit97]{VIT2}
\bgroup\scshape{}C.~Viterbo\egroup{}, Exact {L}agrange submanifolds, periodic
  orbits and the cohomology of free loop spaces,  \emph{Journal of Differential
  Geometry} \textbf{47} no.~3 (1997), 420--468.

\bibitem[Wei95]{WEI}
\bgroup\scshape{}C.~A. Weibel\egroup{}, \emph{An introduction to homological
  algebra}, no.~38, Cambridge university press, 1995.

\end{thebibliography}
\bibliographystyle{aomalpha}
\end{document}